\documentclass[12pt,a4paper]{article}
\usepackage{amsmath,amssymb,amscd,amsthm,latexsym} 
\renewcommand\today{June 3, 2002}

\theoremstyle{plain}  % default
\newtheorem{theorem}{Theorem}[section]
\newtheorem*{theorem*}{Theorem}

\newtheorem{corollary}[theorem]{Corollary}
\newtheorem{lemma}[theorem]{Lemma}
\newtheorem{proposition}[theorem]{Proposition}

\newtheorem{tech-lemma}[theorem]{Technical Lemma}
\newtheorem{definition}[theorem]{Definition}
\theoremstyle{definition}
\theoremstyle{remark}

\newtheorem{remark}[theorem]{Remark}

\newtheorem*{claim*}{Claim}
\numberwithin{equation}{section}
\renewcommand{\leq}{\leqslant}
\renewcommand{\le}{\leqslant}
\renewcommand{\geq}{\geqslant}
\renewcommand{\ge}{\geqslant}
\renewcommand{\setminus}{\smallsetminus}
\newcommand{\lto}{\longrightarrow}
\newcommand{\R}{\mathbb{R}}
\newcommand{\Z}{\mathbb{Z}}
\newcommand{\C}{\mathbb{C}}
\newcommand{\HH}{\mathbb{H}}
\newcommand{\dbar}{\bar{\partial}}
\newcommand{\into}{\hookrightarrow}
\newcommand{\abs}[1]{\lvert#1\rvert}
\newcommand{\norm}[1]{\lVert#1\rVert}

\newcommand{\vol}{\mathrm{vol}}

\newcommand{\balpha}{{\boldsymbol{\alpha}}}
\newcommand{\btau}{{\boldsymbol{\tau}}}
\newcommand{\PU}{\mathrm{PU}}
\newcommand{\PGL}{\mathrm{PGL}} 
\newcommand{\SU}{\mathrm{SU}}
\newcommand{\U}{\mathrm{U}}
\newcommand{\GL}{\mathrm{GL}}
\newcommand{\GCD}{\mathrm{GCD}}

\DeclareMathOperator{\ad}{ad}
\DeclareMathOperator{\Ad}{Ad}
\DeclareMathOperator{\Div}{Div}

\DeclareMathOperator{\tr}{tr}
\DeclareMathOperator{\rk}{rk}
\DeclareMathOperator{\im}{im}
\DeclareMathOperator{\coker}{coker}
\DeclareMathOperator{\Hom}{Hom}
\DeclareMathOperator{\End}{End}
\DeclareMathOperator{\Ext}{Ext}
\DeclareMathOperator{\Id}{Id}
\DeclareMathOperator{\Quot}{Quot}

\begin{document}
\begin{titlepage}
%\mbox
%\bigskip
\noindent 
{\Large\textbf{Surface group representations, Higgs bundles, and 
 holomorphic triples}} 
\bigskip

\noindent 
\textbf{Steven B. Bradlow}\footnotemark[1]$^{,}$\footnotemark[2] \\
%\textit{\nocorr
Department of Mathematics, \\
University of Illinois, \\
Urbana,
IL 61801, 
USA \\
E-mail: \texttt{bradlow@math.uiuc.edu}
%}
\medskip

\noindent 
\textbf{Oscar Garc{\'\i}a--Prada}\footnotemark[1]$^{,}$\footnotemark[3]$^{,}$\footnotemark[5]$^{,}$\footnotemark[6]  \\
%\textit{\nocorr
Departamento de Matem{\'a}ticas, \\
Universidad Aut{\'o}noma de Madrid, \\
28049 Madrid,
Spain \\
E-mail: \texttt{oscar.garcia-prada@uam.es}
%}

\medskip
\noindent 
\textbf{Peter B. Gothen}\footnotemark[1]$^{,}$\footnotemark[4]$^{,}$\footnotemark[5]   \\
%\textit{\nocorr
Departamento de Matem{\'a}tica Pura, \\
Faculdade de Ci{\^e}ncias,
Universidade do Porto, \\
Rua do Campo Alegre 687, 4169-007 Porto,
Portugal \\
E-mail: \texttt{pbgothen@fc.up.pt}
%}
%\date{April 12, 2002}
\footnotetext[1]{
Members of VBAC (Vector Bundles on
Algebraic Curves), which is partially supported by EAGER (EC FP5
Contract no.\ HPRN-CT-2000-00099) and by EDGE (EC FP5 Contract no.\ 
HPRN-CT-2000-00101).}
\footnotetext[2]{Partially supported by 
the National Science Foundation under grant DMS-0072073 }
\footnotetext[3]{Partially supported by 
the Ministerio de Ciencia y Tecnolog\'{\i}a (Spain) under grant BFM2000-0024}
\footnotetext[4]{
Partially supported by the
Funda{\c c}{\~a}o para a Ci{\^e}ncia e a Tecnologia (Portugal) through
the Centro de Matem{\'a}tica da Universidade do Porto and through
grant no.\ SFRH/BPD/1606/2000.}
\footnotetext[5]{
Partially supported by the Portugal/Spain bilateral Programme Acciones
Integradas, 
grant nos.\ HP2000-0015 and AI-01/24}
\footnotetext[6]{
Partially supported by a  British EPSRC grant (October-December 2001)}

\bigskip
\noindent
\textbf{\today}
\vfill

\newpage

\pagestyle{empty}

\paragraph{Abstract.}
Using the $L^2$ norm of the Higgs field as a Morse function, we study
the moduli spaces of $\U(p,q)$-Higgs bundles over a Riemann surface.
We require that the genus of the surface be at least two, but place
no constraints on $(p,q)$.  A key step is the identification of the
function's local minima as moduli spaces of holomorphic triples.  We
prove that these moduli spaces of triples are irreducible and
non-empty.

Because of the relation between flat bundles and fundamental group
representations, we can interpret our conclusions as results about
the number of connected components in the moduli space of semisimple
$\PU(p,q)$-representations.  The topological invariants of the flat
bundles bundle are used to label components. These invariants are
bounded by a Milnor--Wood type inequality. For each allowed value of
the invariants satisfying a certain coprimality condition, we prove
that the corresponding component is non-empty and connected. If the
coprimality condition does not hold, our results apply to the
irreducible representations.
\vfill
\newpage
\tableofcontents
\end{titlepage}
\pagestyle{plain}
%%%%%%%%%%%%%%%%%%%%%
\section{Introduction} \label{sec:introduction}
%%%%%%%%%%%%%%%%%%%%%%

 The core of this paper is a Morse theoretic study of the the moduli
 space of $\U(p,q)$-Higgs bundles over a Riemann surface $X$ of genus $g\ge 2$.
 Our interest in this space comes from two sources. The first is its
 relevance to questions concerning the representation variety for
 representations of $\pi_1 X$ in the real Lie group $\PU(p,q)$. The
 second has to do with the intrinsic geometry revealed by the Morse
 function and the methods we are able to use to carry out our
 analysis. Our main goals are to fully understand the minimal 
 submanifolds of the Morse function and, thereby, to count the number
 of connected components in the representation varieties.
 
A  Higgs bundle consists of a  holomorphic bundle together with a Higgs
 field, i.e.\ a section of a certain associated vector bundle. A
 $\U(p,q)$-Higgs bundle is a special case of the $G$-Higgs bundles
 defined by Hitchin in \cite{hitchin:1987}, where $G$ is a real form
 of a complex reductive Lie group. Such objects provide a natural generalization
 of holomorphic vector bundles, which correspond to the case 
 $G=\U(n)$ and zero Higgs field. In particular, they permit an 
 extension to other groups of the Narasimhan and Seshadri theorem 
 (\cite{narasimhan-seshadri:1965})
  on the relation between unitary representations of $\pi_1 X$ and stable 
 vector bundles.
 
  By embedding $\U(p,q)$ in $\GL(p+q)$ we can give a concrete description 
  of a $\U(p,q)$-Higgs bundle as a pair
 $$(V\oplus W,\Phi=\left(
   \begin{smallmatrix}
     0 & \beta \\
     \gamma & 0
   \end{smallmatrix}
   \right))\
 $$
 \noindent where $V$ and $W$ are holomorphic vector bundles
 of rank $p$ and $q$ respectively, $\beta$ is a section in
 $H^0(\Hom(W,V) \otimes K)$, and $\gamma \in H^0(\Hom(V,W) \otimes K)$,
 so that $\Phi\in H^0(\mathrm{End}(V\oplus W) \otimes K)$. Foremost
 among the key features of such objects is (by the work of Hitchin,
 \cite{hitchin:1987,hitchin:1992} Donaldson \cite{donaldson:1987},
  Corlette \cite{corlette:1988} and Simpson
  \cite{simpson:1988,simpson:1992,simpson:1994a,simpson:1994b}) 
 the existence of moduli spaces of
 polystable objects which can be identified with moduli spaces of
 solutions to natural gauge theoretic equations. Moreover, since the
 gauge theory equations amount to a projective flatness condition,
 these moduli spaces correspond with a moduli spaces of flat
 structures. In the case of
 $\U(p,q)$-Higgs bundles, the flat structures correspond to semisimple
 representations of $\pi_1 X$ into the group $\PU(p,q)$. The Higgs
 bundle moduli spaces can thus be used, in a way which we make precise
 in Sections
 \ref{sec:background} and \ref{sec:higgs-bundles}, as a tool to study
 the representation variety
 $$\mathcal{R}(\PU(p,q)) = \Hom^+(\pi_1 X, \PU(p,q)) / \PU(p,q)\ ,$$
 \noindent where $\Hom^+(\pi_1 X, \PU(p,q))$ denotes the semisimple
 representations into $\PU(p,q)$ and the the quotient is by the
 adjoint action.
 
 This relation between Higgs bundles and surface group representations
 has been successfully exploited by others, going back originally to
 the work of Hitchin and Simpson on complex reductive groups. The use
 of Higgs bundle methods to study representation variety
 $\mathcal{R}(G)$ for real $G$ was pioneered by Hitchin in
 \cite{hitchin:1992}, and further developed in 
 \cite{gothen:1995,gothen:2001} and by Xia and Xia-Markman in 
 \cite{xia:1997, xia:1999, xia:2000,markman-xia:2001}. Where we differ
  from these works is that in none
 of them is the general case of $\PU(p,q)$ considered. What we have in
 common is that in all cases insight into the topology of the space
 $\mathcal{R}(G)$ comes from a natural Morse function on the
 corresponding moduli space of Higgs bundles.
 
 The natural Morse function measures the $L^2$-norm of the Higgs field.
 This turns out to provide a suitably non-degenerate Morse function
 which is, moreover, a proper map. In some cases 
 (cf.\ \cite{hitchin:1987,gothen:1995}) all the critical submanifolds are  
 sufficiently well understood so as to permit the extraction of detailed 
  topological information such as the Poincar\'e polynomial.
 In our case our understanding is confined mostly to the submanifolds
 corresponding to the local minima of the Morse function.
Fortunately, this is sufficient for our purposes, namely to 
understand the number of components of the Higgs moduli spaces, and 
thus of the representation varieties. 
 
 The Morse function is non-negative but cannot always attain its zero
 lower bound. For $\GL(n)$-Higgs bundles, this lower bound is
 attained, with the minimizing points in the moduli space consisting
 of semistable vector bundles with zero Higgs field. However in the
 case of $\U(p,q)$-Higgs bundles, the special form of the
 underlying holomorphic bundle prevents a polystable Higgs bundle 
 from having a vanishing Higgs field. The minimizers of the Morse function 
 thus have a more complicated structure than simply that of a stable bundle. 
 Generalizing the results in \cite{gothen:1995,gothen:2001}, we show that
 any minimizer consists of a pair of bundles together with a morphism 
 between them. That is, the minimizers correspond precisely to a special 
 case of the holomorphic triples introduced in 
 \cite{bradlow-garcia-prada:1996}. 
 
 The holomorphic triples admit moduli spaces of stable objects in
 their own right.  In order to exploit the relation between these
 spaces and the minimal submanifolds in the moduli spaces of
 $\U(p,q)$-Higgs bundles, we need a sufficiently good understanding of
 the triples moduli spaces. A substantial part of this paper is
 devoted to acquiring just such an understanding. The way we acquire
 the needed information is similar in spirit to techniques used by
 Thaddeus in \cite{thaddeus:1994}. The key idea (described fully in sections
 \ref{sec:stable-triples} and \ref{crossing-critical-values}) is that the moduli
 spaces of triples come in discrete families, with the members of the
 families ordered by intervals in the range of a continuously varying
 real parameter. As the parameter moves to the large extreme of its
 range, the structure of the corresponding moduli spaces simplify and
 we can obtain a detailed description. Moreover, as the parameter
 decreases, we can track, albeit somewhat crudely, how the moduli
 spaces change. Combining these pieces of data, we get just enough
 information about the moduli space of relevance to our Higgs bundle
 problem.
 
 We now give a brief summary of the contents and main results of this paper.
 
 In Sections \ref{sec:background} and \ref{sec:higgs-bundles} we give
 some background and describe the basic objects of our study.  In
 Sections \ref{sec:background} we describe the natural invariants
 associated with representations of
 $\pi_1 X$ into $\PU(p,q)$. We also discuss the invariants associated
 with representations of $\Gamma$, the universal central extensions of
 $\pi_1$, into $\U(p,q)$. In both cases, these involve a pair of
 integers $(a,b)$ which can be interpreted respectively as degrees of
 rank $p$ and rank $q$ vector bundles over $X$. In the case of the
 $\PU(p,q)$ representations, the pair is well defined only as a class
 in a quotient
 $\Z\oplus\Z/(p,q)\Z$. This leads us to define subspaces
 $\mathcal{R}[a,b]\subset \mathcal{R}(\PU(p,q))$ and
 $\mathcal{R}_{\Gamma}(a,b)\subset \mathcal{R}_{\Gamma}(\U(p,q))$. For
 fixed $(a,b)$, the space $\mathcal{R}_{\Gamma}(a,b)$ fibers over
 $\mathcal{R}[a,b]$ with connected fibers.
 
 In section \ref{sec:higgs-bundles} we give precise definitions of the
 $\U(p,q)$-Higgs bundles and their moduli spaces and establish their
 essential properties. Thinking of a $\U(p,q)$-Higgs bundle as a pair
 $(V\oplus W,\Phi)$, the parameters $(a,b)$ appear here as the degrees
 of the bundles $V$ and $W$. We denote the moduli space of polystable
 $\U(p,q)$-Higgs bundles with $\deg(V)=a$ and $\deg W=b$ by
 $\mathcal{M}(a,b)$, and identify $\mathcal{M}(a,b)$ with the component
 $\mathcal{R}_{\Gamma}(a,b)$ of $\mathcal{R}_{\Gamma}(\U(p,q))$. This,
 together with the fibration over $\mathcal{R}_{\Gamma}(\U(p,q))$ are
 the crucial links between the Higgs moduli and the surface group
 representation varieties.
 
 Except for the last section, where we translate back to the language
 of representation varieties, the rest of the paper is concerned with
 the spaces $\mathcal{M}(a,b)$.  Fixing $p,q,a$ and $b$, we begin the
 Morse theoretic analysis of $\mathcal{M}(a,b)$ in Section
 \ref{sec:morse-theory}. Using the $L^2$-norm of the Higgs field
 $\Phi=\left(
   \begin{smallmatrix}
     0 & \beta \\
     \gamma & 0
   \end{smallmatrix}
   \right)$ as the Morse function, the basic result we need (cf.\ Proposition
   \ref{prop:topology-exercise}) is that 
 this function has a minimum on each connected component of ${\mathcal{M}}(a,b)$,
and if the subspace of local minima is connected then so is 
$\mathcal{M}(a,b)$. The next step is to identify the local minima, 
the loci of which we denote by 
 $\mathcal{N}(a,b)$ . We prove (cf.\ Propositions
 \ref{prop:thm-minima-1} and \ref{lem:thm-minima-4}) that these
 correspond precisely to the polystable Higgs bundles in which
 $\beta=0$ or $\gamma=0$. The data defining a Higgs bundle with
 $\beta=0$ can thus be written as the triple $(W\otimes K, V,\gamma)$.
 Similarly, the $\gamma=0$ minima correspond to triples
 $(V\otimes K, W,\beta)$. This brings us to the theory of such
 holomorphic triples.
 
 In sections \ref{sec:stable-triples}-\ref{sec:n_1=n_2}
 we develop the theory we need concerning holomorphic triples and
 their moduli spaces. While only triples of a specific special kind
 correspond to the minima on the $\U(p,q)$-Higgs moduli, we develop
 the theory for the general case in which a holomorphic triple is
 specified by the set $T=(E_{1},E_{2},\phi)$, where
 $E_{1}$ and $E_{2}$ are holomorphic bundles on $X$ and $\phi \colon
 E_{2} \to E_{1}$ is holomorphic (see \cite{garcia-prada:1994} and
 \cite{bradlow-garcia-prada:1996}).  There is a notion of stability for
 triples which depends on a real parameter $\alpha$ and there are
 moduli spaces of $\alpha$-polystable triples, which are shown in
 \cite{bradlow-garcia-prada:1996}  (see also \cite{garcia-prada:1994})
 to be projective varieties. 
In order for $\mathcal{N}_\alpha$ to be non-empty,
one  must have $\alpha\geq \alpha_m$ with $\alpha_m=d_1/n_1-d_2/n_2\geq 0$. In the
case  $n_1\neq n_2$ there
is also a finite  upper bound $\alpha_M$.
When the parameter $\alpha$ varies, the
 nature of the $\alpha$-stability condition only changes for a
 discrete  number of so-called \emph{critical values} of $\alpha$ (see
 section \ref{sec:triples-definitions} for the precise statements).
 We denote by \begin{displaymath}
   \mathcal{N}_\alpha = \mathcal{N}_\alpha(n_1,n_2,d_1,d_2)
 \end{displaymath}
 the moduli space of
 $\alpha$-polystable triples with $\rk(E_{i}) = n_{i}$ and $\deg(E_{i}) =
 d_{i}$ for $i=1,2$.  The subspace of $\alpha$-stable triples inside
 $\mathcal{N}_\alpha$,  denoted by  $\mathcal{N}_\alpha^s$,  is
 a quasi-projective variety.  
In Theorem \ref{thm:smoothdim} we show that
 
 \begin{theorem*} $\mathcal{N}_\alpha^s$ is smooth for all values
 of $\alpha$ greater than or equal to $2g-2$. 
 \end{theorem*}
 
 \noindent We show furthermore that the triples which appear in
 $\mathcal{N}(a,b)$ are $\alpha$-polystable with $\alpha=2g-2$. We
 must thus understand a moduli space of $\alpha$-stable triples, with
 $\alpha$ on the boundary of the range in which the moduli spaces are smooth.
 We do this indirectly, by obtaining a description of
 $\mathcal{N}_\alpha^s$ when $\alpha$ is large and then examining how
 the moduli space changes as $\alpha$ decreases.
 
 In  section \ref{crossing-critical-values} we examine how the
 moduli spaces differ for values of $\alpha$ on opposite sides of a
 critical value. If $\mathcal{N}_{\alpha_c^{\pm}}$
 denote the moduli spaces for values of $\alpha$ above and below a
 critical value $\alpha_c$, we denote the loci along which they differ
 by $\mathcal{S}_{\alpha_c^{\pm}}$ respectively. Our main results are
 
 \begin{theorem*}[Theorem \ref{thm:codim}]
Let $\alpha_c\in (\alpha_m,\alpha_M)$  be a critical value for 
triples of type $(n_1,n_2,d_1,d_2)$.  If $\alpha_c> 2g-2$ then 
the loci $\mathcal{S}_{\alpha_c^{\pm}} 
\subset\mathcal{N}^s_{\alpha_c^{\pm}}$ are contained in subvarieties of 
codimension at least $g-1$. In particular, they are contained in 
subvarieties of strictly positive codimension if $g\ge 2$. If 
$\alpha_c = 2g-2$ then the same is true for 
$\mathcal{S}_{\alpha_c^+}$.
  \end{theorem*}

 \begin{theorem*}[Corollary \ref{cor:birationality}] Let $\alpha_1$ and 
$\alpha_2$ be any two values in $(\alpha_m,\alpha_M)$ 
such that $\alpha_m <\alpha_1 <\alpha_2 <\alpha_M$ and $2g-2 
\le\alpha_1$. Then 
\begin{itemize}
\item The moduli spaces $\mathcal{N}^s_{\alpha_1}$ and 
$\mathcal{N}^s_{\alpha_2}$ have the same number of connected 
components, and 
\item The moduli space $\mathcal{N}^s_{\alpha_1}$ is irreducible
 if and only if
$\mathcal{N}^s_{\alpha_2}$ is.
\end{itemize} 
 \end{theorem*}
Where by $\alpha_M$ we denote the upper bound for $\alpha$ if $n_1\neq n_2$, or $\infty$
if $n_1=n_2$. 
 
 In sections \ref{sect:special-alpha}-\ref{sec:n_1=n_2}
 we look more closely at how the parameter $\alpha$ affects the nature of $\alpha$-stable triples.
 There are three cases to consider, namely $n_1<n_2$, $n_1>n_2$ and
 $n_1=n_2$. However, using a duality result, it is enough to consider
 $n_1>n_2$ and $n_1=n_2$. In the first case, as mentioned above, there is a bounded interval
 $[\alpha_m,\alpha_M]$ outside of which $\mathcal{N}_\alpha$ is empty.
 Within this interval we identify a number of special values beyond
 which the structure of $\alpha$-stable triples simplify; by Corollary
 \ref{injective-triples} the map $\phi:E_2\rightarrow E_1$ is
 injective if $\alpha>\alpha_0$, by Proposition
 \ref{prop:torsion-degree-bound} the cokernel is torsion free. Finally, for
 the largest values of $\alpha$, i.e.\ for values greater than a bound which
 we denote by $\alpha_L$, we show (cf.\ Proposition
 \ref{prop:E2-stability}) that $\alpha$-stable triples have the form
 $$
 0 \lto E_2 \overset{\phi}{\lto} E_1 \lto F \lto 0,
 $$
 with  $F$ locally free, and $E_2$ and $F$ semistable. This leads to 
 a description of $\mathcal{N}_{\alpha}$ for any $\alpha$ in the 
 range $\alpha_L<\alpha<\alpha_M$. Denoting this moduli space by 
 $\mathcal{N}_L$, we get
 
 \begin{theorem*}[Theorem \ref{thm:largealpha}]
 Let $n_1>n_2$ and $d_1/n_1> d_2/n_2$. 
 
The moduli space 
$\mathcal{N}^s_L(n_1,n_2,d_1,d_2)$ is  smooth, and is birationally 
equivalent to a $\mathbb{P}^N$-fibration over 
$M^s(n_1-n_2,d_1-d_2) \times M^s(n_2,d_2)$, where $M^s(n,d)$ denotes
the moduli space of stable bundles of degree $n$ and rank $d$, and 
the fiber dimension is 
$N=n_2d_1-n_1d_2+n_1(n_1-n_2)(g-1)-1$. In particular, 
$\mathcal{N}_L^s(n_1,n_2,d_1,d_2)$ is non-empty and irreducible. If 
$\GCD(n_1-n_2,d_1-d_2)=1$ and $\GCD(n_2,d_2)=1$, the birational 
equivalence is an isomorphism. 

Moreover, $\mathcal{N}_L(n_1,n_2,d_1,d_2)$ is irreducible and hence 
birationally equivalent to $\mathcal{N}_L^s(n_1,n_2,d_1,d_2)$. 
 \end{theorem*}
 
\begin{theorem*} [Theorem 
\ref{thm:irreducibility-moduli-stable-triples}, 
Corollary \ref{cor:gcd=1}] Let $\alpha$ be any value in the range 
$2g-2\leq\alpha< 
\alpha_M$. Then 
$\mathcal{N}^s_\alpha$ is birationally equivalent to 
$\mathcal{N}^s_L$.  In particular it is non-empty and irreducible. 

Let $(n_1,n_2,d_1,d_2)$ be such that $\GCD(n_2,n_1+n_2,d_1+d_2)=1$. 
If  $\alpha$ is generic then  $\mathcal{N}_\alpha$ is birationally 
equivalent to $\mathcal{N}_L$, and in particular it is irreducible. 
\end{theorem*}

 The case $n_1=n_2$ differs from the $n_1>n_2$ case in two ways. The
 range for $\alpha$ is unbounded above, and in general there is no way
 to avoid torsion in the cokernel of the map $\phi$. The range for
 $\alpha$ presents no difficulties since (cf.\ Theorem
 \ref{thm:stabilization}) beyond a finite bound there are no changes
 in the moduli spaces. It thus still makes sense to identify a `large 
 $\alpha$ moduli space, $\mathcal{N}_L$'.  
We prove the following.

 \begin{theorem*}[Theorem \ref{thm:moduli-n1=n2}]
The moduli space $\mathcal{N}_L(n,n,d_1,d_2)$ is non-empty and irreducible. Moreover, it is
birationally  equivalent to a 
$\mathbb{P}^N$-fibration  over
$M^s(n,d_2)\times \Div^d(X)$, where the fiber dimension is
$N=n(d_1-d_2)-1$.
\end{theorem*}

\begin{theorem*}[Theorem \ref{thm:irreducibility-moduli-stable-triples-n1=n2}] 
If $ \alpha\ge 2g-2$ then the moduli space  
$\mathcal{N}^s_\alpha(n,n,d_1,d_2)$ is  
 birationally equivalent 
to $\mathcal{N}_L(n,n,d_1,d_2)$ and hence non-empty and  irreducible. 
Moreover, $\mathcal{N}_{\alpha}(n,n,d_1,d_2)$ is birationally 
equivalent to $\mathcal{N}_L(n,n,d_1,d_2)$, and hence irreducible, if 
also
\begin{itemize}
\item [$\bullet$] $\GCD(n,2n,d_1+d_2)=1$ and $\alpha\ge 2g-2$ is  generic, or
\item [$\bullet$] $d_1-d_2<\alpha$,
\end{itemize}
\end{theorem*}
 
 In section \ref{sec:main-results} we apply our results to the 
 moduli spaces $\mathcal{M}(a,b)$, and hence to the
 components  $\mathcal{R}_\Gamma(a,b)$ and 
 $\mathcal{R}[a,b]$ of the representation varieties  
 $\mathcal{R}_\Gamma(\U(p,q))$  and $\mathcal{R}(\PU(p,q))$, respectively.
Some of the results depend only on the combination
 $$\tau=\tau(a,b)=2\frac{aq-bp}{p+q}\ ,$$
 known as the Toledo invariant. Indeed, $(a,b)$ is constrained by the bounds 
 $0\le|\tau|\le\tau_M$, 
 where $\tau_M=2\min\{p,q\}(g-1)$. Originally proved by Domic and 
 Toledo in \cite{domic-toledo:1987}, these bounds emerge naturally 
 from our point of view (cf.\ Corollary \ref{cor:toledo} and Remark 
 \ref{remark:MWbound}).  After a discussion (in Section \ref{subs:ab-tau}) 
 of the relation between $(a,b)$ and $\tau$, and 
 (in section \ref{subs: 2g-2}) of the significance of the coprime
 condition  $\GCD(p+q,a+b)=1$, we assemble (in section 
 \ref{sec:results-U(p,q)-Higgs}) our results for the Higgs 
 Moduli spaces. Summarizing the results of section 
 \ref{sec:results-U(p,q)-Higgs} into one Theorem, we get

\begin{theorem*}

Let $(a,b)$ be such that $|\tau(a,b)|\le\tau_M$. Unless further 
restrictions are imposed, let $(p,q)$ be any pair of positive 
integers. 

\begin{itemize}
  
\item[$(1)$] If either of the following sets of conditions apply, 
then the moduli space $\mathcal{M}^s(a,b)$ is a non-empty, smooth 
manifold of the expected dimension, with connected closure 
$\bar{\mathcal{M}}^s(a,b)$: 
\begin{enumerate}
\item[$(i)$] $0<|\tau(a,b)|<\tau_M$ ,
\item[$(ii)$] $|\tau(a,b)|=\tau_M$ and $p=q$
\end{enumerate}

\item[$(2)$] If any one of the 
following sets of conditions apply, then the moduli space 
$\mathcal{M}(a,b)$ is non-empty and connected: 
\begin{enumerate}
\item[$(i)$] $\tau(a,b)=0$,
  
\item[$(ii)$] $|\tau(a,b)|=\tau_M$ and $p\neq q$. , 
  
\item[$(iii)$] $(p-1)(2g-2) < |\tau| \leq \tau_M = p(2g-2)$ and 
$p=q$,

\item[$(iv)$]  $\GCD(p+q,a+b)=1$ 
\end{enumerate}

\item[$(3)$] If  $|\tau(a,b)|=\tau_M$ and $p\neq q$ then
 any element in $\mathcal{M}(a,b)$ is 
strictly semistable (i.e.\ $\mathcal{M}^s(a,b)$ is empty).  If 
$p< q$, then any such representation decomposes as a direct sum of a 
$\U(p,p)$-Higgs bundle  with maximal Toledo invariant and a polystable 
vector bundle of rank $q-p$.  Thus, if $\tau = p(2g-2)$ then there is 
an isomorphism 
  \begin{displaymath}
    \mathcal{M}(p,q,a,b) \cong \mathcal{M}(p,p,a,a - p(2g-2))
    \times M(q-p,b-a + p(2g-2)),
  \end{displaymath}
  where the notation $\mathcal{M}(p,q,a,b)$ indicates
  the moduli space of  $\U(p,q)$-Higgs bundles  with
  invariants $(a,b)$, and $M(n,d)$ is  the moduli space of
semistable vector bundles of rank $n$ and degree $d$. (A similar 
result  holds if $p>q$, as well as if $\tau=-p(2g-2)$). 

\item[$(4)$] If  $\GCD(p+q,a+b)=1$ then $\mathcal{M}(a,b)$
is a smooth manifold of the expected dimension. 
\end{itemize}

\end{theorem*}

Translating this into the language of representations  of $\Gamma$ 
and the fundamental group, we get the following.

\begin{theorem*}[Theorem \ref{thm:results-R-Gamma}]
  
Let $(a,b)$ be such that $|\tau(a,b)|\le\tau_M$. Unless further 
restrictions are imposed, let $(p,q)$ be any pair of positive 
integers. 

\begin{itemize}
  
\item[$(1)$] If either of the following sets of conditions apply, 
then the moduli space $\mathcal{R}_\Gamma^*(a,b)$ of irreducible 
semi-simple representations, is a non-empty, smooth manifold of the 
expected dimension, with connected closure 
$\bar{\mathcal{R}}_\Gamma^*(a,b)$:
\begin{enumerate}
\item[$(i)$] $0<|\tau(a,b)|<\tau_M$ ,
\item[$(ii)$] $|\tau(a,b)|=\tau_M$ and $p=q$
\end{enumerate}

\item[$(2)$] If any one of the 
following sets of conditions apply, then the moduli space 
$\mathcal{R}_\Gamma(a,b)$ of all semi-simple representations is 
non-empty and connected: 
\begin{enumerate}
\item[$(i)$] $\tau(a,b)=0$,
  
\item[$(ii)$] $|\tau(a,b)|=\tau_M$ and $p\neq q$. , 
  
\item[$(iii)$] $(p-1)(2g-2) < |\tau| \leq \tau_M = p(2g-2)$ and 
$p=q$,

\item[$(iv)$]  $\GCD(p+q,a+b)=1$ 
\end{enumerate}

\item[$(3)$] If  $|\tau(a,b)|=\tau_M$ and $p\neq q$ 
then any representation  in $\mathcal{R}_\Gamma(a,b)$ is reducible  
(i.e. $\mathcal{R}^*_\Gamma(a,b)$ is empty).  If $p < q$, then any 
such representation decomposes as a direct sum of a semisimple 
representation of 
$\Gamma$ in $\U(p,p)$ with maximal Toledo invariant and a semisimple 
representation in  $\U(q-p)$.  Thus, if $\tau = p(2g-2)$ then there 
is an isomorphism 
  \begin{displaymath}
    \mathcal{R}_\Gamma(p,q,a,b) \cong \mathcal{R}_\Gamma(p,p,a,a - p(2g-2))
    \times R_{\Gamma}(q-p,b-a + p(2g-2)),
  \end{displaymath}
  where the notation $\mathcal{R}_\Gamma(p,q,a,b)$ indicates
  the moduli space of representations of $\Gamma$ in $\U(p,q)$ with
  invariants $(a,b)$, and $R_{\Gamma}(n,d)$ denotes the moduli space of
  degree $d$ representations of $\Gamma$ in $\U(n)$. (A similar result
  holds if $p>q$, as well as if $\tau=-p(2g-2)$).
\item[$(4)$] If  $\GCD(p+q,a+b)=1$ then $\mathcal{R}_\Gamma(a,b)$ 
is a smooth manifold of the expected dimension. 
\end{itemize}
\end{theorem*}

\begin{theorem*}[Theorem \ref{thm:results-R}]

Let $(a,b)$ be such that $|\tau(a,b)|\le\tau_M$. Unless further 
restrictions are imposed, let $(p,q)$ be any pair of positive 
integers. 

\begin{itemize}
  
\item[$(1)$] If 
either of the following sets of conditions apply, then the moduli 
space $\mathcal{R}^*[a,b]$ of irreducible semi-simple 
representations, is non-empty, with connected closure 
$\bar{\mathcal{R}}^*[a,b]$: 
\begin{enumerate}
\item[$(i)$] $0<|\tau(a,b)|<\tau_M$ ,
\item[$(ii)$] $|\tau(a,b)|=\tau_M$ and $p=q$
\end{enumerate}

\item[$(2)$] If any one of the 
following sets of conditions apply, then the moduli space 
$\mathcal{R}[a,b]$ of all semi-simple representations is non-empty 
and connected: 
\begin{enumerate}
\item[$(i)$] $\tau(a,b)=0$,
  
\item[$(ii)$] $|\tau(a,b)|=\tau_M$ and $p\neq q$. , 
  
\item[$(iii)$] $(p-1)(2g-2) < |\tau| \leq \tau_M = p(2g-2)$ and 
$p=q$,

\item[$(iv)$]  $\GCD(p+q,a+b)=1$ 
\end{enumerate}

\item[$(3)$] If  $|\tau(a,b)|=\tau_M$ and $p\neq q$ 
then any representation  in $\mathcal{R}[a,b]$ is reducible  (i.e.\ 
$\mathcal{R}^*[a,b]$ is empty).  If $p < q$, then any such representation
reduces to a semisimple  representation of 
$\pi_1 X$  in  
  $\mathrm{P}(\U(p,p) \times \U(q-p))$, such that the representation in
  $\PU(p,p)$ induced via projection on the first factor has maximal
  Toledo invariant.  (A similar result
  holds if $p>q$, as well as if $\tau=-p(2g-2)$).
\end{itemize}
\end{theorem*}

Statement (3) in the previous theorem is a generalization to 
arbitrary $(p,q)$ of a  result of D. Toledo \cite{toledo:1989}  when 
$p=1$ and L. Hern\'andez  
\cite{hernandez:1991} when $p=2$. This rigidity phenomenon  for the moduli space
of representations for the largest value of the Toledo invariant turns out
to be of significance in relation to Hitchin's Teichm\"uller components 
for the real split form of a complex group \cite{hitchin:1992} (this will be
discussed somewhere else \cite{garcia-ramanan:inpreparation}).

We note, finally, that our methods clearly have wider applicability 
than to the $\U(p,q)$-Higgs bundles and representations into 
$\PU(p,q)$. A careful scrutiny of the Lie algebra properties used in our 
proofs suggests a generalization to any real group $G$ for which 
$G/K$ is hermitian symmetric, where $K\subset G$ is a maximal compact subgroup. 
This will be addressed in a future publication.

The main results proved in this paper were announced in the note 
\cite{bradlow-garcia-gothen:2001}.  In that note we claim that the 
connectedness results hold for the moduli spaces $\mathcal{R}(a,b)$ 
and $\mathcal{R}[a,b]$, whether or not the coprimality condition 
$\GCD(p+q,a+b) = 1$ is satisfied (and similarly for the corresponding moduli 
of triples and $\U(p,q)$-Higgs bundles).  While we expect this
to be true, we have not so far been able to prove it.  We hope to 
come back to this question in a future publication.

\subparagraph{Acknowledgements.} We thank the mathematics departments
of the University of Illinois at Urbana-Champaign, the Universidad
Aut{\'o}noma de Madrid and the University of Aarhus, the Department of
Pure Mathematics of the University of Porto, the Mathematical Sciences
Research Institute of Berkeley and the Mathematical Institute of the
University of Oxford for their hospitality during various stages of
this research.  We thank Ron Donagi, Bill Goldman, Tom\'as G\'omez,
Rafael Hern\'andez, Nigel Hitchin, Alastair King, Eyal Markman,
Vicente Mu\~noz, Peter
Newstead, S. Ramanan, Domingo Toledo, and Eugene Xia, for many
insights and patient explanations.

%%%%%%%%%%%%%%%%%%%%%%%%%%%%%%%%%%%%%%%%%%%%%%%%%%%%%%%%%%%%%%
\section{Representations of surface groups}\label{sec:background}
%%%%%%%%%%%%%%%%%%%%%%%%%%%%%%%%%%%%%%%%%%%%%%%%%%%%%%%%%%%%%%
%%%%%%%%%%%%%%%%%%%%%%%%%
\subsection{Definitions}
%%%%%%%%%%%%%%%%%%%%%%%%%

Let $X$ be a closed oriented surface of genus $g \geq 2$ and let $G$ 
be either $\U(p,q)$ or $\PU(p,q)$ where $p$ and $q$ are any positive 
integers.  We think of $\U(p,q)$ as the subgroup of $\GL(n)$ (with 
$n=p+q$) which leaves invariant a hermitian form of signature 
$(p,q)$. It is a non-compact real form of 
$\GL(n)$ with center $S^1$ and maximal compact subgroup $\U(p)\times\U(q)$.
The quotient $\U(p,q)/\U(p)\times\U(q)$ is a hermitian symmetric 
space. The adjoint form $\PU(p,q)$ is given by the exacts sequence of 
groups 
$$
1\longrightarrow \U(1)\longrightarrow 
\U(p,q)\longrightarrow\PU(p,q)\longrightarrow 1.
$$
\noindent By a representation of $\pi_1 X$ in $G$ we mean a homomorphism 
$\rho \colon \pi_1 X \to G$. Fixing $\PU(p,q)\subset \PGL(n)$, we say a 
representation of $\pi_1 X$ in $\PU(p,q)$ is semi-simple if it 
defines a semi-simple $\PGL(n)$ representation. The group $\PU(p,q)$ 
acts on the set of representations via conjugation. Restricting to 
the semi-simple representations, we get the character variety, 
\begin{equation}\label{eqn:RGdef}
  \mathcal{R}(\PU(p,q)) = \Hom^+(\pi_1 X, \PU(p,q)) / \PU(p,q).
\end{equation}
This can be described as follows: from the standard presentation 
\begin{displaymath}
  \pi_1 X = \langle A_{1},B_{1}, \ldots, A_{g},B_{g} \;|\;
  \prod_{i=1}^{g}[A_{i},B_{i}] = 1 \rangle
\end{displaymath}
we see that $\Hom^{+}(\pi_1 X, \PU(p,q))$ can be embedded in 
$\PU(p,q)^{2g}$ via 
\begin{align*}
  \Hom^{+}(\pi_1 X, \PU(p,q)) &\to \PU(p,q)^{2g} \\
  \rho &\mapsto (\rho(A_1), \ldots \rho(B_g)).
\end{align*}
We give $\Hom^{+}(\pi_1 X, \PU(p,q))$ the subspace topology and 
$\mathcal{R}(\PU(p,q))$ the quotient topology; this is Hausdorff because we
have restricted attention to semi-simple representations. 
We can similarly define 
 \begin{equation}\label{eq:RgGdef}
  \mathcal{R}_\Gamma(\U(p,q)) = \Hom^+(\Gamma, \U(p,q)) / \U(p,q),
\end{equation}
where $\Gamma$ is the central extension 
\begin{equation}\label{eq:gamma}
0\longrightarrow\mathbb{Z}\longrightarrow\Gamma\longrightarrow\pi_1 
X\longrightarrow 1\ 
\end{equation} 
\noindent defined (as in \cite{atiyah-bott:1982})
by the generators $A_{1},B_{1}, \ldots, A_{g},B_{g}$ and a central 
element $J$ subject to the relation $\prod_{i=1}^{g}[A_{i},B_{i}] = 
J$.  Regarding $\U(p,q)$ as a subset of $\GL(n)$, the representations 
in $\Hom^+(\Gamma, \U(p,q))$ are direct sums of irreducible 
representations on $\C^n$. 
The first step in the study of the topological properties of 
$\mathcal{R}(G)$ is to identify the appropriate topological invariants
of a representation $\rho \colon \pi_1 X \to G$. To do that, one uses
the correspondence between representations of $\pi_1 X$ in $G$ and
flat principal $G$-bundles on $X$. We  start with $G=\PU(p,q)$.
Let $\rho \colon \pi_1 X \to \PU(p,q)$ be a representation.  The 
corresponding flat principal $\PU(p,q)$-bundle is 
\begin{displaymath}
  P_{\rho} = \tilde{X} \times_{\rho} \PU(p,q),
\end{displaymath}
where $\tilde{X}$ is the universal cover of $X$. Since $X$ has real 
dimension two, any $\PU(p,q)$-bundle lifts to a $\U(p,q)$-bundle. 
Moreover, a 
$\PU(p,q)$-bundle with a flat connection can be lifted to a 
$\U(p,q)$-bundle with a projectively flat connection, i.e.\ with
a connection with constant central curvature.  Now, for any smooth 
(not necessarily flat) $\U(p,q)$-bundle there is a reduction of  
structure group to the maximal compact subgroup 
$\U(p) \times \U(q)$. Taking the standard representation on $\mathbb{C}^p\oplus
\mathbb{C}^q$, we get an associated vector bundle of the form $V\oplus W$, where
$V$ and $W$ are rank $p$ and $q$ complex vector bundles respectively.
Such bundles over a Riemann surface are topologically classified by a 
pair of integers 
\begin{displaymath}
  (a,b) = (\deg(V),\deg(W)).
\end{displaymath}
The lift to a $\U(p,q)$-bundle, and therefore the pair $(a,b)$, is 
however not uniquely determined. If we twist the associated vector 
bundle (plus projectively flat connection) by a line bundle $L$ with 
a connection with constant curvature, then after projectivizing we 
obtain the same flat 
$\PU(p,q)$-bundle. If the degree of $L$ is $l$ then the invariant 
associated to the twisted bundle is $(a+pl,b+ql)$.  In order to 
obtain a well defined invariant for the representation $\rho$ we must 
thus take the quotient of $(\Z 
\oplus 
\Z)$ by the 
$\Z$-action $l\cdot (a,b)=(a+pl,b+ql)$, i.e.\ we must pass to 
the quotient $(\Z \oplus \Z) / (p,q)\Z$ in the exact sequence
$$
0\rightarrow\Z \rightarrow \Z \oplus \Z 
\rightarrow(\Z \oplus \Z) / (p,q)\Z\rightarrow 0\ .
$$
\noindent Since the $\PU(p,q)$-orbits in $\Hom(\pi_1 X, \PU(p,q))$ 
under the conjugation action correspond to isomorphism classes of 
flat $\PU(p,q)$-bundles, the above construction defines a map 
\begin{equation}\label{eq:invariants}
c: \mathcal{R}(\PU(p,q))\longrightarrow (\Z \oplus \Z) / (p,q).
\end{equation} 
\noindent The map is continuous and is thus constant on connected 
components of $\mathcal{R}(\PU(p,q))$. 
\begin{remark} This map can be seen from a different point of view,
from which it seen that the target space is $\pi_{1} \PU(p,q)$. We 
begin with the observation that the flat bundle 
$P_{\rho}$ is described by locally constant transition functions. 
Thus the isomorphism class of this bundle is represented by a class 
in the (non-abelian) cohomology set $H^1(X, \PU(p,q))$, where, by 
abuse of notation, we denote the sheaf of locally constant maps into 
$\PU(p,q)$ on $X$ by the same symbol $\PU(p,q)$.  Let 
$\widetilde{\PU}(p,q)$ be the universal cover of $\PU(p,q)$.  The 
short exact sequence of groups 
\begin{displaymath}
  \pi_1 \PU(p,q) \to \widetilde{\PU}(p,q) \to \PU(p,q)
\end{displaymath}
induces a sequence of cohomology sets and, since $\pi_1 
\PU(p,q)$ is Abelian, the coboundary map 
\begin{displaymath}
  \delta\colon H^1(X, \PU(p,q)) \to H^2(X, \pi_1 \PU(p,q))
\end{displaymath}
can be defined.  The obstruction to lifting the flat 
$\PU(p,q)$-bundle $P_{\rho}$ to a flat $\widetilde{\PU}(p,q)$-bundle is exactly the
image of the cohomology class of $P_{\rho}$ under $\delta$.  We 
denote this class by $ c(\rho) \in H^2(X, \pi_1 \PU(p,q)) 
\cong \pi_1 \PU(p,q)$.
Next we recall the calculation of $\pi_{1} \PU(p,q)$. The maximal 
compact subgroup of $\U(p,q)$ is $\U(p) \times \U(q)$ and the 
inclusion $\U(p) \times \U(q) \into \U(p,q)$ is a homotopy 
equivalence.  The determinant gives an isomorphism of fundamental 
groups $\pi_1 \U(p) \xrightarrow{\cong} \pi_1 \U(1) \cong \Z$.  Hence 
the map $ \U(p,q) \to \U(1) \times \U(1)$ defined by 
\begin{equation}
  \label{eq:upq-det}
  \begin{split}
  \U(p,q) &\to \U(1) \times \U(1) \\
  \begin{pmatrix}
    x & y \\
    z & w
  \end{pmatrix}
  &\mapsto (\det(x),\det(w))
  \end{split}
\end{equation}
gives an isomorphism $\pi_1 \U(p,q) \xrightarrow{\cong} \Z \oplus 
\Z$. Furthermore, the composition of the standard inclusion $\U(1) 
\into \U(p,q)$ and the map given in \eqref{eq:upq-det} is the map
$\lambda \mapsto (\lambda^{p},\lambda^{q})$ from 
$\U(1) \to \U(1) \times \U(1)$. The induced map on fundamental groups is 
$n \mapsto (pn,qn)$. The short exact sequence 
\begin{equation}\label{ses:S1UPU}
  \U(1) \to \U(p,q) \to \PU(p,q)
\end{equation}
is a fibration, so we see that $\pi_1 \U(p,q) \to \pi_1 \PU(p,q)$ is 
surjective.  It follows that we have a commutative diagram 
\begin{displaymath}
  \begin{CD}
    \pi_1\U(1) @>>> \pi_1\U(p,q) @>>> \pi_1\PU(p,q) \\
    @V\cong VV @V\cong VV @V\cong VV \\
    \Z @>(p\cdot,q\cdot)>> \Z \oplus \Z @>>> (\Z \oplus \Z) / (p,q)\Z
  \end{CD}\,
\end{displaymath}
and hence $c(\rho)$ defines a class 
$[a,b] \in (\Z \oplus \Z)/(p,q)\Z$. This is the same class as that defined by the map
(\ref{eq:invariants}) though, since we will not make use of this, we 
omit the proof. 
\end{remark}
%%%%%%%%%%%%%%%%%%%%
\subsection{Invariants on $\mathcal{R}_{\Gamma}(\U(p,q))$ and 
relation to $\mathcal{R}(\PU(p,q))$}
%%%%%%%%%%%%%%%%%
Putting together (\ref{ses:S1UPU}) and (\ref{eq:gamma}) we get the 
commutative diagram 
\begin{displaymath}
  \begin{CD}
   \U(1) @>>> \U(p,q) @>>>\PU(p,q)\\
    @AAA @A\tilde{\rho}AA@A\rho AA \\
    \Z @>>> \Gamma @>>> \pi_1 X
  \end{CD}\ .
\end{displaymath}
\noindent {}From this we get a surjection 
$\pi:\mathcal{R}_{\Gamma}(\U(p,q))\rightarrow\mathcal{R}(\PU(p,q))$. We can 
understand the fibers of this map as follows. By the same argument as 
in \cite{atiyah-bott:1982}\footnotemark
\footnotetext{While \cite{atiyah-bott:1982} gives the argument for
$\U(n)$ and $\PU(n)$, there are no essential changes to be made in order
to adapt for the case of $\U(p,q)$ and $\PU(p,q)$.}, 
$\mathcal{R}_{\Gamma}(\U(p,q))$ can be identified as the moduli space of $\U(p,q)$-bundles on $X$ with projectively flat structures. Taking the reduction to
the maximal compact $\U(p)\times\U(q)$, we thus associate to each 
class $\tilde{\rho}\in \mathcal{R}_{\Gamma}(\U(p,q))$ a vector bundle 
of the form $V\oplus W$, where $V$ and $W$ are rank $p$ and $q$ 
respectively, and thus a pair of integers $(a,b) =(\deg(V),\deg(W))$. 
The map $\tilde{c}:\tilde{\rho}\mapsto (a,b)$ fits in a commutative 
diagram 
\begin{displaymath}
  \begin{CD}
    \mathcal{R}_{\Gamma}(\U(p,q)) @> \pi >> \mathcal{R}(\PU(p,q))  \\
    @V \tilde{c} VV @V c VV  \\
    \Z\oplus\Z  @>>> (\Z \oplus \Z) / (p,q)\Z
  \end{CD}\ .
\end{displaymath}
\noindent 
\noindent We can now define the subspaces 
\begin{align*}
\mathcal{R}_\Gamma(a,b) :&= \tilde{c}^{-1}(a,b)\\
&= \{ \tilde{\rho} \in \mathcal{R}_\Gamma(\U(p,q))\;\; | \;\;
\tilde{c}(\tilde{\rho})=(a,b) \in \Z \oplus \Z\}, \\
\mathcal{R}[a,b] :&= c^{-1}[a,b]\\
&=  \{ \rho \in \mathcal{R}(\PU(p,q))\;\; |\;\;
c(\rho)=[a,b] \in \Z \oplus\Z/(p,q)\Z \}\ .
\end{align*}
\noindent Clearly we have surjective maps 
\begin{equation}
  \label{eq:principal-jac}
  \mathcal{R}_\Gamma(a,b) \to \mathcal{R}[a,b].
\end{equation}
\noindent Moreover, the pre-image 
\begin{equation}
\pi^{-1}(\mathcal{R}[a,b])=\bigcup_{(a,b)}\mathcal{R}_\Gamma(a,b) 
\end{equation}
\noindent where the union is over all $(a,b)$ in the class
$[a,b]\in \Z \oplus\Z/(p,q)\Z$. As mentioned above, 
tensoring by line bundles with constant curvature 
connections of degree $l$ gives an isomorphism 
\begin{displaymath}
  \mathcal{R}_\Gamma(a,b)
  \xrightarrow{\cong}
  \mathcal{R}_\Gamma(a+pl,b+ql)\ .
\end{displaymath}
Notice that if the invariant $c(\rho)$ of a representation 
$\rho\in\mathcal{R}(\PU(p,q))$ can be represented by the pair 
$(a,-a)$,  then the associated $\U(p,q)$-bundle has degree zero and the 
projectively flat connection is actually flat.  Thus 
$\rho$ defines to a representation of $\pi_1 X$ in $\U(p,q)$. Under 
the correspondence between $\mathcal{R}(\PU(p,q))$ and 
$\mathcal{R}_{\Gamma}(\U(p,q))$, $\rho$ corresponds to a $\Gamma$ 
representation in which the central element $J$ acts trivially. 
Furthermore the subspaces 
\begin{displaymath}
\mathcal{R}(a) = \{ \rho \in \mathcal{R}(\U(p,q))\;\;|\;\;
c(\rho)=(a,-a)\;\;\mbox{with}\;\; a\in \Z \} 
\end{displaymath}
can be identified with the subspaces 
$\mathcal{R}_\Gamma(a,-a)\subset \mathcal{R}_{\Gamma}(\U(p,q))$. 
Finally, we observe that $\mathrm{Jac}(X)$, the moduli space of flat 
degree zero line bundles, acts by tensor product of bundles on 
$\mathcal{R}_\Gamma(a,b)$.  Since $\mathrm{Jac}(X)$ is isomorphic to 
the torus $\U(1)^{2g}$, we get the following relation between 
connected components. 

\begin{proposition}\label{prop:principal-jac} The map 
$ \mathcal{R}_\Gamma(a,b) \to \mathcal{R}[a,b]$ given in 
\eqref{eq:principal-jac}
defines a $\U(1)^{2g}$-fibration which, if the total space and base 
are smooth manifolds, is a smooth principal bundle. Thus the subspace 
$\mathcal{R}[a,b] \subseteq \mathcal{R}(\PU(p,q))$ is 
connected if $\mathcal{R}_\Gamma(a,b)$ is connected. 
\hfill\qed 
\end{proposition}

We will study 
$\mathcal{R}_\Gamma(a,b)$ by choosing a complex structure on 
$X$ and identifying this space with a certain moduli space of Higgs bundles.
This is carried out in the next section.  In the rest of the paper, 
the subspaces of irreducible representations are denoted by 
$\mathcal{R}^*$. 
 
%%%%%%%%%%%%%%%%%%%%%%%%%%%%%%%%%%%%%%%%%%%%%%%%
\section{Higgs bundles and flat connections}
\label{sec:higgs-bundles}
%%%%%%%%%%%%%%%%%%%%%%%%%%%%%%%%%%%%%%%%%%%%%%%%
%%%%%%%%%
\subsection{$\GL(n)$-Higgs bundles}
%%%%%%%%%
Give $X$ the structure of a Riemann surface.  We recall (from
\cite
{donaldson:1987,corlette:1988,hitchin:1987,simpson:1988,simpson:1994a,simpson:1994b})
the following basic facts about $\GL(n)$-Higgs bundles. 
\begin{definition}
\label{defn:GLnHiggs}
\begin{enumerate}
\item A $\GL(n)$-Higgs bundle on
$X$ is a pair $(E,\Phi)$, where $E$  is a rank $n$ holomorphic vector
bundle over $X$ and $\Phi \in H^0(\End(E) \otimes K)$ is a 
holomorphic endomorphism of $E$ twisted by the canonical bundle $K$ 
of $X$. 
\item The $\GL(n)$-Higgs bundle $(E,\Phi)$ is \emph{stable} if the  slope 
stability condition 
\begin{equation}\label{eq:stability}
\mu(E') < \mu(E)
\end{equation} 
\noindent holds for all proper 
$\Phi$-invariant subbundles $E'$ of $E$.  Here the slope is defined by
$\mu(E)=\deg(E)/\rk(E)$ and $\Phi$-invariance means that $\Phi(E')\subset E'\otimes K$. 
\item Semistability is defined by replacing the above strict inequality with a weak 
inequality. A Higgs bundle is called polystable if is the direct sum 
of stable Higgs bundles with the same slope. 
\item Given a hermitian metric $H$ on $E$, let $A$ denote the unique 
connection compatible with the holomorphic structure and unitary with 
respect to $H$. Hitchin's equations on $(E,\Phi)$ are
\begin{equation}
  \label{eq:hitchin1}
  \begin{aligned}
  F_A + [\Phi,\Phi^*] &= -\sqrt{-1}\mu \text{Id}_E \omega, \\
  \dbar_{A} \Phi &=0, 
  \end{aligned}
\end{equation}
where $\omega$ is the K\"ahler form on $X$, $\text{Id}_E$ is the
identity on $E$, $\mu = \mu(E)$ and $\dbar_A$ is the antiholomorphic 
part of the covariant derivative $d_A$ . 
\end{enumerate}
\end{definition}
\begin{proposition}
\cite{donaldson:1987,hitchin:1987,simpson:1988,simpson:1994a,simpson:1994b}
 \label{prop:HK}
\begin{enumerate}
\item Let $(E,\Phi)$ be a $\GL(n)$-Higgs bundle. Then $(E,\Phi)$ is 
polystable if and only if it admits a hermitian metric such that 
Hitchin's equation (\ref{eq:hitchin1}) is satisfied 
\item There is a moduli space of  rank $n$ 
degree $d$ polystable Higgs bundles which is a quasi-projective 
variety of complex dimension $2(d+n^2(g-1))$. 
\item If we define a Higgs connection (as in \cite{simpson:1994a}) by 
\begin{equation}\label{eqn:Dhiggs} 
D=d_A+\theta
\end{equation}
\noindent where $\theta=\Phi+\Phi^* $, then Hitchin's equations are 
equivalent to the conditions
\begin{equation}
  \label{eq:harmonic1}
  \begin{aligned}
  F_D =& -\sqrt{-1}\mu \text{Id}_E \omega, \\
  &d_A \theta =0, \\
  &d_A^* \theta =0.
  \end{aligned}
\end{equation}
\item In particular, since $X$ is a Riemann surface, if $A$ satisfies 
(\ref{eq:hitchin1}) then $D$ is a projectively flat connection. If 
$\deg(E)=0$ then $D$ is actually flat. It follows that in this 
case the pair $(E,D)$ defines a representation of $\pi_1 X$ in 
$\GL(n)$. If $\deg(E)\ne 0$, then the pair $(E,D)$ defines a 
representation of $\pi_1 X$ in 
$\PGL(n)$, or equivalently, a representation of $\Gamma$ in $\GL(n)$.
By the theorem of Corlette (\cite{corlette:1988}), every semisimple
representation of $\Gamma$ (and therefore all semisimple
representation of $\pi_1 X$) arise in this way.
\item There is thus a bijective correspondence between the moduli space of 
polystable Higgs bundles of rank $n$ and the moduli space of 
(conjugacy classes of) semisimple representations of $\Gamma$ in 
$\GL(n)$. If the  degree of the Higgs bundle is zero, then the first
moduli space corresponds actually to the representation variety for 
representations of $\pi_1 X$ in $\GL(n)$. 
\end{enumerate}
\end{proposition}
%%%%%%%%%
\subsection{$\U(p,q)$-Higgs bundles}
%%%%%%%%%
If we fix integers $p,q$ such that $n=p+q$, then we can isolate a 
special class of $\GL(n)$-Higgs bundles by the requirements that
\begin{equation} \label{upq-higgs-bundle}
  \begin{aligned}
  E &= V \oplus W \\
  \Phi &=
  \left(
  \begin{smallmatrix}
    0 & \beta \\
    \gamma & 0
  \end{smallmatrix}
  \right)
  \end{aligned}
\end{equation}
where $V$ and $W$ are holomorphic vector bundles on $X$ with 
$\rk(V) = p$, $\rk(W) = q$, $\deg(V) = a$,  $\deg(W) = b$,
$\beta \in H^0(\Hom(W,V) \otimes K)$, and 
$\gamma \in H^0(\Hom(V,W) \otimes K)$.
We can describe such Higgs bundles more intrinsically as follows. Let 
$P_{\GL(p)}$ and $P_{\GL(q)}$ be the principal frame bundles for $V$ and $W$ respectively.  Let $P=P_{\GL(p)}\times P_{\GL(q)}$ be the fiber product,
and let $\Ad P=P\times_{\Ad}\mathfrak{gl}(n)$ be the adjoint bundle, 
where 
$\GL(p)\times\GL(q)\subset\GL(n)$ acts by the Ad-action on the the 
Lie algebra of $\GL(n)$.  Let 
$(\mathfrak{gl}(p)\oplus\mathfrak{gl}(q))^{\perp}\subset\mathfrak{gl}(n)$
be the orthogonal complement with respect to the usual inner product. 
This defines a subbundle 
\begin{equation}\label{eq:Ppq}
P_{p,q}:= P\times_{\Ad}(\mathfrak{gl}(p)\oplus\mathfrak{gl}(q))^{\perp}
\subset \Ad P\ .
\end{equation}
\noindent We can then make the following definition.
\begin{definition} 
\label{defn:upq}
A $\U(p,q)$-Higgs bundle\footnotemark
\footnotetext{The reason for the name is 
explained by the following remarks and by Lemma \ref{lemma:upq}} on 
$X$ is a pair $(P,\Phi)$ where 
$P$ is a holomorphic principal $\GL(p)\times\GL(q)$ bundle, and 
$\Phi$ is a holomorphic section of the vector bundle 
$P_{p,q}\otimes K$ (where $P_{p,q}$ is the bundle
defined in (\ref{eq:Ppq}).
\end{definition}
\begin{remark} We can always write $P=P_{\GL(p)}\times P_{\GL(q)}$. 
If we let $V$ and $W$ be the standard vector bundles associated to 
$P_{\GL(p)}$ and $P_{\GL(q)}$ respectively, then any 
$\Phi\in H^0(P_{p,q}\otimes K)$ can be written as
in (\ref{upq-higgs-bundle}). We will usually adopt the vector bundle 
description of $\U(p,q)$-Higgs bundles. 
\end{remark} 
\begin{remark}\label{rem:HiggsG}  Definition \ref{defn:upq} is compatible
with the definitions in \cite{hitchin:1992} and in 
\cite{gothen:1995}. There they define a $G$-Higgs bundle for any real 
form of a complex reductive Lie group. The bundle in their definition 
is a principal 
$H^{\mathbb{C}}$-bundle, where $H\subset G$ is a maximal compact subgroup 
and $H^{\mathbb{C}}$ is its complexification. Thus in the case that 
$G=\U(p,q)$, we get that  $H^{\mathbb{C}}=\GL(p)\times\GL(q)$. The Higgs 
field is precisely a section of the bundles which appears in 
Definition \ref{defn:upq}. {}From a different perspective, Definition 
\ref{defn:upq} defines an example of a principal pair in the sense of 
\cite{banfield} and 
\cite{mundet}. Strictly speaking, since the canonical bundle $K$ 
plays the role of a fixed `twisting bundle', what we get is a 
principal pair in the sense of \cite{BGM}. The defining data for the 
pair are then (i) the principal 
$\GL(p)\times\GL(q)\times\GL(1)$-bundle 
$P_{\GL(p)}\times P_{\GL(q)}\times P_K$, where $P_K$ is the frame bundle
for $K$ and (ii) the associated vector bundle $P_{p,q}\otimes K$.
\end{remark}
\begin{lemma}\label{lemma:upq} Let $(V\oplus W,\Phi)$ be a $\U(p,q)$-Higgs bundle
 with a Hermitian 
metric $H=H_V\oplus H_W$, i.e.\ such that  $V\oplus W$ is a unitary
decomposition. Let $A$ be a unitary connection with respect to $H$, 
and let $D=d_A+\theta$ be the corresponding Higgs connection, where 
$\theta=\Phi+\Phi^*$. Then $D$ is a $\U(p,q)$-connection, i.e.\ in any 
unitary local frame the connection 1-form takes its values in the Lie 
algebra of $\U(p,q)$. 
\end{lemma}
\begin{proof} Fix a local unitary frame (with respect to $H=H_V\oplus 
H_W$). Then $D=d+A+\theta$, where $A$ takes its values in $\mathfrak{u}(p)\oplus 
\mathfrak{u}(q)\subset \mathfrak{u}(p,q)\subset \mathfrak{u}(n)$, while $\theta$ takes its 
values in 
$(\mathfrak{u}(p)\oplus \mathfrak{u}(q))^\perp \cap \mathfrak{u}(p,q)$.
 \end{proof}
\begin{definition}\label{defn:U(pq)Higgs}
Let $(E,\Phi)$ be a $\U(p,q)$-Higgs bundle
with $E=V\oplus W$ and $\Phi=\left(
  \begin{smallmatrix}
    0 & \beta \\
    \gamma & 0
  \end{smallmatrix}\right)$. 
We say $(E,\Phi)$ is a {\it stable} $\U(p,q)$-Higgs bundle if the  
slope stability condition (\ref{eq:stability}), i.e.\ 
$\mu(E') < \mu(E)$, is satisfied for all 
$\Phi$-invariant subbundles of the form $E'=V'\oplus W'$, i.e.\ for all
subbundles $V'\subset V$ and 
$W'\subset W$ such that 
\begin{align}
\beta &:W'\longrightarrow V'\otimes K\\
\gamma &:V'\longrightarrow W'\otimes K\ .
\end{align}
\noindent Semistability for $\U(p,q)$-Higgs bundles is defined by replacing the 
above strict inequality with a weak inequality, and polystability 
means a direct sum of stable $\U(p,q)$-Higgs bundles all with the same 
slope. In particular a polystable $\U(p,q)$-Higgs bundle is the 
direct sum of (lower rank) 
$\U(p',q')$-Higgs bundles.  
\noindent We shall say that a polystable $\U(p,q)$-Higgs bundle which is not 
stable is \emph{reducible}. 
\noindent Two $\U(p,q)$-Higgs bundles $(V\oplus W,\Phi)$ and $(V'\oplus W',\Phi')$
are isomorphic if there are isomorphisms $g_V:V\rightarrow V'$ and 
$g_W:W\rightarrow W'$ which intertwine $\Phi$ and $\Phi'$, i.e.\ 
such that $(g_V\oplus g_W)\otimes I_K\circ\Phi=\Phi'\circ(g_V\oplus 
g_W)$ where $I_K$ is the identity on $K$.
\end{definition}
\begin{remark}\label{rem:stability-equiv} The  stability condition for a $\U(p,q)$-Higgs bundle is {\em a priori} weaker than the stability condition 
given in Definition \ref{defn:GLnHiggs} for $\GL(n)$-Higgs bundles. 
Namely, the slope condition has to be  satisfied only for all proper 
non-zero $\Phi$-invariant subbundles which respect the decomposition 
$E = V \oplus W$, that is, subbundles of the form $E' = V' \oplus W'$
with $V' \subseteq V$ and $W' \subseteq W$.  However, it is shown in  
\cite[Section 2.3]{gothen:2001} that the weaker condition is in fact 
equivalent to the ordinary stability of $(E,\Phi)$. 
\end{remark}
\begin{proposition}
 \label{prop:U(pq)HK}
Let $(E,\Phi)$ be a $\U(p,q)$-Higgs bundle with $E=V\oplus W$ and 
$\Phi=\left( 
  \begin{smallmatrix}
    0 & \beta \\
    \gamma & 0
  \end{smallmatrix}\right)$. Then $(E,\Phi)$ is $\U(p,q)$-polystable 
  if and only if it admits a compatible hermitian metric $H=H_V\oplus H_W$ such 
  that Hitchin's equation (\ref{eq:hitchin1}) is satisfied
\end{proposition}
\begin{proof} Though not explicitly proved there, this is a special case of
the correspondence invoked in \cite{hitchin:1992} for $G$-Higgs 
bundles where $G$ is a real form of a reductive Lie group. By Remark 
\ref{rem:HiggsG}  it can also be seen as a special case of the 
Hitchin--Kobayashi correspondence for principal pairs (cf.\ 
\cite{banfield} and 
\cite{mundet} and \cite{BGM}).  We note finally that in one direction the result
follows immediately from Remark \ref{rem:stability-equiv} : if 
$(V\oplus W, \Phi)$ supports a compatible metric such that 
(\ref{eq:hitchin1}) is satisfied, then it is polystable as a 
$\GL(n)$-Higgs bundle, and hence it is $\U(p,q)$-polystable. 
\end{proof}
\begin{remark}This correspondence allows us, via the next theorem, to use 
$\U(p,q)$-Higgs bundles to study representations of the surface 
groups $\pi_1 X$ and $\Gamma$ into $\U(p,q)$ and $\PU(p,q)$
\end{remark}

\begin{definition}Fix integers $a$ and $b$. Let $\mathcal{M}(a,b)$ 
denote the moduli space of isomorphism classes of polystable 
$\U(p,q)$-Higgs bundles with $\deg(V)=a$ and $\deg W=b$. 
\end{definition}

\begin{proposition}\label{prop:R=M} 
The moduli space $\mathcal{M}(a,b)$ can 
be identified with the moduli space of  
$\U(p,q)$-Higgs bundles which admit solutions to Hitchin's equations. 
It is a quasi-projective variety which is smooth away from the points 
representing reducible $\U(p,q)$-Higgs bundles. There is an 
homeomorphism between $\mathcal{R}_\Gamma(a,b)$ and 
$\mathcal{M}(a,b)$. This restricts to give a homeomorphism between
the subspace $\mathcal{R}^*_\Gamma(a,b)$ of  irreducible elements in 
$\mathcal{R}_\Gamma(a,b)$ and the  subspace $\mathcal{M}^s(a,b)$ of 
stable Higgs bundles in $\mathcal{M}(a,b)$. 
\end{proposition}

\begin{proof} The first statement is a direct consequence of 
Proposition \ref{prop:U(pq)HK}. The construction of 
$\mathcal{M}(a,b)$ is essentially the same as in section 
\S 9 of \cite{simpson:1994a}. There the moduli space of $G$-Higgs 
bundles is constructed for any reductive group $G$. We take 
$G=\GL(p)\times \GL(q)$. The difference between a $\U(p,q)$-Higgs bundle and a $\GL(p)\times \GL(q)$-Higgs bundle is entirely in the nature of 
the Higgs fields. Taking the standard embedding of $\GL(p)\times 
\GL(q)$ in $\GL(p+q)$ we see that in a $\GL(p)\times \GL(q)$-Higgs 
bundle the Higgs field $\Phi$ takes its values in the subspace
$(\mathfrak{gl}(p)\oplus\mathfrak{gl}(q))\subset\mathfrak{gl}(p+q)$, while 
in a $\U(p,q)$-Higgs bundle the Higgs field $\Phi$ takes its values 
in the complementary subspace 
$(\mathfrak{gl}(p)\oplus\mathfrak{gl}(q))^{\perp}$. Since both subspaces are 
invariant under the adjoint action of 
$\GL(p)\times\GL(q)$, the same method of construction works for the moduli spaces
of both types of Higgs bundle. 
 
Suppose that $(E=V\oplus W,\Phi)$ represents a point in 
$\mathcal{M}(a,b)$, i.e.\ suppose that it is a 
$\U(p,q)$-polystable Higgs bundle, and suppose that with metric $H=H_V\oplus H_W$
Hitchin's equation (\ref{eq:hitchin1}) is satisfied. Rewriting the 
equations in terms of the Higgs connection $D=d_A+\theta$, where 
$A$ is the metric connection determined by $H$ and $\theta=\Phi+\Phi^*$,
we see that $D$ is projectively flat. By Lemma 
\ref{lemma:upq} it  is a projectively flat $\U(p,q)$-connection, and 
thus defines a point in $\mathcal{R}_\Gamma(a,b)$. Conversely by 
Corlette's theorem \cite{corlette:1988}, every representation in 
$\Hom^+(\pi_1 X,\PU(p,q))$, or equivalently every representation 
in $\Hom^+(\Gamma,\U(p,q))$, arises in this way.
\end{proof}
\begin{remark}\label{remark:coprime} If $\GCD(p+q,a+b)=1$ then for 
purely numerical reasons
there are no strictly semistable $\U(p,q)$-Higgs bundles in 
$\mathcal{M}(a,b)$. In this case $\mathcal{M}^s(a,b)=\mathcal{M}(a,b)$.
\end{remark}
\begin{proposition}With $n=p+q$ and $d=a+b$, let $\mathcal{M}(d)$ denote the 
moduli space of polystable $\GL(n)$-Higgs bundles of degree $d$. If 
$p\ne q$ or $a\ne b$ then $\mathcal{M}(a,b)$ embeds as a closed subvariety in  
$\mathcal{M}(d)$. If $p=q$ and $a=b$, then there is a finite morphism from 
$\mathcal{M}(a,a)$ to $\mathcal{M}(d)$.
\end{proposition}
\begin{proof}Let $[V\oplus W,\Phi]_{p,q}$ denote the point in 
$\mathcal{M}(a,b)$ represented by the $\U(p,q)$-Higgs bundle 
$(V\oplus W,\Phi)$. Then $(E=V\oplus W,\Phi)$ is a polystable 
$\GL(n)$-Higgs bundle and the map 
$\mathcal{M}(a,b)\rightarrow\mathcal{M}(d)$ is defined by
$$[V\oplus W,\Phi]_{p,q}\mapsto [E,\Phi]_n\ ,$$
\noindent where $[,]_n$ denotes the isomorphism class in 
$\mathcal{M}(d)$. The only question is whether this map is injective.
Suppose that $(E=V\oplus W,\Phi)$ and $(E'=V'\oplus W',\Phi')$ are 
isomorphic as $\GL(n)$-Higgs bundles. Let the isomorphism be given by 
complex gauge transformation $g:E\rightarrow E'$. Since we can regard 
the smooth splitting of $E$ as fixed, we see that unless $V\cong W'$ 
and $W\cong V'$, the gauge transformation must already be of the form  
$(\begin{smallmatrix}
    g_V & 0 \\
    0 & g_W
  \end{smallmatrix})$, i.e.\ 
$[V\oplus W,\Phi]_{p,q}=[V'\oplus W',\Phi']_{p,q}$. But in order to 
have $V\cong W'$ and $W\cong V'$ we require $p=q$ and $a=b$. In that case, if $V$ and $W$ are non-isomorphic, then $[V\oplus W 
,(\begin{smallmatrix} 
    0 & \beta \\
    \gamma & 0
  \end{smallmatrix})]_n=[W\oplus V 
,(\begin{smallmatrix} 
    0 & \gamma \\
    \beta & 0
  \end{smallmatrix})]_n$ but the Higgs bundles are not isomorphic as 
  $\U(p,q)$-Higgs bundles.
\end{proof}
%%%%%%%%%%%%%%%%%%%%%%%%%%%%%%%%%%%%%%%%%%%%%%%%
\subsection{Deformation theory of Higgs bundles}
\label{sec:higgs-deformation-theory}
%%%%%%%%%%%%%%%%%%%%%%%%%%%%%%%%%%%%%%%%%%%%%%%%
A main tool in the study of the topology of moduli spaces of
Higgs bundles is given by the Morse theoretic techniques introduced by
Hitchin \cite{hitchin:1992}. In order to use these methods, we first need to 
recall  the deformation theory of Higgs bundles. We  refer to 
the paper by Biswas and Ramanan  \cite{biswas-ramanan:1994}
for details.
Let $(E,\Phi)$ be a stable  $\U(p,q)$-Higgs bundle in ${\mathcal{M}}^s(a,b)$.
Some times it will be convenient to use the following notation:
  \begin{align}\label{u-notation}
  U & = \End(E)\notag \\ 
U^+ &= \End(V) \oplus \End(W), \notag \\ 
U^- &= \Hom(W,V)  \oplus \Hom(V,W). 
 \end{align}
Clearly, $U=U^+\oplus U^-$. Note that $\Phi \in H^0(U^-\otimes K)$ 
and that 
$\ad(\Phi)$ interchanges $U^+$ and $U^-$.
As it is shown in \cite{biswas-ramanan:1994},   the Zariski tangent space to
${\mathcal{M}}^s(a,b)$ at the point defined by $(E,\Phi)$ can be identified
with  the first hypercohomology of the complex of sheaves
\begin{equation}
  \label{eq:tangentspace}
  C^{\bullet} : 
U^+  \xrightarrow{\ad(\Phi)}  U^- \otimes K.  
\end{equation}
One has the long exact sequence
\begin{equation}
  \label{eq:long-exact-tangent}
\begin{array}{ccccccccccc}
  0 &\lto  \mathbb{H}^0(C^{\bullet}) &\lto& H^0(U^+) &\lto&
  H^0(U^-\otimes K) &\lto&  \mathbb{H}^1(C^{\bullet}) & &  \\ 
    &    &  \lto & H^1(U^+) &\lto&  H^1(U^-\otimes  K) &\lto &  
\mathbb{H}^2(C^{\bullet})&\lto& 0,
\end{array}
\end{equation}
from which one  obtains the following.
\begin{proposition}\label{prop:smoothness-higgs}
The moduli space  of  stable  $\U(p,q)$-Higgs bundles
is a  smooth complex variety of dimension $1+(p+q)^2 (g-1)$.
\end{proposition}
\begin{proof}
Let $(E,\Phi)$ be a stable $\U(p,q)$-Higgs bundle.
Then $(E,\Phi)$ is simple, that is, its only automorphisms are the
non-zero scalars.  Thus, if $(E,\Phi)$ is stable,
$$
\ker\bigl(\ad(\Phi) \colon H^0(U) \to
H^0(U\otimes K) \bigr) = \C.
$$ 
Since
$U = U^+ \oplus U^-$ and
$\ad(\Phi)$ interchanges these two summands it follows that
\begin{align}
  \label{eq:ker-phi-u}
  \ker\bigl(\ad(\Phi) \colon H^0(U^+) \to
  H^0(U^-\otimes K)\bigr) &= \C \\
  \label{eq:ker-phi-u-perp}
  \ker\bigl(\ad(\Phi) \colon H^0(U^-) \to
  H^0(U^+\otimes K)\bigr) &= 0 
\end{align}
Hence, if $(E,\Phi)$ is stable, \eqref{eq:ker-phi-u} shows that
$\mathbb{H}^0(C^{\bullet}) = \C$.  
To show that the moduli space is smooth at a neighbourhood of $(E,\Phi)$
we need to show that $\mathbb{H}^2(C^{\bullet}) = 0$.
But we have natural
$\ad$-invariant isomorphisms $U^+ \cong (U^+)^*$ and
$U^- \cong (U^-)^*$.
% and, if $x \in \lie{u}_{\C}^{\perp}$ then $\ad(x)$ is symmetric 
% with respect to these isomorphisms, i.e.\ the diagram 
% \begin{displaymath}
%   \begin{CD}
%     \lie{u}_{\C} @>\cong>> \lie{u}_{\C}^* \\
%     @VV\ad(x)V              @VV\ad(x)^tV \\
%     \lie{u}_{\C}^{\perp} @>\cong>> (\lie{u}_{\C}^{\perp})^*
%   \end{CD}
% \end{displaymath}
% commutes.  
Thus
$$
\ad(\Phi) \colon H^1(U^+) \to 
H^1(U^-\otimes K)
$$
is Serre dual to $\ad(\Phi) \colon H^0(U^-) \to
H^0(U^+\otimes K)$.  Thus \eqref{eq:ker-phi-u-perp} 
shows that $\mathbb{H}^2(C^{\bullet}) = 0$.
The dimension of the moduli space is hence
\begin{eqnarray}
\dim \mathbb{H}^1(C^{\bullet}) &=& 1-\chi(U^+)+\chi(U^-\otimes K)\nonumber\\
                               &=& 1+(p^2+q^2)(g-1)+ 2pq (g-1)\nonumber\\
                               &=& 1+(p+q)^2 (g-1).\nonumber
\end{eqnarray}
\end{proof} 
\begin{remark}
Notice that the dimension of the moduli space of stable $\U(p,q)$-Higgs
bundles is half the dimension of the moduli space of stable 
$\GL(p+q,\C)$-Higgs bundles.
\end{remark} 
\begin{remark}
As pointed out previously, if $\GCD(p+q,a+b)=1$, then there are no 
strictly semistable elements in $\mathcal{M}(a,b)$ and hence 
$\mathcal{M}(a,b)$ is smooth. 
\end{remark}

%%%%%%%%%%%%%%%%%%%%%%%%%%%%%%%%%%%%%%%%%%%%%%%%%
\subsection{Bounds on the topological invariants}
\label{sec:topological-bounds}
%%%%%%%%%%%%%%%%%%%%%%%%%%%%%%%%%%%%%%%%%%%%%%%%%
In this section we shall show how the Higgs bundle point of view
provides an easy proof of a result of Domic and Toledo
\cite{domic-toledo:1987} which allows us to bound the topological
invariants $\deg(V)$ and $\deg(W)$ for which $\U(p,q)$-Higgs bundles
may exist.  The lemma is a slight variation
on the results of \cite[Section 3]{gothen:2001} (cf.\ also Lemma 3.6
of Markman and Xia \cite{markman-xia:2001}).
\begin{lemma}
  \label{lem:slope-bound}
  Let $(E,\Phi)$ be a semistable $\U(p,q)$-Higgs bundle.  Then
  \begin{align}
    \label{eq:muE1-mu}
    p(\mu(V) - \mu(E)) &\leq \rk(\gamma)(g-1), \\
    \label{eq:muE2-mu}
    q(\mu(W) - \mu(E)) &\leq \rk(\beta)(g-1).
  \end{align}
  If equality occurs in \eqref{eq:muE1-mu} then either $(E,\Phi)$ is
  strictly semistable or $p = q$ and $\gamma$ is an isomorphism. 
  If equality occurs in \eqref{eq:muE2-mu} then either $(E,\Phi)$ is
  strictly semistable or $p = q$ and $\beta$ is an isomorphism. 
\end{lemma}
\begin{proof}
  If $\gamma = 0$ then $V$ is $\Phi$-invariant and so, by stability,
  $\mu(V) \leq \mu(E)$ where equality can only occur if 
  $(E,\Phi)$ is strictly semistable.  This proves \eqref{eq:muE1-mu} 
  in the case $\gamma = 0$ and we may, therefore, assume that $\gamma 
  \neq 0$.   Let $N \subseteq V$ be the vector bundle
  associated to $\ker(\gamma)$ and let $I \subseteq W$ be the vector
  bundle associated to $\im(\gamma) \otimes K^{-1}$.  Then
  \begin{equation}
    \rk(N) + \rk(I) = p
    \label{eq:rkG+rkH}
  \end{equation}
  and, since $\gamma$ induces a non-zero section of $\det((V/N)^* \otimes
  I \otimes K)$,
  \begin{equation}
    \deg(N) + \deg(I) + \rk(I)(2g-2) \geq \deg(V).
    \label{eq:degG+degH-1}
  \end{equation}
  The bundles $N$ and $V \oplus I$ are $\Phi$-invariant subbundles 
  of $E$ and hence we obtain by semistability that $\mu(N) \leq \mu(E)$ 
  and $\mu(V \oplus I) \leq \mu(E)$ or, equivalently, that
  \begin{align}
    \label{eq:degG}
    \deg(N) &\leq \mu(E) \rk(N), \\
    \label{eq:degH}
    \deg(I) &\leq \mu(E)(p + \rk(I)) - \deg(V).
  \end{align}
  Adding \eqref{eq:degG} and \eqref{eq:degH} and using 
  \eqref{eq:rkG+rkH} we obtain
  \begin{equation}
    \deg(N) + \deg(I) \leq 2\mu(E)p - \deg(V).
    \label{eq:degG+degH-2}
  \end{equation}
  Finally, combining \eqref{eq:degG+degH-1} and 
  \eqref{eq:degG+degH-2} we get 
  \begin{displaymath}
    \deg(V) - \rk(I)(2g-2) \leq 2 \mu(E)p - \deg(V),
  \end{displaymath}
  which is equivalent to \eqref{eq:muE1-mu} since $\rk(\gamma) = \rk(I)$.  
  Note that equality can only occur if we have equality in 
  \eqref{eq:degG} and \eqref{eq:degH} and thus either $(E,\Phi)$ is 
  strictly semistable or neither of the subbundles $N$ and $V \oplus 
  I$ is proper and non-zero.  In the latter case, clearly $N = 0$ and 
  $I = W$ and therefore $p = q$; furthermore we must also 
  have equality in \eqref{eq:degG+degH-1} implying that $\gamma$ is an 
  isomorphism.
  An analogous argument applied to $\beta$ proves \eqref{eq:muE2-mu}.
\end{proof}
\begin{remark}
  The proof also shows that if we have equality in, say,
  \eqref{eq:muE1-mu} then $\gamma \colon V/N \to I\otimes K$ is an
  isomorphism.  In particular, if $p <q$ and $\mu(V) -
  \mu(E) = g-1$ then $\gamma \colon V \xrightarrow{\cong} I \otimes K$.
\end{remark}
We can reformulate Lemma~\ref{lem:slope-bound} to obtain 
the following corollary.
\begin{corollary}
  \label{cor:mu-muE_i}
  Let $(E,\Phi)$ be a semistable $\U(p,q)$-Higgs bundle.  Then
  \begin{align}
    \label{eq:mu-mu_2}
    q(\mu(E) - \mu(W)) &\leq \rk(\gamma)(g-1), \\
    \label{eq:mu-mu_1}
    p(\mu(E) - \mu(V)) &\leq \rk(\beta)(g-1).
  \end{align}
\end{corollary}
\begin{proof}
  To see that \eqref{eq:mu-mu_2} is equivalent to \eqref{eq:muE1-mu}
  one simply notes that $\mu(W) - \mu(E) =
  \frac{p}{q}\bigl( \mu(E) - \mu(V)\bigr)$.  Similarly
  \eqref{eq:mu-mu_1} is equivalent to \eqref{eq:muE2-mu}.
\end{proof} 
An important corollary of the lemma above is the 
following Milnor--Wood type inequality for $\U(p,q)$-Higgs bundles (due
to Domic and Toledo \cite{domic-toledo:1987}, improving on a bound
obtained by Dupont \cite{dupont:1978} in the case $G = \SU(p,q)$). 
This result gives bounds on the possible values of the topological
invariants $\deg(V)$ and $\deg(W)$.
\begin{corollary}
  \label{cor:toledo}
  Let $(E,\Phi)$ be a semistable $\U(p,q)$-Higgs bundle.  Then
  \begin{equation}
    \label{eq:toledo}
    \frac{pq}{p+q} \abs{\mu(V) - \mu(W)} 
    \leq \min\{p,q\}(g-1).
  \end{equation}
\end{corollary}
\begin{proof}
  Since $\mu(E) = \frac{p}{p+q}\mu(V) + \frac{q}{p+q}\mu(W)$ 
  we have $\mu(V) - \mu(E) = \frac{q}{p+q}(\mu(V) - 
  \mu(W)$ and therefore \eqref{eq:muE1-mu} gives
  \begin{displaymath}
    \frac{pq}{p+q} (\mu(V) - \mu(W)) \leq 
    \rk(\gamma)(g-1).
  \end{displaymath}
  A similar argument using \eqref{eq:muE2-mu} shows that 
  \begin{displaymath}
    \frac{pq}{p+q} (\mu(W) - \mu(V)) \leq 
    \rk(\beta)(g-1).
  \end{displaymath}
  But, obviously, $\rk(\beta)$ and $\rk(\gamma)$ are both less than or equal to 
  $\min\{p,q\}$ and the result follows.
\end{proof}

\begin{definition}\label{defn:toledo}
Let $a=\deg(V)$  and $b=\deg(W)$.  The number
  \begin{equation}
\tau=\tau(a,b)    = 2\frac{qa - pb}{p + q}
\label{toledo-invariant}
  \end{equation}
  is known as the {\bf Toledo invariant} of the representation 
  corresponding to $(E,\Phi)$.  
 \end{definition}

\begin{remark}
\label{rem:tau-M}
Since
$$
    \tau = 2\frac{pq}{p+q} (\mu(V) - \mu(W)),
$$
the inequality \eqref{eq:toledo} can 
  thus be written 
$$\abs{\tau} \leq \min\{p,q\}(2g-2)\ .$$
We denote $\tau_M=\min\{p,q\}(2g-2)$.
\end{remark}

%%%%%%%%%%%%%%%%%%%%%%%%%%%%%%%%%%%%%%%%%%%%%%%%%%%%%%%%%%%%%%%%%
\subsection{Moduli space for $p=q$ and $|\tau|=\tau_M$}
\label{subs:moduli-p=q}
%%%%%%%%%%%%%%%%%%%%%%%%%%%%%%%%%%%%%%%%%%%%%%%%%%%%%%%%%%%%%%%%%

Suppose $p=q$.  Then $\tau=\tau(p,p,a,b)=a-b$.
In this section we  give an alternative (more explicit)
description of the moduli space $\mathcal{M}(a,b)$ in the case 
in which the Toledo invariant is maximal, i.e.\ 
$|\tau|=|a-b|=\tau_M=p(2g-2)$. 

Before  doing this, we need to review briefly the notion of 
$L$-twisted Higgs pairs. Let $L$ be a line bundle. An $L$-twisted 
Higgs pair $(V,\theta)$ consists of a holomorphic vector bundle
$V$ and an $L$-twisted homomorphism $\theta:V\lto V\otimes L$.
The notions of stability, semistability and polystability 
are defined as for Higgs bundles. The moduli space of 
semistable $L$-twisted
Higgs pairs has been constructed by Nitsure using GIT 
\cite{nitsure:1991}. Let $\mathcal{M}_L(n,d)$ be the moduli
space of polystable $L$-twisted Higgs pairs of rank $n$ and
degree $d$.

\begin{proposition} \label{prop:p=q-toledo-max}
Let $p=q$ and $|a-b|=p(2g-2)$. Then
$$
\mathcal{M}(a,b)\cong \mathcal{M}_{K^2}(p,a)\cong \mathcal{M}_{K^2}(p,b).
$$
\end{proposition}

\begin{proof}
Let $(E=V\oplus W, \Phi)\in \mathcal{M}(a,b)$. Suppose for 
definiteness  that $b-a=p(2g-2)$. {}From (\ref{eq:muE2-mu}) it follows 
that $\gamma: V\lto W\otimes K$ is an isomorphism. We can then 
compose  $\beta: W\lto V\otimes K$ with 
$\gamma\otimes \Id_K: V\otimes K\lto W\otimes K^2$
to obtain a $K^2$-twisted Higgs pair
$\theta_W: W\lto W\otimes K^2$. 
Similarly, twisting $\beta: W\lto V\otimes K$ with $K$ and composing
with $\gamma$, we obtain a $K^2$-twisted Higgs pair 
$\theta_V: V\lto V\otimes K^2$. Conversely, given an isomorphism
$\gamma: V\lto W\otimes K$, we can recover $\beta$ from $\theta_V$
as well as from $\theta_W$.
It is clear that the (poly)stability of $(E,\Phi)$ is equivalent
to the (poly)stability of $(V,\theta_V)$ and to the (poly)stability of 
$(W,\theta_W)$, proving the claim.
\end{proof}

\begin{remark}
The moduli space $\mathcal{M}_{K^2}(p,a)$ contains an open (irreducible)
 subset consisting
of a rank $N$ vector bundle over  $M^s(p,a)$. This is because the
stability of $V$ implies the stability of any  $K^2$-twisted Higgs pair
$(V,\theta_V)$, and $H^1(\End V\otimes K^2)=0$. The rank $N$
is determined by the Riemann--Roch Theorem.
\end{remark}

%%%%%%%%%%%%%%%%%%%%%%%%%%%%%%%%%%%%%%%%%%%%%%%%%%%%%%%%%%%%%%%%%
\subsection{Rigidity for extreme values of the Toledo 
invariant}\label{subs:rigidity}
%%%%%%%%%%%%%%%%%%%%%%%%%%%%%%%%%%%%%%%%%%%%%%%%%%%%%%%%%%%%%%%%%

{}From the bounds in Section \ref{sec:topological-bounds}
it follows that if $p< q$ (a similar  result holds  for $p>q$)
and  $(a,b)$ such  that $|\tau|=\tau_M$ there are no stable 
$\U(p,q)$-Higgs bundles and every element in $\mathcal{M}(a,b)$ is in 
fact reducible. In particular the moduli space has smaller dimension 
than expected exhibiting a certain kind of rigidity. This phenomenon (for 
large Toledo invariant) has been studied from the point of view of 
representations of the fundamental group by D. Toledo 
\cite{toledo:1989} when  
$p=1$ and L. Hern\'andez \cite{hernandez:1991} when $p=2$. We deal here with the general case which, as far as we know,  has not 
appeared  previously in the literature. To state our result, we use  
the more precise notation   
$\mathcal{M}(p,q,a,b)$ for the moduli
space of $\U(p,q)$-Higgs bundles such that $\deg(V)=a$, and
$\deg(W)=b$, and write  the Toledo invariant as
  \begin{equation}
\tau= \tau(p,q,a,b) 
    = 2\frac{qa - pb}{p + q}.
\label{toledo-invariant-a-b}
\end{equation}

\begin{proposition}\label{prop:rigidity}
Let $(p,q,a,b)$ with  $p<q$ and $|\tau(p,q,a,b)|=p(2g-2)$.
Then every element in  $\mathcal{M}(p,q,a,b)$ is strictly semistable and
 decompose as the direct sum of a polystable $\U(p,p)$-Higgs bundle with
maximal Toledo invariant and a polystable vector  bundle of rank $(q-p)$. 
To be precise, let $\tau=p(2g-2)$, then
  \begin{equation}
    \label{rigidity}
    \mathcal{M}(p,q,a,b)\cong \mathcal{M}(p,p,a,a-p(2g-2))
\times M(q-p, b-a +p(2g-2)).
  \end{equation}
In particular, the dimension at a smooth point in 
$\mathcal{M}(p,q,a,b)$ is $2+ (p^2+5q^2-2pq)(g-1)$, and it is
hence strictly smaller than the expected dimension. (A similar result
holds if $p>q$, as well as if $\tau=-p(2g-2)$).
\end{proposition}
\begin{proof}
Suppose that $\tau(p,q,a,b)=p(2g-2)$.
Let  $(E,\Phi)\in\mathcal{M}(p,q,a,b)$. Then 
$\mu(V)-\mu(E)=g-1$ and $\mu(E)-\mu(W)=\frac{p}{q}(g-1)$. 
Since $\rk(\beta)$ and $\rk(\gamma)$ are at most $p$, it follows from
(\ref{eq:muE1-mu}) and (\ref{eq:mu-mu_1}) that 
$\rk(\beta)=\rk(\gamma)=p$. Let 
$W_\gamma=\im(\gamma)\otimes K^{-1}$ and let $W_\beta=\ker(\beta)$. 
One has that $V\oplus W_\gamma$ is a $\Phi$-invariant subbundle of 
$V\oplus W$, and $\mu(V\oplus W_\gamma)=\mu(E)$. We see that $(E,\Phi)$
is strictly semistable (as we already knew from Lemma 
\ref{lem:slope-bound}). Since it is polystable it must split as 
$$
(V\oplus W_\gamma,\Phi)\oplus (0\oplus W/W_\gamma,0).
$$
In fact, using 
\begin{displaymath}
  0 \lto \ker(\Phi) \lto V\oplus W \lto (V\oplus W)\otimes K \lto 0.
\end{displaymath}
we see that $\ker\Phi=(0\oplus W_\beta,0)$ and 
$W/W_{\gamma}=W_\beta$. It is clear that $(V\oplus W_\gamma,\Phi)\in 
\mathcal{M}(p,p,a,a-p(2g-2))$ and $W_\beta \in M(q-p, b-a +p(2g-2))$. 
Moreover, 
$(V\oplus W_\gamma,\Phi)$ has maximal Toledo invariant, that is,
$\tau(p,p,a,a-p(2g-2))=2p(g-1)$.
To complete the proof we observe that 
\begin{align*}
\dim \mathcal{M}^s(p,p,a,a-p(2g-2))+\dim  M^s(q-p, b-a +p(2g-2))\\
=1+(2p)^2(g-1)+1+(q-p)^2(g-1)=2+(p^2+5q^2-2pq)(g-1),
\end{align*}
which, since $q>1$,  is smaller than $1+(p+q)^2(g-1)$, the dimension of
$\mathcal{M}(p,q,a,b)$.
\end{proof}

\begin{corollary}
Let $(p,q,a,b)$ with  $p<q$ and $\tau(p,q,a,b)=p(2g-2)$.
Then 
$$
    \mathcal{M}(p,q,a,b)\cong \mathcal{M}_{K^2}(p,a-p(2g-2))
\times M(q-p, b-a +p(2g-2)).
$$
\end{corollary}
\begin{proof}
It follows from Propositions \ref{prop:rigidity} and \ref{prop:p=q-toledo-max}.
\end{proof}

%%%%%%%%%%%%%%%%%%%%%%%%%%%%
\section{Morse theory}
\label{sec:morse-theory}
%%%%%%%%%%%%%%%%%%%%%%%%%%%%

%%%%%%%%%%%%%%%%%%%%%%%%%%%%%%%
\subsection{The Morse function}\label{subs:Morse}
%%%%%%%%%%%%%%%%%%%%%%%%%%%%%%%
An extremely efficient tool for studying topological properties of
Higgs bundle moduli is considering the $\C^*$-action on the moduli
space given by multiplying the Higgs field $\Phi$ by a non-zero
scalar.  In order to explain this, it is convenient to consider
the moduli space from the gauge theory point of view.  
The $\U(1)$-action $\Phi \mapsto e^{i\theta} \Phi$ on
$\mathcal{M}(a,b)$ preserves solutions to Hitchin's equations
\eqref{eq:hitchin1} (but the full $\C^*$-action does not preserve
solutions to the first equation).  Restricted to the smooth locus of
$\mathcal{M}(a,b)$ this turns out to be a Hamiltonian circle
action with moment map (up to multiplication by a constant)
\begin{displaymath}
  f(A,\Phi) = \int_{X}\abs{\Phi}^2,
\end{displaymath}
the $L^2$-norm squared of the Higgs field.  Thus, if the moduli space
is smooth, $f$ is a perfect Bott-Morse function (by a theorem of
Frankel \cite{frankel:1959}; the essential point is that the Morse
indices are even) and the critical points of $f$ are exactly the fixed
points of the circle action.
In particular, if ${\mathcal{M}}(a,b)$ is smooth, then its number of 
connected components is the same as the number of connected 
components of the subspace of local minima of $f$.
However, even if ${\mathcal{M}}(a,b)$ is not smooth, $f$ can be
used to obtain information about the connected components of
${\mathcal{M}}(a,b)$ in the following way.  It was shown by Hitchin
\cite{hitchin:1987}, using Uhlenbeck's weak compactness theorem, that
$f$ is proper. We have the following general result.
\begin{proposition}
\label{prop:real-topology-exercise}
Let $Z$ be a Hausforff space and let $f \colon Z \to \R$ be proper and
bounded below.  Then $f$ attains a minimum on each connected
component of $Z$ and, furthermore, if the subspace of local minima of
$f$ is connected then so is $Z$. \hfill\qed
\end{proposition}

In particular this applies to our situation.

\begin{proposition}
\label{prop:topology-exercise}
Let $f \colon \mathcal{M}(a,b) \to \R$ be defined as above.  Then $f$
has a minimum on each connected component of ${\mathcal{M}}(a,b)$.
Moreover, if the subspace of local minima of $f$ is connected then so
is $\mathcal{M}(a,b)$.  \hfill\qed
\end{proposition}
%%%%%%%%%%%%%%%%%%%%%%%%%%%%%%%%%%%%%%%%%%%%%%%%%%%
\subsection{Critical points of the Morse function}
\label{sec:critical-points}
%%%%%%%%%%%%%%%%%%%%%%%%%%%%%%%%%%%%%%%%%%%%%%%%%%%
Next we recall Hitchin's method \cite{hitchin:1987,hitchin:1992} for
determining the local minima of $f$.  A point $(E,\Phi)$ is a fixed
point of the circle action if and only if it is a system of Hodge
bundles, that is, 
\begin{equation}
  E = F_1 \oplus \cdots \oplus F_m
  \label{eq:variation-of-hodge}
\end{equation}
for holomorphic
vector bundles $F_i$ such that the restriction 
$$
\Phi_i := \Phi_{|F_i} \in H^0(\Hom(F_i,F_{i+1})\otimes K).
$$
To see why this is true, note that, if $(A,\Phi)$ represents a fixed
point then there must be a 1-parameter family of gauge transformations
$g(\theta)$ taking $(A,\Phi)$ to $(A,e^{i\theta}\Phi)$ and this gives
an infinitesimal $\U(p)\times \U(q)$-gauge transformation $\psi =
\dot{g}$ which is covariantly constant (i.e.\ $d_{A}\psi = 0$) and
such that $[\psi,\Phi] = i\Phi$.  It follows that we can decompose $E$
in holomorphic subbundles $F_{\lambda}$ on which $\psi$ acts as
$i\lambda$ and furthermore that $\Phi$ maps $F_{\lambda}$ to
$F_{\lambda + 1} \otimes K$.  If $(E,\Phi)$ is stable each of these
components of $\Phi$ is non-zero since, otherwise, $(E,\Phi)$ would be
reducible.  We can therefore write $E = F_{1} \oplus \cdots \oplus
F_{m}$ where the weight of $\psi$ on $F_{k+1}$ is one plus the weight
of $\psi$ on $F_{k}$.  In general (for reducible Higgs bundles)
$(E,\Phi)$ is the direct sum of such chains.  Furthermore, note that
$\psi$ preserves $V$ and $W$ and, therefore, the $F_k$ are direct sums
of bundles contained in $V$ and $W$.
The decomposition $E = F_1 \oplus \cdots \oplus F_m$ gives a
corresponding decomposition of the bundle $U=\End(E)$ into eigenbundles for the adjoint action of
$\psi$:
\begin{displaymath} 
U  = \bigoplus_{k=-m+1}^{m-1} U_k,
\end{displaymath}
where
\begin{math}
  U_k = \bigoplus_{i-j = k} \Hom(F_j,F_i).
\end{math}
We now show how to find the local minima of $f$ which are represented 
by stable Higgs bundles.  In fact, one can do much better: it is 
possible to calculate the Morse index at the critical point.
It follows from Hitchin's calculations in \cite[\S 8]{hitchin:1992} 
(see also \cite[Section 2.3.2]{gothen:1995})
that the subspace of the tangent space at $(E,\Phi)$ on which the
Hessian of $f$ has eigenvalue $-k$ is $\mathbb{H}^1$ of the following
complex:
\begin{equation}
  \label{eq:weight-k-space}
  C^{\bullet}_{k} : 
  U_{k}^{+}
  \xrightarrow{\ad(\Phi)}
  U_{k+1}^{-} \otimes K,
\end{equation}
where we use the notation
\begin{align*}
  U_k^{+} &= U_{k} \cap U^+, \\
  U_k^{-} &= U_{k} \cap U^-,
\end{align*}
with $U^+$ and $U^-$ as defined  in (\ref{u-notation}).
This description of the eigenspace of the Hessian of $f$ gives rise to
the long exact sequence
\begin{equation}
  \label{eq:long-exact-k}
\begin{array}{ccccccccccc}
  0 &\lto  \mathbb{H}^0(C^{\bullet}_k) &\lto& H^0(U^+_k) &\lto&
  H^0(U^-_{k+1}\otimes K) &\lto&  \mathbb{H}^1(C^{\bullet}_k) & &  \\ 
    &    &  \lto & H^1(U^+_k) &\lto&  H^1(U^-_{k+1}\otimes  K) &\lto &  
\mathbb{H}^2(C^{\bullet}_k)&\lto& 0.
\end{array}
\end{equation}
Note that if some $F_{k}$ has a non-zero component
in both $V$ and $W$ this leads a decomposition of $(E,\Phi)$ and so,
if $(E,\Phi)$ is stable, the $F_k$ are alternately contained in $V$
and $W$.  It follows that
\begin{equation}
      \label{eq:u-even-odd}
U^+= \bigoplus_{\text{$k$ even}} U_k; \quad
U^- = \bigoplus_{\text{$k$ odd}} U_k.
\end{equation}
Thus we see that all the
eigenvalues of the Hessian of $f$ are even and that
\begin{equation}
  \label{eq:weight-2k-space}
  C^{\bullet}_{2k} : 
  U_{2k} \xrightarrow{\ad(\Phi)} U_{2k+1} \otimes K.
\end{equation}
Hence we have the long exact sequence
\begin{multline}
  \label{eq:long-exact-2k}
  0 
  \to \mathbb{H}^0(C^{\bullet}_{2k}) \to H^0(U_{2k}) 
  \to H^0(U_{2k+1} \otimes K) \to \mathbb{H}^1(C^{\bullet}_{2k}) \\ 
  \to H^1(U_{2k}) \to H^1(U_{2k+1} \otimes K) \to 
  \mathbb{H}^2(C^{\bullet}_{2k}) \to 0.  
\end{multline}
Note that when $(E,\Phi)$ is a stable critical point we have
$\mathbb{H}^0(C^{\bullet}_{2k}) = \mathbb{H}^2(C^{\bullet}_{2k}) = 0$
for $k \neq 0$, while $\mathbb{H}^2(C^{\bullet}_{0}) = 0$ and
$\mathbb{H}^0(C^{\bullet}_{0}) = \C$.  We can therefore use the
exact sequence \eqref{eq:long-exact-2k} and the
Riemann--Roch formula to calculate the dimension of
$\mathbb{H}^1(C^{\bullet}_{2k})$:
\begin{multline}
  \label{eq:dim(H^1(C_2k))}
  \dim(\mathbb{H}^1(C^{\bullet}_{2k})) = \chi(U_{2k+1} \otimes K)
  - \chi(U_{2k}) \\
  = (g-1)\bigl(\rk(U_{2k+1}) + \rk(U_{2k})\bigr) + \deg(U_{2k+1}) -
  \deg(U_{2k})
\end{multline}
for $k \neq 0$, while for $k=0$ we obtain
\begin{equation}
  \label{eq:dim(H^1(C_0))}
  \dim(\mathbb{H}^1(C^{\bullet}_{0})) = 1 + (g-1)\bigl(\rk(U_{1}) +
  \rk(U_{0})\bigr) + \deg(U_{1}) - \deg(U_{0}).
\end{equation}
In particular we can calculate the Morse index of $f$ at a smooth
critical point:
\begin{equation}
  \label{eq:morse-index}
\dim\biggl(\mathbb{H}^1\Bigl(\bigoplus_{k \geq 1}
  C^{\bullet}_{2k}\Bigl)\biggr)= 
\sum_{k = 2}^{m-1} \bigl((g-1) \rk (U_k) + (-1)^{k+1} \deg(U_k)\bigr).
\end{equation}
Note that this is the complex dimension; thus the Morse index is twice
this number.
{}From this discussion we see that the following proposition holds.
\begin{proposition}
  A stable $\U(p,q)$-Higgs bundle of the form
  \eqref{eq:variation-of-hodge} corresponds to a local minimum of $f$
  if and only if
\begin{displaymath}
\mathbb{H}^1( C^{\bullet}_{2k}) = 0
\end{displaymath}
for all ${k \geq 1}$.
\end{proposition}

\begin{definition}
Let 
\begin{equation}
 \mathcal{N}(a,b) =\{ (E,\Phi) \in \mathcal{M}(a,b)\;\;  |
\;\; \beta= 0 \;\;\mbox{or}\;\;\gamma=0\}. \label{minima}
\end{equation}
\end{definition}

The following is an important  result for our approach that 
 characterizes the local minima of $f$.
\begin{theorem}
  \label{thm:minima}
Let $(E,\Phi)$ be  a polystable $\U(p,q)$-Higgs bundle 
in ${\mathcal{M}}(a,b)$. Then $(E,\Phi)$ is a local 
minimum of $f \colon  {\mathcal{M}}(a,b) \to \R$ if and only 
if $(E,\Phi)$ belongs to  $\mathcal{N}(a,b)$.
\end{theorem}
\begin{proof}
This follows directly from Propositions 
\ref{prop:thm-minima-1},
\ref{prop:thm-minima-2},
\ref{lem:thm-minima-3}, and
\ref{lem:thm-minima-4},
 which are given in the following three sections. 
\end{proof}
\begin{remark}
  \label{rem:minima-for-lower-ranks}
  This Theorem was already known to hold when $p,q \leq 2$ (by
  the results of \cite{gothen:2001}, Hitchin \cite{hitchin:1987},
  and Xia \cite{xia:2000}), and also when $p=q$ and $(p-1)(2g-2)<|\tau|\leq p(2g-2)$
by Markman-Xia \cite{markman-xia:2001}.
\end{remark}

Which section actually vanishes for a minimum is given by the following.

\begin{proposition} 
\label{prop:vanishing}
Let $(E,\Phi)\in \mathcal{N}(a,b)$. Then 
  \begin{enumerate}
    \item[(1)] $\gamma = 0\Leftrightarrow a/p \leq b/q
    \Leftrightarrow\tau \leq 0$,
    \item[(2)] $\beta =0\Leftrightarrow a/p \geq b/q
  \Leftrightarrow\tau \geq 0$,
\end{enumerate}
In particular, $\beta = \gamma = 0$ if and only if $a/p = b/q$ (i.e.
$\tau =0$).
  
\end{proposition}
\begin{proof}
The relation between the conditions on $\tau$ and the conditions on
$a/p-b/q$ follows directly from the definition of 
$\tau$ (cf.\ 
\ref{toledo-invariant}). By Lemma~\ref{lem:slope-bound} and  we have 
\begin{align*}
\beta =&0\Rightarrow a/p \geq b/q\\
\gamma =&0\Rightarrow a/p \leq b/q.
\end{align*}
Since we assume that at least one of $\beta$ or $\gamma$ must vanish, 
we can use the contrapositives of these implications to infer further 
that 
\begin{align*} 
a/p > b/q &\Rightarrow \gamma=0\\ a/p < b/q & \Rightarrow \beta=0\\ 
\end{align*} 
It remains to check that $\beta=\gamma=0$ if $a/p = b/q$ and 
$(E,\Phi)$ is polystable. We know that at least one of $\beta$ and $\gamma$
must vanish. Suppose that $\beta=0$. Then $\Phi(W)=0\subset W\otimes 
K$, i.e.\ $W$ is $\Phi$-invariant. By stability, we get that 
$\mu(W)\le \mu(E)$, with equality if and only if $(E,\Phi)$ splits as 
$(V,\Phi)\oplus (W,0)$. But $\mu(W)=\mu(E)$ since $a/p = b/q$. Thus  
$(E,\Phi)$ splits as indicated, and hence $V$ must be 
$\Phi$-invariant. Since $\Phi(V)=\gamma(V)\subset W\otimes K$, it follows that 
$\gamma=0$. A similar argument shows that $\beta=0$ if $\gamma=0$.

This can also conveniently be seen from the gauge theory point of 
view: Hitchin's equations show that the value of $f$ on a 
$\U(p,q)$-Higgs bundle such that $\beta =0$ is $f = \norm{\gamma}^{2} 
= q(\mu(W) - 
\mu(E))$, while $f = \norm{\beta}^{2} = p(\mu_{V}- \mu(E))$ if 
$\gamma = 0$. 
\end{proof}
\begin{corollary}\label{minima-toledo-0}
If $a/p = b/q$ then $\mathcal{N}(a,b)\cong M(p,a)\times M(q,b)$.
\end{corollary}

\begin{proof} If  $a/p = b/q$, then any $(E,\Phi)\in \mathcal{N}(a,b)$
has $E=V\oplus W$ and $\Phi=0$. The polystability of $(E,\Phi)$ is 
thus equivalent to the polystability of $V$ and $W$. 
\end{proof}

%%%%%%%%%%%%%%%%%%%%%%%%%%%%%%%%%%%%%%%%%%%%%%%%%%%
\subsection{Local minima and the adjoint bundle}
\label{sec:local-minima-adjoint}
%%%%%%%%%%%%%%%%%%%%%%%%%%%%%%%%%%%%%%%%%%%%%%%%%%%
In this section we find a criterion for $(E,\Phi)$ to be a local
minimum in terms of the adjoint bundle.  
%This may also be useful for groups other than $\U(p,q)$.  
We use the notation introduced in Section \ref{sec:critical-points}. 
Consider the complex $C^{\bullet}_{k}$ defined in
\eqref{eq:weight-k-space} and let
\begin{displaymath}
  \chi(C^{\bullet}_{k}) = \dim \mathbb{H}^{0}(C^{\bullet}_{k})
  - \dim \mathbb{H}^{1}(C^{\bullet}_{k})
  + \dim \mathbb{H}^{2}(C^{\bullet}_{k}).
\end{displaymath}
\begin{proposition}
  \label{prop:adjoint-minima}
  Let $(E,\Phi)$ be a stable $\U(p,q)$-Higgs bundle which is a fixed
  point of the $S^1$-action on $\mathcal{M}(a,b)$.
  Then $\chi(C^{\bullet}_k) \leq 0$ and equality holds if and only if 
  \begin{displaymath}
    \ad(\Phi) \colon U_{k}^{+} \to U_{k+1}^{-} \otimes K
  \end{displaymath}
  is an isomorphism.
\end{proposition}
\begin{proof}
  For simplicity we shall adopt the notation
  \begin{displaymath}
    \Phi_{k}^{\pm}
     = \Phi_{| U_{k}^{\pm}} \colon U_{k}^{\pm} \lto U_{k+1}^{\mp}
     \otimes K.
  \end{displaymath}
  The key fact we need is that there is a natural $\ad$-invariant
  isomorphism $U \cong
  U^{*}$ under which we have $U^+ \cong
  (U^+)^*$, $U^- \cong (U^-)^*$ and $U_{k}^{\pm} \cong (U_{-k}^{\pm})^{*}$. 
  Since $\ad(\Phi)^{t} = \ad(\Phi) \otimes 1_{K^{-1}}$ under
  this isomorphism we have
  \begin{equation}
    \label{eq:Phi-symmetric}
    (\Phi_{k}^{\pm})^{t} = \Phi_{-k-1}^{\mp} \otimes 1_{K^{-1}}.
  \end{equation}
  We have the short exact sequence
  \begin{displaymath}
    0 \lto \ker(\Phi_{k}^{+}) \lto (U_{k+1}^{-} \otimes K)^* \lto
    \im(\Phi_{k}^{+}) \lto 0.
  \end{displaymath}
 {}From \eqref{eq:Phi-symmetric} we have $\ker(\Phi_{k}^{+,t}) \cong
  \ker(\Phi_{-k-1}^{-}) \otimes K^{-1}$. Thus, tensoring the above
  sequence by $K$, we obtain the short exact sequence
  \begin{displaymath}
    0 \lto \ker(\Phi_{-k-1}^{-}) \lto (U_{k+1}^{-})^* \lto
    \im(\Phi_{k}^{+}) \otimes K \lto 0.
  \end{displaymath}
  It follows that 
  \begin{displaymath}
    \deg(\im(\Phi_k^{+})) \leq \deg(U_{k+1}^{-}) +
    (2g-2)\rk(\Phi_{k}^{+}) + \deg(\ker(\Phi_{-k-1}^{-})).
  \end{displaymath}
  Combining this inequality with the fact that
  \begin{equation}
    \label {eq:degU_k^+}
      \deg(U_{k}^{+}) \leq
      \deg(\ker(\Phi_{k}^{+})) + \deg(\im(\Phi_{k}^{+})), 
  \end{equation}
  we obtain
  \begin{equation}
    \label{eq:deg-im-phi_k}
    \deg(U_k^{+}) \leq \deg(U_{k+1}^{-}) +
    (2g-2)\rk(\Phi_{k}^{+}) + \deg(\ker(\Phi_{-k-1}^{-}))
    + \deg(\ker(\Phi_{k}^{+})). 
  \end{equation}
  Since $(E,\Phi)$ is semistable, so is the Higgs bundle
  $(\End(E),\ad(\Phi))$.  Clearly $\ker(\Phi_{k}^{\pm}) \subseteq
  \End(E)$ is $\Phi$-invariant and hence, from semistability,
  \begin{displaymath}
    \deg(\ker(\Phi_{k}^{\pm})) \leq 0,
  \end{displaymath}
  for all $k$.  Substituting this inequality in
  \eqref{eq:deg-im-phi_k}, we obtain
  \begin{equation}
    \label{eq:deg-U_k^+}
    \deg(U_k^{+}) \leq \deg(U_{k+1}^{-}) +
    (2g-2)\rk(\Phi_{k}^{+}).
  \end{equation}
{}From the long exact sequence \eqref{eq:weight-k-space} and the
  Riemann--Roch formula we obtain
  \begin{align*}
    \chi(C^{\bullet}_k)
    &= \chi(U_{k}^{+}) - \chi(U_{k+1}^{-} \otimes K) \\
    &= (1-g)\bigl(\rk(U_{k}^{+}) + \rk(U_{k+1}^{-}) \bigr)
       +\deg(U_{k}^{+}) - \deg(U_{k+1}^{-}).
  \end{align*}
  Using this identity and the inequality \eqref{eq:deg-U_k^+} we see
  that 
  \begin{displaymath}
    \chi(C^{\bullet}_k) \leq
    (g-1)\bigl(2 \rk(\Phi_{k}^{+}) - \rk(U_{k}^{+}) - \rk(U_{k+1}^{-})
    \bigr). 
  \end{displaymath}
  Hence $\chi(C^{\bullet}_k) \leq 0$.  Furthermore, if equality holds
  we have 
  \begin{displaymath}
    \rk(\Phi_{k}^{+}) = \rk(U_{k}^{+}) = \rk(U_{k+1}^{-})
  \end{displaymath}
  and also equality must hold in \eqref{eq:deg-U_k^+} and so
  $\deg(\im(\Phi_{k}^{+})) = \deg(U_{k+1}^{-} \otimes K)$, showing
  that $\Phi_{k}^{+}$ is an isomorphism as claimed.
\end{proof}
%%%%%%%%%%%%%%%%%%%%%%%%%%%%%%%%%%%%%%%%%%%%%%%%%%%%%%%%%%%%%%%%%%%%%%
\subsection{Stable Higgs bundles}
\label{sec:stable-minima}
%%%%%%%%%%%%%%%%%%%%%%%%%%%%%%%%%%%%%%%%%%%%%%%%%%%%%%%%%%%%%%%%%%%%%%
In this section we shall prove Theorem~\ref{thm:minima} for stable 
Higgs bundles. The reducible (polystable) ones will be dealt with 
in the next section.  We continue to use the notation of Section
\ref{sec:critical-points}.
\begin{proposition}
  \label{prop:thm-minima-1}
  Let $(E,\Phi)$ be a stable $\U(p,q)$-Higgs bundle with $\beta =0$ 
  or $\gamma=0$.  Then $(E,\Phi)$ is a local minimum of $f$.
\end{proposition}
\begin{proof}
  If $\beta =\gamma=0$ then clearly $(E,\Phi)$ is an absolute minimum of $f$.
  Otherwise such a Higgs bundle is a Hodge bundle of
  length $2$, that is, $E = F_{1} \oplus F_{2}$ with $F_{1}= V$ and
  $F_{2} = W$ (if $\beta =0$) or vice-versa (if $\gamma=0$).  Hence $\End(E) =
  U_{-1} \oplus U_{0} \oplus U_{1}$, in other words, $U_{k} = 0$ for
  $\abs{k} > 1$.  It follows that the complex $C^{\bullet}_{k}$ is
  zero for any $k>0$ and hence all eigenvalues of the Hessian of $f$
  are positive.
\end{proof}
The hard part is to show that any other critical point is not a local 
minimum of $f$.
\begin{proposition}
  \label{prop:thm-minima-2}
  Let $(E,\Phi) = (F_{1} \oplus \cdots \oplus F_{m}, \Phi)$ be a stable
  $\U(p,q)$-Higgs bundle representing a critical point of $f$ such that 
  $m \geq 3$.  Then $(E,\Phi)$ is not a local minimum of $f$.
\end{proposition}
\begin{proof}
  Note that $U_{k} = 0$ for $\abs{k} \geq m$; in particular $U_{m} =
  0$.  Note also that, since $(E,\Phi)$ is stable,
  $\mathbb{H}^0(C^{\bullet}_{m-1}) = \mathbb{H}^2(C^{\bullet}_{m-1}) =
  0$ (cf.\ the discussion at the end of Section
  \ref{sec:higgs-deformation-theory}) and therefore
  \begin{equation}
    \label{eq:hyper1=-chi}
    \mathbb{H}^1(C^{\bullet}_{m-1}) = - \chi(C^{\bullet}_{m-1}).
  \end{equation}
  We shall consider the cases when
  $m$ is odd and even separately.
  \emph{The case $m$ odd.}
  In this case $m-1$ is even and therefore \eqref{eq:u-even-odd} shows
  that $U_{m-1}^{+} = U_{m-1} \neq 0$ while $U_{m}^{-} \subseteq U_m =
  0$.  It therefore follows from Proposition \ref{prop:adjoint-minima}
  that $\chi(C^{\bullet}_{m-1}) < 0$.  Hence we have
  \begin{displaymath}
    \mathbb{H}^1(C^{\bullet}_{m-1}) = - \chi(C^{\bullet}_{m-1}) > 0,
  \end{displaymath}
  showing that $(E,\Phi)$ is not a local minimum of $f$.
  
  \emph{The case $m$ even.}  In this case we shall show that
  $\mathbb{H}^1(C^{\bullet}_{m-2}) \neq 0$.  First note that
  \eqref{eq:u-even-odd} shows that $U_{m-2}^{+} = U_{m-2}$ and
  $U_{m-1}^{-} = U_{m-1}$.  Thus, using \eqref{eq:hyper1=-chi} and
  Proposition \ref{prop:adjoint-minima}, we get that
  $\mathbb{H}^1(C^{\bullet}_{m-2})$ vanishes if and only if
  \begin{displaymath}
    \ad(\Phi) \colon U_{m-2} \to U_{m-1} \otimes K
  \end{displaymath}
  is an isomorphism.  Note that
  \begin{align*}
    U_{m-1} &= \Hom(F_{1},F_{m}) \\
    U_{m-2} &= \Hom(F_{1},F_{m-1}) \oplus \Hom(F_{2},F_{m}).
  \end{align*}
  If $\ad(\Phi) \colon U_{m-2} \to U_{m-1} \otimes K$ is an
  isomorphism, the same is true for its restriction to any fibre.
  But, by stability of $(E,\Phi)$, the bundles $F_1$, $F_2$, $F_{m-1}$
  and $F_m$ are all non-zero and, therefore, Lemma
  \ref{lem:linear-algebra} below shows that $\ad(\Phi) \colon U_{m-2}
  \to U_{m-1} \otimes K$ cannot be an isomorphism.
\end{proof}
\begin{lemma}
  \label{lem:linear-algebra}
  Let $c' \colon \mathbb{V}'_2 \to \mathbb{V}'_1$ and $c'' \colon \mathbb{V}''_2 \to \mathbb{V}''_1$ be
  linear maps between finite dimensional vector spaces.
  Assume that $\mathbb{V}'_1 \oplus \mathbb{V}'_2 \neq 0$ and
  $\mathbb{V}''_1 \oplus \mathbb{V}''_2 \neq 0$.
  Define
  \begin{align*}
    C \colon \Hom(\mathbb{V}''_1,\mathbb{V}'_1) \oplus \Hom(\mathbb{V}''_2,\mathbb{V}'_2) &\lto
    \Hom(\mathbb{V}''_2,\mathbb{V}'_1) \\
    (\psi_1,\psi_2) &\longmapsto c'\psi_2 - \psi_1c''.
  \end{align*}
  If $C$ is an isomorphism, then exactly one of the following
  alternatives must occur:
  \begin{itemize}
  \item[$(1)$] $\mathbb{V}'_1 = \mathbb{V}''_2 =0$ and $c'=c''=0$. 
  \item[$(2)$] $\mathbb{V}''_1=0$, $\mathbb{V}'_1, \mathbb{V}'_2,
  \mathbb{V}''_2 \neq 0$ and $c'\colon \mathbb{V}'_2
  \overset{\cong}{\longrightarrow} \mathbb{V}'_1$.  \item[$(3)$]
  $\mathbb{V}'_2=0$, $\mathbb{V}'_1, \mathbb{V}''_1, \mathbb{V}''_2
  \neq 0$ and $c''\colon \mathbb{V}''_2
  \overset{\cong}{\longrightarrow} \mathbb{V}''_1$.
  \end{itemize}
  In particular, if $\mathbb{V}'_1$, $\mathbb{V}'_2$, $\mathbb{V}''_1$
  and $\mathbb{V}''_2$ are all non-zero then $C$ cannot be an
  isomorphism.
\end{lemma}
\begin{proof}
  If $(c',c'') = (0,0)$ then $C = 0$ and therefore 
  \begin{displaymath}
    \Hom(\mathbb{V}''_2,\mathbb{V}'_1) 
    = \Hom(\mathbb{V}''_1,\mathbb{V}'_1) 
    = \Hom(\mathbb{V}''_2,\mathbb{V}'_2) = 0.
  \end{displaymath}
  If $\mathbb{V}'_1 \neq 0$ then $\mathbb{V}''_1 = \mathbb{V}''_2 = 0$
  which is absurd, hence $\mathbb{V}'_1 = 0$.  Similarly one sees that
  $\mathbb{V}''_2 \neq 0$ and thus alternative $(1)$ occurs.
  Henceforth assume that $(c',c'') \neq (0,0)$.  Let $r'_i = \dim
  \mathbb{V}'_i$ and $r''_i = \dim \mathbb{V}''_i$ for $i=1,2$.  If
  $C$ is an isomorphism then $r''_1 r'_1 + r''_2 r'_2 = r''_2 r'_1$
  from which it follows that
  \begin{align*}
    r''_2(r'_1 - r'_2) &= r''_1 r'_1, \\
    r'_1(r''_2-r''_1)  &= r''_2 r'_2.
  \end{align*}
  Hence
  \begin{align}
    r'_1 &\geq r'_2, \label{eq:r'_1-geq-r'_2} \\
    r''_2 &\geq r''_1. \label{eq:r''_2-geq-r''_1}
  \end{align}
  Assume that we have strict inequality in \eqref{eq:r'_1-geq-r'_2}
  and \eqref{eq:r''_2-geq-r''_1}.  Then, in particular, $\im(c')$
  and $\ker(c'')$ must both be non-zero.  Choose a complement to
  $\im(c')$ in $\mathbb{V}'_1$ so that 
  \begin{displaymath}
    \mathbb{V}'_1 = \im(c') \oplus \im(c')^{\perp}.
  \end{displaymath}
  We then have an inclusion
  \begin{displaymath}
    \Hom(\ker(c''),\im(c')^{\perp}) \into \Hom(\mathbb{V}''_2,\mathbb{V}'_1).
  \end{displaymath}
  Let $\psi = (\psi_1,\psi_2) \in \Hom(\mathbb{V}''_1,\mathbb{V}'_1) \oplus
  \Hom(\mathbb{V}''_2,\mathbb{V}'_2)$ and $x \in \ker(c'')$, then
  \begin{displaymath}
    C(\psi)(x) = c' \psi_2 (x) - \psi_1 c'(x) = c' \psi_2 (x),
  \end{displaymath}
  which belongs to $\im(c')$.  Hence $\im(C)$ and
  $\Hom(\ker(c''),\im(c')^{\perp})$ have trivial intersection and,
  therefore, $C$ cannot be an isomorphism, which is absurd.  It
  follows that equality must hold in at least one of the inequalities
  \eqref{eq:r'_1-geq-r'_2} and \eqref{eq:r''_2-geq-r''_1}.
  Suppose that equality holds in \eqref{eq:r'_1-geq-r'_2} then $r''_1
  r'_1 = 0$.  Suppose first that $r'_1 =0$, i.e.\ $\mathbb{V}'_1 = 0$,
  then $\mathbb{V}'_2 \neq 0$ and, since $c'=0$, we also have $c''
  \neq 0$.  It follows that $\mathbb{V}'_2 \neq 0$, $\mathbb{V}''_2
  \neq 0$ and that $C(\psi_1,\psi_2) = -\psi_1 c''$.  Since $C$ is an
  isomorphism this shows that alternative $(2)$ occurs.
  In a similar manner one sees that if equality holds in
  \eqref{eq:r''_2-geq-r''_1} then alternative $(3)$ occurs.  Obviously
  the three alternatives are mutually exclusive.
\end{proof}
%%%%%%%%%%%%%%%%%%%%%%%%%%%%%%%%%%%%%%%%%%%%%%%%%%%%%%%%%%%%%%%%%%%%%%%%%%%
\subsection{Reducible Higgs bundles}
\label{sec:reducible-minima}
%%%%%%%%%%%%%%%%%%%%%%%%%%%%%%%%%%%%%%%%%%%%%%%%%%%%%%%%%%%%%%%%%%%%%%%%%%%
 
In this section we shall finally conclude the proof of 
Theorem~\ref{thm:minima} by showing that it also holds for reducible 
Higgs bundles.  First we shall show that a reducible Higgs bundle 
which is not of the form given in Theorem~\ref{thm:minima} cannot be 
a local minimum of $f$; for this we use an argument similar to the one
given by Hitchin \cite[\S 8]{hitchin:1992} for the case of
$\mathrm{PSL}(n,\R)$. 
\begin{proposition}
  \label{lem:thm-minima-3}
  Let $(E,\Phi)$ be a reducible $\U(p,q)$-Higgs bundle.  
  If $\beta \neq 0$ and $\gamma \neq 0$ then $(E,\Phi)$ is not a local minimum 
  of $f$.
\end{proposition}
\begin{proof}
  We know that $(E,\Phi)$ is a direct sum of stable $\U(p',q')$-Higgs
  bundles of lower rank.  Since $f(E,\Phi)$ is the sum of the values
  of $f$ on each of the summands (on the corresponding lower rank
  moduli space), it follows that each of these stable summands is a
  local minimum of $f$ on the lower rank moduli space.  In particular,
  $(E,\Phi)$ is a fixed point of the circle action.  Therefore, by
  Proposition \ref{prop:thm-minima-1}, each stable direct summand has 
$\beta=0$
  or $\gamma=0$.  Hence we can choose two stable direct summands
  $(E',\Phi')$ and $(E'',\Phi'')$ such that $\gamma' \neq 0$ and 
$\beta'' \neq
  0$ and $\beta'= \gamma''=0$.  It is clearly sufficient to show that $(E'
  \oplus E'',\Phi' \oplus \Phi'')$ is not a local minimum of $f$ on
  the corresponding moduli space and we can therefore assume that
  $(E,\Phi) = (E' \oplus E'',\Phi' \oplus \Phi'')$ without loss of
  generality.
  
  Let $\mathcal{C}$ denote the configuration space of all solutions to
  Hitchin's equations and let $(A_0,\Phi_0) \in \mathcal{C}$ be the
  gauge theory object corresponding to $(E,\Phi)$.  As in
  Section~\ref{sec:critical-points} we denote by $\psi = \dot{g}$ the
  infinitesimal gauge transformation counteracting the circle action.
  We can write $\psi = \psi' + \psi''$, where $\psi'$ and $\psi''$ are
  infinitesimal gauge transformations of $(E',\Phi')$ and
  $(E'',\Phi'')$ respectively; recall that we may assume that
  $\tr(\psi') = \tr{\psi}'' = 0$.  The calculations of \cite[\S
  8]{hitchin:1992} show that if $\alpha(t) = (A_{t},\Phi_{t})$ is a
  smooth curve in $\mathcal{C}$ such that $(A_{0},\Phi_{0})$
  represents a fixed point of the circle action then $(f \circ
  \alpha)'(0) = 0$.  Furthermore, if $\alpha'(0) =
  (\dot{A},\dot{\Phi})$ is a tangent vector to $\mathcal{C}$ at
  $(A_{0},\Phi_{0})$ and $\psi$ acts on $\dot{A}$ with weight
  $\lambda$ and on $\dot{\Phi}$ with weight $\nu$ then
  \begin{displaymath}
    (f \circ \alpha)''(0) = -\lambda c_{1} - (\nu - 1) c_{2}
  \end{displaymath}
  for strictly positive constants $c_{1}$ and $c_{2}$.
  
  Suppose that we have a family $(E_t,\Phi_t)$ of polystable Higgs
  bundles such that $(E_0,\Phi_0) = (E,\Phi)$ and let $\alpha(t) =
  (A_t, \Phi_t)$ be the corresponding family of solutions to Hitchin's
  equations.  Then the holomorphic structure on $E_t$ is given by the
  $(0,1)$-part of $A$, while the holomorphic $\Phi$ coincides with its
  gauge theory counterpart and hence the weight of $\psi$ on
  $(\dot{A},\dot{\Phi})$ coincides with its weight on
  $(\dot{E},\dot{\Phi})$.
  
  This provides a criterion for proving that a given fixed point is
  not a minimum: it will be sufficient to find a 1-parameter family
  $(E_{t},\Phi_{t})$ of polystable Higgs bundles such that
  $(E_0,\Phi_0) = (E,\Phi)$ and such that $(\dot{E},\dot{\Phi})$ lies
  in a direct sum of strictly positive weight spaces for $\psi$. In
  fact we shall construct a family such that $\dot{\Phi} = 0$ and thus
  we simply need to consider the tangent vector $\dot{E} \in
  H^1(\End(E))$ at $0$ to the family $E_t$.
  Recall from Section~\ref{sec:critical-points}
  that we have decompositions $E' = \bigoplus F_{k}'$ and $E'' =
  \bigoplus F_{k}''$ into eigenspaces of $\psi$.    Clearly we have
  \begin{align*}
    F_{1}' &= V'\ ,  & F_{2}' &= W'\ , \\
    F_{1}'' &= W''\ , & F_{2}'' &= V'' \ .
  \end{align*}
  Let $\lambda_{V}'$ and $\lambda_{W}'$ be the weights of the action
  of $\psi$ on $V'$ and $W'$ respectively, and analogously for $E''$.
  We then have that 
  \begin{align*}
    \lambda_{W}'  &= \lambda_{V}' + 1\ , &
    \lambda_{V}'' &= \lambda_{W}'' + 1\ .  
  \end{align*}
  and, since $\tr(\psi)' = \tr(\psi)'' = 0$,
  \begin{align*}
    \lambda_{V}' p' + \lambda_{W}' q' &= 0\ , \\
    \lambda_{V}'' p'' + \lambda_{W}'' q'' &= 0\ ,    
  \end{align*}
  where $p' = \rk(V')$, $q' = \rk(W')$, $p'' = \rk(V'')$ and $q'' =
  \rk(W'')$. 
{}From these equations we conclude that 
  \begin{align}
    \lambda_{W}' - \lambda_{W}'' &=
      \frac{p'}{p' + q'} + \frac{p''}{p'' + q''} > 0\ , \\
    \lambda_{V}'' - \lambda_{V}' &=
      \frac{q''}{p'' + q''} + \frac{q'}{p' + q'} > 0\ .
  \end{align}
  Hence the subspaces $H^1(\Hom(W'',W'))$ and $H^1(\Hom(V',V''))$ of
  $H^1(\End(E))$ have strictly positive weights and it will suffice
  for us to construct a family $(E_t,\Phi_t)$ as described above such
  that $\dot{E}$ lies in the direct sum of these two spaces and
  $\dot{\Phi} = 0$.  For this argument to be valid it is of course
  essential for this direct sum to be non-zero; in fact it is proved
  in Lemma~\ref{lem:non-vanishing} below that both $H^1(\Hom(W'',W'))$
  and $H^1(\Hom(V',V''))$ are non-vanishing.  In order to find such a
  family we adapt the construction of such a family used in the proof
  of Proposition~4.3 of \cite{gothen:2001} (for $\SU(2,2)$-Higgs
  bundles):
  Let $\eta \in H^1(\Hom(V',V''))$ and $\sigma \in H^1(\Hom(W'',W))$
  be non-zero; we can then define a
  deformation of $(E,\Phi)$ by using that $\eta$ defines an
  extension
  \begin{displaymath}
    0 \lto V'' \lto V^{\eta} \lto V' \lto 0\ ,
  \end{displaymath}
  while $\sigma$ defines an extension
  \begin{displaymath}
    0 \lto W' \lto W^{\sigma} \lto W'' \lto 0\ .
  \end{displaymath}
  Let $E^{(\eta,\sigma)} = V^{\eta} \oplus W^{\sigma}$ and define
  $\Phi^{(\eta,\sigma)}$ by the compositions
  \begin{align*}
    b^{(\eta,\sigma)} &\colon
      W^{\sigma} \lto W'' \overset{b''}{\lto} V'' \to V^{\eta}\ , \\
    c^{(\eta,\sigma)} &\colon
      V^{\eta} \lto V' \overset{c'}{\lto} W' \lto W^{\sigma}.  \\
  \end{align*}
  Note that $(E^0,\Phi^0) = (E,\Phi)$ (the Higgs fields agree since
  $b' = c'' = 0$).  Now define the family $(E_{t},\Phi_{t}) =
  (E^{(\eta t,\sigma t)}, \Phi^{(\eta t,\sigma t)})$.  It is then
  clear that $\dot{E} = (\eta,\sigma)$.  In order to see that
  $\dot{\Phi} = 0$ we simply note that our definition of
  $\Phi^{(\eta,\sigma)}$ did not change the Higgs field but only the
  holomorphic structure on the bundles: thus the reason why the Higgs
  field stays holomorphic with respect to the deformed holomorphic
  structure is that $\beta' = \gamma'' = 0$.
  
  In order to conclude the proof we still need to show that $(E_t,
  \Phi_t)$ is a polystable Higgs bundle for each $t$.  It is in fact
  easy to see that $(E^{\eta,\sigma},\Phi^{\eta,\sigma})$ is stable:
  the essential point is that the destabilizing subbundles $V'$ and
  $W''$ of $(E,\Phi)$ are not subbundles of the deformed Higgs bundle.
  We leave the details to the reader.
\end{proof}
  \begin{lemma}
    \label{lem:non-vanishing}
  The cohomology groups
      $H^1(\Hom(V',V''))$ and 
      $H^1(\Hom(W'',W'))$
  are both non-vanishing.
  \end{lemma}
  \begin{proof}
    Since $\gamma'' = 0$, $V''$ is a 
    $\Phi$-invariant subbundle of $E''$ and therefore $\mu(V'') <
    \mu(E'') = \mu(E')$.  Using the Riemann--Roch formula we therefore
    obtain
  \begin{align*}
    h^0(\Hom(V',V'') - h^1(\Hom(V',V'')
    &= p'p''(1-g + \mu(V'') - \mu(V')) \\
    &< p'p''(1-g + \mu(E') - \mu(V')). \\
  \end{align*}
  Since $\rk(\beta') \leq p'$ the inequality \eqref{eq:mu-mu_1} of 
  Corollary~\ref{cor:mu-muE_i} shows that $\mu(E') - \mu(V') \leq 
  g-1$ and we therefore deduce that 
  \begin{displaymath}
    h^0(\Hom(V',V'') - h^1(\Hom(V',V'') < 0,
  \end{displaymath}
  from which it follows that $H^1(\Hom(V',V'') \neq 0$.  
  Similarly one sees that $H^1(\Hom(W'',W')) \neq 0$.
\end{proof}
In order to finish the proof of Theorem~\ref{thm:minima} we only need 
to show that any reducible $\U(p,q)$-Higgs bundle with $\beta=0$ 
or $\gamma=0$ is a local minimum of $f$.
\begin{proposition}
  \label{lem:thm-minima-4}
  Let $(E,\Phi)$ be a reducible $\U(p,q)$-Higgs bundle with
  $\beta=0$ or $\gamma=0$.  Then $(E,\Phi)$ is a local minimum of $f$.
\end{proposition}
\begin{proof}

  Recall from (\ref{minima}) that $\mathcal{N}(a,b)$ is 
the subspace of polystable
  $\U(p,q)$-Higgs bundles  in  $\mathcal{M}(a,b)$  with $\beta=0$ or 
$\gamma=0$.  {}From Proposition~\ref{prop:vanishing}
we know  that if $a/p > b/q$ then $\beta
  = 0$ and $f(\mathcal{N}(a,b)) = q(\mu(W) - \mu(E))$. 
  Since we know that $\mathcal{M}(a,b) \setminus
  \mathcal{N}(a,b)$ does not contain any local minima of $f$
  and that $f$ has a global minimum on $\mathcal{M}(a,b)$ it
  follows that this global minimum is exactly
  $f(\mathcal{N}(a,b))$ and, therefore,
  $f(\mathcal{N}(a,b))$ consists of local minima of $f$.  Of
  course a similar argument applies if $a/p < b/q$.
  
  Finally, if $a/p = b/q$ then we know from
  Proposition~\ref{prop:vanishing} that $\Phi = 0$ and hence
  $f(\mathcal{N}(a,b)) = 0$, showing that
  $\mathcal{N}(a,b)$ consists of (global) minima of the
  positive function $f$.
\end{proof}

%%%%%%%%%%%%%%%%%%%%%%%%%%%%%%%%%%%%%%%%%%%
\subsection{Local minima and connectedness}
\label{sec:minima-connectedness}
%%%%%%%%%%%%%%%%%%%%%%%%%%%%%%%%%%%%%%%%%%%

Denote by $\mathcal{M}^{s}(a,b) \subseteq \mathcal{M}(a,b)$ and
$\mathcal{N}^{s}(a,b) \subseteq \mathcal{N}(a,b)$ the subspaces
consisting of stable $\U(p,q)$-Higgs bundles, and denote by
$\bar{\mathcal{M}}^{s}(a,b)$ and $\bar{\mathcal{N}}^{s}(a,b)$ their
respective closures. In this section we explain how to obtain
connectedness results on $\mathcal{M}^{s}(a,b)$ and
$\bar{\mathcal{M}}^{s}(a,b)$.

The invariants $(a,b)$ will be fixed in the following and we shall
occasionally drop them from the notation and write $\mathcal{M} =
\mathcal{M}(a,b)$, etc.

\begin{proposition}
  \label{prop:closure-N^s}
  The closure of $\mathcal{N}^{s}$ in $\mathcal{M}$ coincides with
  $\bar{\mathcal{N}}^{s}$ and
  \begin{displaymath}
    \bar{\mathcal{N}}^{s} = \bar{\mathcal{M}}^{s} \cap \mathcal{N}\ .
  \end{displaymath}
\end{proposition}

\begin{proof}
  Clear.
\end{proof}

Now consider the restriction of the Morse function to
$\bar{\mathcal{M}}^{s}$,
\begin{displaymath}
  f \colon \bar{\mathcal{M}}^{s} \to \R\ .
\end{displaymath}
since $\bar{\mathcal{M}}^{s}$ is closed in $\mathcal{M}$ the
restriction of $f$ remains proper.

\begin{proposition}
  \label{prop:minima-closure}
  The restriction of $f$ to $\bar{\mathcal{M}}^{s}$ is proper and the
  subspace of local minima of this function coincides with
  $\bar{\mathcal{N}}^{s}$.
\end{proposition}

\begin{proof}
  Properness of the restriction follows from properness of $f$ and the
  fact that $\bar{\mathcal{M}}^{s}$ is closed in $\mathcal{M}$.

  We know that $f$ is constant on $\mathcal{N}$ and that its
  value on this subspace is the global minimum of $f$ on
  $\mathcal{M}$.  Thus $\bar{\mathcal{N}}^{s}$ is contained in the
  subspace of local minima of $f$.
  
  It remains to see that there are no other local minima of the
  restriction of $f$ to $\bar{\mathcal{M}}^{s}$.  We already know that
  the subspace of local minima on $\mathcal{M}^{s}$ is
  $\mathcal{N}^{s}$.  Now, $\mathcal{M}^{s}$ is open in
  $\bar{\mathcal{M}^{s}}$ so there cannot be any additional local
  minima on $\mathcal{M}$.
  
  Thus all we we need to prove is that there are no local minima of
  $f$ in $(\bar{\mathcal{M}^{s}} \setminus \mathcal{M}^{s}) \setminus
  \bar{\mathcal{N}}^{s}$.  So let $(E,\Phi)$ be a strictly poly-stable
  $\U(p,q)$-Higgs bundle representing a point in this space, then
  from Proposition~\ref{prop:closure-N^s} we see that
  $\beta \neq 0$ and $\gamma \neq 0$.  In the proof of
  Proposition~\ref{lem:thm-minima-3} we constructed a family
  $(E_t,\Phi_t)$ of $\U(p,q)$-Higgs bundles such that $(E,\Phi) =
  (E_0,\Phi_0)$ and $(E_t,\Phi_t)$ is stable for $t \neq 0$.
  Furthermore we showed that the restriction of $f$ to this family
  does not have a local minimum at $(E_0,\Phi_0)$.  It follows that
  $(E,\Phi)$ is not a local minimum of $f$ on $\bar{\mathcal{M}^{s}}$.
\end{proof}

\begin{proposition}
  \label{prop:N-M-list}
  \begin{itemize}
  \item[$(1)$] If $\mathcal{N}(a,b)$ is connected, then so is
    $\mathcal{M}(a,b)$.
  \item[$(2)$] If $\mathcal{N}^s(a,b)$ is connected, then so is
    $\bar{\mathcal{M}}^s(a,b)$.
  \end{itemize}
\end{proposition}

\begin{proof}
  (1) In view of Proposition~\ref{prop:topology-exercise}, this follows
      from Theorem~\ref{thm:minima}.
      
  (2) If $\mathcal{N}^s(a,b)$ is connected, then so is its closure
      $\bar{\mathcal{N}}^s(a,b)$. But from
      Proposition~\ref{prop:minima-closure},
      $\bar{\mathcal{N}}^s(a,b)$ is the subspace of local minima
      of the proper positive map $f \colon
      \bar{\mathcal{M}}^s(a,b) \to \R$.  Hence the result follows
      from
      Proposition~\ref{prop:real-topology-exercise}.
\end{proof}

%%%%%%%%%%%%%%%%%%%%%%%%%%%%%%
\section{Stable triples}
\label{sec:stable-triples}
%%%%%%%%%%%%%%%%%%%%%%%%%%%%%%
%%%%%%%%%%%%%%%%%%%%%%%%%%%%%%%%%%%%%%%%
\subsection{Definitions and basic facts}
\label{sec:triples-definitions}
%%%%%%%%%%%%%%%%%%%%%%%%%%%%%%%%%%%%%%%%
We briefly recall the relevant definitions for holomorphic triples as
studied in \cite{bradlow-garcia-prada:1996} and \cite{garcia-prada:1994}; 
we refer to these papers for details.
A \emph{holomorphic triple} on $X$,
$T = (E_{1},E_{2},\phi)$ consists of two holomorphic
vector bundles $E_{1}$ and $E_{2}$ on $X$ and a
holomorphic map $\phi \colon E_{2} \to E_{1}$.  
A homomorphism from $T' = (E_1',E_2',\phi')$ 
to $T = (E_1,E_2,\phi)$   is a commutative diagram 
\begin{displaymath}
  \begin{CD}
    E_2' @>\phi'>> E_1' \\
    @VVV @VVV  \\
    E_2 @>\phi>> E_1.
  \end{CD}
\end{displaymath}
$T'=(E_1',E_2',\phi')$ is a subtriple of $T = (E_1,E_2,\phi)$   
if the homomorphisms of sheaves $E_1'\to E_1$ and $E_2'\to E_2$ 
are injective.
For any $\alpha \in \R$ the \emph{$\alpha$-degree} and
\emph{$\alpha$-slope} of $T$ are 
defined to be
\begin{align*}
  \deg_{\alpha}(T)
  &= \deg(E_{1}) + \deg(E_{2}) + \alpha
  \rk(E_{2}), \\ 
  \mu_{\alpha}(T)
  &=
  \frac{\deg_{\alpha}(T)}
  {\rk(E_{1})+\rk(E_{2})} \\ 
  &= \mu(E_{1} \oplus E_{2}) +
  \alpha\frac{\rk(E_{2})}{\rk(E_{1})+
    \rk(E_{2})}.
\end{align*}
The triple $T = (E_{1},E_{2},\phi)$ is
\emph{$\alpha$-stable} if
\begin{displaymath}
  \mu_{\alpha}(T')
  < \mu_{\alpha}(T)
\end{displaymath}
for any proper sub-triple $T' = (E_{1}',E_{2}',\phi')$. 
Sometimes it is convenient to use
\begin{equation}\label{delta-stability}
\Delta_\alpha(T')=\mu_\alpha(T')-\mu_\alpha(T),
\end{equation}
in terms of which the $\alpha$-stability of $T$ is equivalent
to $\Delta_\alpha(T')<0$ for any proper sub-triple $T'$.
We define  \emph{$\alpha$-semistability} by replacing the above 
strict inequality with a weak inequality. A triple is called
\emph{$\alpha$-polystable} if it is the direct sum of $\alpha$-stable
triples of the same $\alpha$-slope.
Write $\mathbf{n}=(n_1,n_2)$ and
$\mathbf{d} = (d_1,d_2)$.  We denote by
\begin{displaymath}
  \mathcal{N}_\alpha
  = \mathcal{N}_\alpha(\mathbf{n},\mathbf{d})
  = \mathcal{N}_\alpha(n_1,n_2,d_1,d_2)
\end{displaymath}
the moduli space of $\alpha$-polystable triples $T =
(E_{1},E_{2},\phi)$ with $\rk(E_i)=n_i$ and $\deg(E_i) = d_i$ for
$i=1,2$. The subspace of  $\alpha$-stable triples is denoted by
$ \mathcal{N}_\alpha^s$. We refer to $(\mathbf{n},\mathbf{d})=(n_1,n_2,d_1,d_2)$ as the 
{\em type} of the triple.
Like in the Higgs bundle case, the stability condition for 
triples arises in relation to  some  gauge-theoretic equations, known
as the {\em vortex equations}. Namely,
given a triple $T=(E_1,E_2,\phi)$, one is looking for hermitian metrics $H_1$ and
$H_2$ on $E_1$ and $E_2$, respectively, such that
  \begin{equation}
    \label{eq:coupled-vortex}
  \begin{split}
    \sqrt{-1}\Lambda F(E_1) + \phi\phi^{*}
      &= \tau_{1} \Id_{E_1},  \\
    \sqrt{-1}\Lambda F(E_2)  - \phi^{*}\phi
      &= \tau_{2} \Id_{E_2}, 
  \end{split}
  \end{equation}
where $\Lambda$ is contraction by the K\"ahler form of a metric on $X$  
normalized so that $\vol(X)=2\pi$, $F(E_i)$ is the curvature of the unique connection on $E_i$ compatible
with $h_i$ and the holomorphic structure of $E_i$, and $\tau_1$ and $\tau_2$ are
real parameters satisfying $d_1+d_2=n_1\tau_1+n_2\tau_2$.
A  solution to \eqref{eq:coupled-vortex} exists if and only if $T$ is
$\alpha$-polystable for $\alpha=\tau_1-\tau_2$ (\cite{bradlow-garcia-prada:1996}).
There are certain necessary conditions in order for   $\alpha$-semistable 
triples to exist.
Let $\mu_i=d_i/n_i$ for $i=1,2$. We define
\begin{align} 
 \alpha_m= &\mu_1-\mu_2, \label{alpha-bounds-m} \\
      \alpha_M = & (1+ \frac{n_1+n_2}{|n_1 - n_2|})(\mu_1 - \mu_2), \;\;
  n_1\neq n_2. \label{alpha-bounds-M} 
\end{align}
One has the  following (\cite{bradlow-garcia-prada:1996,garcia-prada:1994}).  
\begin{proposition}
  \label{prop:alpha-range} 
The moduli space   $\mathcal{N}_\alpha(n_1,n_2,d_1,d_2)$ is a complex analytic
variety, which is projective when $\alpha$ is rational.
A necessary condition for $\mathcal{N}_\alpha(n_1,n_2,d_1,d_2)$ 
to be non-empty is 
\begin{enumerate}
\item[]
$0\leq \alpha_m \leq \alpha \leq \alpha_M$ \ \  if \ \  $n_1\neq n_2$,
\item[]
$0\leq \alpha_m \leq \alpha$ \ \    if  \ \ $n_1= n_2$.
\end{enumerate}
\end{proposition}
The moduli space of triples for $\alpha=\alpha_m$ is given by the
following.
\begin{proposition}\label{moduli-alpha_m}
A triple $T=(E_1,E_2,\phi)$ is $\alpha_m$-polystable
if and only if $\phi=0$ and $E_1$ and $E_2$ are polystable. We thus have
$$
\mathcal{N}_{\alpha_m}(n_1,n_2,d_1,d_2)\cong M(n_1,d_1) \times M(n_2,d_2).
$$
\end{proposition}
\begin{proof}
Consider  equations \eqref{eq:coupled-vortex} on $T$.
If $\alpha=\alpha_m$ then $\tau_1=\mu_1$ and $\tau_2=\mu_2$ and hence 
in order to have solutions of \eqref{eq:coupled-vortex} 
we must have  $\phi=0$. In this 
case,  \eqref{eq:coupled-vortex} say that the  hermitian metrics on 
$E_1$ and $E_2$ have constant central curvature. But this is 
equivalent to the polystability of $E_1$ and $E_2$ by the theorem of 
Narasimhan and Seshadri 
\cite{narasimhan-seshadri:1965}.
\end{proof}
\begin{remark}
If  $\alpha_m=0$ and  $n_1\neq n_2$ then
$\alpha_m=\alpha_M=0$ and the moduli space
of $\alpha$ stable triples is empty unless $\alpha=0$.
\end{remark}
Given a triple  $T=(E_1,E_2,\phi)$ one can define  the dual triple
$T^*=(E_2^*,E_1^*,\phi^*)$, where $E_i^*$ is the dual of $E_i$ and
$\phi^*$ is the transpose of $\phi$. 
It is not difficult to prove  (\cite{bradlow-garcia-prada:1996}) that the 
$\alpha$-(semi)stability of $T$ is equivalent to the 
$\alpha$-(semi)stability of $T^*$.
The map $T\mapsto T^*$ defines then an isomorphism
$$
\mathcal{N}_\alpha(n_1,n_2,d_1,d_2) = \mathcal{N}_\alpha(n_2,n_1,-d_2,-d_1).
$$
This can be used to restrict  our study to $n_1\geq n_2$ and appeal to duality to deal with the
other case.
A triple $T=(E_1,E_2,\phi)$ is strictly
$\alpha$-semistable if and only if it has a subtriple 
$T'=(E_1',E_2',\phi')$ such that 
\begin{equation}
  \label{eq:strict-alpha-ss}
  \mu(E'_1 \oplus E'_2) + \alpha \frac{n'_2}{n_1'+n_2'} 
  =  \mu(E_1 \oplus E_2) + \alpha \frac{n_2}{n_1+n_2}\ .
\end{equation}
There are two ways in which this can happen: the first one is if there
exists a subtriple $T'$ such that 
\begin{align*}
  \frac{n'_2}{n_1'+n_2'} &= \frac{n_2}{n_1+n_2}\ , \\
  \mu(E'_1 \oplus E'_2) &= \mu(E_1 \oplus E_2)\ .
\end{align*}
In this case the terms containing $\alpha$ drop from
\eqref{eq:strict-alpha-ss} and $T$ is strictly
$\alpha$-semistable for all values of $\alpha$.  We refer to this
phenomenon as \emph{$\alpha$-independent semistability}. This cannot 
happen if $\GCD(n_2,n_1+n_2,d_1+d_2)=1$. The other way in which 
strict $\alpha$-semistability can happen is if equality holds in 
\eqref{eq:strict-alpha-ss} but 
\begin{equation}
  \label{eq:not-alpha-indep}
  \frac{n'_2}{n_1'+n_2'} \neq \frac{n_2}{n_1+n_2}\ .
\end{equation}
The values of $\alpha\in [\alpha_m,\infty)$ for which this
happens are called \emph{critical values}.
For such an $\alpha$  there exists $(n'_1,n'_2,d'_1,d'_2)$
such that
\begin{displaymath}
  \frac{d'_1+d'_2}{n'_1+n'_2} + \alpha\frac{n'_2}{n'_1+n'_2}
  = \frac{d_1+d_2}{n_1+n_2} + \alpha\frac{n_2}{n_1+n_2}.
\end{displaymath}
In other words,
$$
\alpha=\frac{(n_1+n_2)(d_1'+d_2')-(n_1'+n_2')(d_1+d_2)}{n_1'n_2-n_1n_2'}
$$
with 
$n'_i \leq n_i$,  $(n'_1,n'_2,d'_1,d'_2) \neq (n_1,n_2,d_1,d_2)$,
$(n'_1,n'_2) \neq (0,0)$ and  $n_1'n_2\neq n_1n_2'$.
We say that $\alpha$ is {\em generic} if it is not critical. We thus 
have the following (cf.\ \cite{bradlow-garcia-prada:1996}). 
\begin{proposition}\label{triples-critical-range}
  \begin{itemize}
\item[(1)]
 There is only a discrete number
of critical values of $\alpha\in [\alpha_m,\infty)$ for given 
$(n_1,n_2,d_1,d_2)$.
\item[(2)] If $n_1\neq n_2$ the number of critical values is  finite and  lies in the interval $[\alpha_m,\alpha_M]$.   
\item[(3)] The stability
criteria  for two values of  $\alpha$ lying between two consecutive 
critical values
are equivalent; thus the corresponding moduli spaces are isomorphic.
\item[(4)]  If $\alpha$ is generic and  $\GCD(n_2,n_1+n_2,d_1+d_2) = 1$ 
 then $\alpha$-semistability is equivalent to $\alpha$-stability.
\end{itemize}
\end{proposition}
Let $\alpha_m^+=\alpha_m+\epsilon$, with
$\epsilon$ such that the interval $(\alpha_m, \alpha_m^+]$ does not
contain any critical value (sometimes we refer to this value of 
$\alpha$ as {\em small}. The following is important in the construction
of the moduli space for small $\alpha$.
\begin{proposition}[\cite{bradlow-garcia-prada:1996}]\label{moduli-small}If a triple $T=(E_1,E_2,\phi)$ is  $\alpha_m^+$-semistable
triple, $E_1$ and $E_2$ are semistable. In the converse direction, if $E_1$ or 
$E_2$ is stable and the other is semistable, $T=(E_1,E_2,\phi)$ is  
$\alpha_m^+$-stable.
\end{proposition}
\begin{corollary}
If $\GCD(n_1,d_1)=1$ and $\GCD(n_2,d_2)=1$,
the moduli space $\mathcal{N}_{\alpha_m^+}^s(n_1,n_2,d_1,d_2)$
is isomorphic to the projectivization of a Picard sheaf
over $M(n_1,d_1) \times M(n_2,d_2)$.
\end{corollary}
\begin{proof}
Let $\mathbb{E}_1$ and $\mathbb{E}_2$ the universal bundles over
$X\times M(n_1,d_1)$ and $X \times M(n_2,d_2)$, respectively.
Consider the canonical projections 
$\pi:X \times M(n_1,d_1)\times M(n_2,d_2)\to 
M(n_1,d_1)\times M(n_2,d_2)$;
$\hat{\pi}:X \times M(n_1,d_1)\times M(n_2,d_2)\to X$;
 $\pi_1:X \times M(n_1,d_1)\times M(n_2,d_2)\to 
X\times M(n_1,d_1)$; and 
 $\pi_2:X \times M(n_1,d_1)\times M(n_2,d_2)\to 
X\times M(n_2,d_2)$. {}From Proposition \ref{moduli-small} we deduce that
$$
\mathcal{N}_{\alpha_m^+}^s(n_1,n_2,d_1,d_2)=
\mathbb{P}
(R^1\pi_*(\pi_1^*\mathbb{E}_1\otimes \pi_2^*\mathbb{E}_2^*
\otimes\hat{\pi}^*K)^*).
$$
\end{proof}
It is important for us to have criteria to rule out
strict $\alpha$-semistability when $\alpha$ is an integer. 
\begin{lemma}
\label{critical-integer}
Let $m$ be an integer such that $\GCD(n_1+n_2,d_1+d_2-mn_1)=1$. Then 
  \begin{itemize}
  \item[(1)] $\alpha=m$ is not a critical value,
  \item[(2)] there are no $\alpha$-independent semistable triples.
  \end{itemize}
\end{lemma}
\begin{proof}
To prove (1), suppose that $\alpha=m$ is a critical value. There exist
then a triple $T$ and a proper subtriple $T'$ so that 
$$
(d_1'+d_2'+mn_2')(n_1+n_2)=(d_1+d_2+m n_2)(n_1'+n_2').
$$
This means in particular that $n_1+n_2$ divides 
$(d_1+d_2+ m n_2)(n_1'+n_2')$ and
since $n_1+n_2>n_1'+n_2'$ we must have that $\GCD(n_1+n_2,d_1+d_2+m 
n_2)>1$. But, since $d_1+d_2+mn_2=d_1+d_2-mn_1+m(n_1+n_2)$, we have  
$\GCD(n_1+n_2,d_1+d_2-m n_1)>1$, in contradiction with the 
hypothesis. To prove (2), we will show that 
$\GCD(n_2,n_1+n_2,d_1+d_2)=1$, from which the  result follows by  (4) 
in Proposition \ref{triples-critical-range}. Suppose that 
$\GCD(n_2,n_1+n_2,d_1+d_2)\neq 1$. Then there is 
$(n_2',n',d')$ such that $ \frac{n_2'}{n_1'+n_2'}= 
\frac{n_2}{n_1+n_2}$ and $\frac{d_1'+d_2'}{n_1'+n_2'} 
= \frac{d_1+d_2}{n_1+n_2}$.  
%\begin{align*}
%  \frac{n_2'}{n_1'+n_2'} &= \frac{n_2}{n_1+n_2}\ , \\
%  \frac{d_1'+d_2'}{n_1'+n_2'} &= \frac{d_1+d_2}{n_1+n_2}.
%\end{align*}
%Thus 
%\begin{align*}
%  \frac{d_1'+d_2'-mn_1' +m(n_1'+n_2')}{n_1'+n_2'} 
%&=\frac{d_1'+d_2'+mn_2'}{n_1'+n_2'}\\
%&=\frac{d_1+d_2+mn_2}{n_1+n_2}\\
%&=  \frac{d_1+d_2-mn_1 +m(n_1+n_2)}{n_1+n_2}.
%\end{align*}
%This implies that
It follows that 
$$ 
\frac{d_1'+d_2'-mn_1'}{n_1'+n_2'}=\frac{d_1+d_2-mn_1}{n_1+n_2},
$$
and hence $\GCD(n_1+n_2,d_1+d_2-mn_1)\neq 1$, in contradiction with 
the hypothesis. 
\end{proof}
 
%%%%%%%%%%%%%%%%%%%%%%%%%%%%%%%%%%%%%%%%
\subsection{Minima as triples}
\label{sec:minima-as-triples}
%%%%%%%%%%%%%%%%%%%%%%%%%%%%%%%%%%%%%%%%
  
Let $(E,\Phi)$ be a $\U(p,q)$-Higgs bundle with $\gamma = 0$.  We 
can then define a holomorphic triple 
$T = (E_{1},E_{2},\phi)$ by 
\begin{align*}
  E_{1} &= V\otimes K \\
  E_{2} &= W, \\
  \phi &= \beta;
\end{align*}
and given a holomorphic triple we can define an associated
$\U(p,q)$-Higgs bundle with $\gamma = 0$.  Similarly, there is a 
bijective correspondence between $\U(p,q)$-Higgs bundles with 
$\beta = 0$ and holomorphic triples given by
\begin{align*}
  E_{1} &= W\otimes K, \\
  E_{2} &= V \\
  \phi &= \gamma.
\end{align*}
The link between the stability conditions for holomorphic triples and
$\U(p,q)$-Higgs bundles is given by the following result.
\begin{proposition}
  \label{prop:triple-higgs-stability}
  A $\U(p,q)$-Higgs bundle $(E,\Phi)$ with $\beta =0$ or $\gamma=0$ is
  \mbox{(semi)}stable if and only if the corresponding holomorphic
  triple $T = (E_{1},E_{2},\phi)$ is
  $\alpha$-(semi)stable for $\alpha = 2g-2$.
\end{proposition}
\begin{proof}
  For definiteness we shall assume that $\gamma=0$ (of course, the same 
  argument applies if $\beta =0$).  We then have $E_{1} = V 
  \otimes K$ and $E_{2} = W$ and, hence,
  \begin{displaymath}
    \deg(E_{1}) = \deg(V) + p(2g-2).
  \end{displaymath}
  Since $p = \rk(E_{1})$ and $q =
  \rk(E_{2})$ it follows that
  \begin{align*}
    \mu_{\alpha}(T)
    &=\frac{\deg(V) + \deg(W) + p(2g-2) + \alpha 
    q}{p + q} \\
    &= \mu(E) + \frac{p}{p+q}(2g-2) + 
    \frac{q}{p+q}\alpha.
  \end{align*}
  If we set $\alpha = 2g-2$ we therefore have
  \begin{equation}
    \label{eq:mu-alpha}
    \mu_{\alpha}(T) = \mu(E) + 2g-2.
  \end{equation}
  Clearly the correspondence between holomorphic triples and
  $\U(p,q)$-Higgs bundles gives a correspondence between sub-triples
  $T' = (E_{1}',E_{2}',\phi')$ and $\Phi$-invariant subbundles of $E$ which
  respect the decomposition $E = V \oplus W$ (i.e., subbundles $E' =
  V' \oplus W'$ with $V' \subseteq V$ and $W' \subseteq W$).  It was
  shown in Section 2.3 of \cite{gothen:2001} that $(E,\Phi)$ is
  (semi)stable if and only if the (semi)stability condition holds
  for $\Phi$-invariant subbundles which respect the decomposition $E =
  V \oplus W$.  On the other hand, it follows from \eqref{eq:mu-alpha}
  that
  \begin{displaymath}
    \mu(E') < \mu(E)
  \end{displaymath}
  if and only if
  \begin{displaymath}
    \mu_{\alpha}(T') <
    \mu_{\alpha}(T)
  \end{displaymath}
  (and similarly for semistability), thus concluding the proof.
\end{proof}

Recall that 
${\mathcal{M}}(a,b)$ is the space  of polystable $\U(p,q)$-Higgs
bundles with $\deg(V) = a$ and $\deg(W) = b$. 
We have the following important 
characterization of the subspace of local minima of $f$ on 
$\mathcal{M}(a,b)$.
\begin{theorem}
  \label{thm:minima=triple-moduli}
  Let $\mathcal{N}(a,b)$ be the subspace of local minima of
  $f$ on $\mathcal{M}(a,b)$.
  \begin{itemize}
  \item[$(1)$] If $a/p \leq b/q$ then
  \begin{displaymath}
    \mathcal{N}(a,b)
    \cong \mathcal{N}_{2g-2}(p,q,a + p(2g-2),b),
  \end{displaymath}
  the moduli space of $\alpha$-polystable triples with
$n_1 = p$,    $n_2 = q$,    $d_1 = a + p(2g-2)$,     $d_2 = b$, 
 and $\alpha = 2g-2$.
  \item[$(2)$] If $a/p \geq b/q$ then
  \begin{displaymath}
    \mathcal{N}(a,b)
    \cong \mathcal{N}_{2g-2}(q,p,b + q(2g-2),a),
  \end{displaymath}
   the moduli space of $\alpha$-polystable triples with
   $n_1 = q$, $n_2 = p$,    $d_1 = b + q(2g-2)$, 
    $d_2 = a$,    and $\alpha = 2g-2$.
   \item[$(3)$]
  If $a/p = b/q$ then
  \begin{displaymath}
    \mathcal{N}(a,b)
    \cong M(p,a) \times M(q,b),
  \end{displaymath}
  where $M(p,a)$ is the moduli space of polystable bundles of rank 
  $p$ and degree $a$, and $M(q,b)$ is the moduli space of 
  polystable bundles of rank $q$ and degree $b$.
  \end{itemize}
\end{theorem}
\begin{proof}
It is obtained by combining 
Theorem~\ref{thm:minima}, Proposition~\ref{prop:vanishing}, 
Proposition~\ref{prop:triple-higgs-stability} and 
Corollary~\ref{minima-toledo-0}. 
\end{proof}
\begin{remark}
  If  $a/p = b/q$ then $2g-2$ is equal to the extreme value, $\alpha_m$,  
  for triples of type
  $(p,q,a + p(2g-2),b)$ and also for triples of type $(q,p,b + q(2g-2),a)$.
  Thus by Proposition \ref{moduli-alpha_m} the descriptions of
  $\mathcal{N}(a,b)$ in (1) and (2) of Theorem 
  \ref{thm:minima=triple-moduli} coincide with each other
  and agree with that in (3).
\end{remark}

\begin{theorem}
  \label{thm:connected-triples-connected-higgs}

  \begin{itemize}
    
  \item[$(1)$] Let $a/p \leq b/q$.
    \begin{itemize}
    \item     If $\mathcal{N}_{2g-2}(p,q,a + p(2g-2),b)$
    is connected then $\mathcal{M}(a,b)$ is connected.
  \item If $\mathcal{N}^{s}_{2g-2}(p,q,a + p(2g-2),b)$ is connected
    then $\bar{\mathcal{M}}^s(a,b)$ is connected.
    \end{itemize}
  \item[$(2)$] Let $a/p \geq b/q$.
    \begin{itemize}
    \item    If $\mathcal{N}_{2g-2}(q,p,b + q(2g-2),a)$ is connected
    then $\mathcal{M}(a,b)$ is connected.
  \item     If $\mathcal{N}^{s}_{2g-2}(q,p,b + q(2g-2),a)$ is connected
    then $\bar{\mathcal{M}}^s(a,b)$ is connected.
    \end{itemize}
   \item[$(3)$]
    If $a/p = b/q$ then $\mathcal{M}(a,b)$ and
    $\bar{\mathcal{M}}^s(a,b)$ are connected.
  \end{itemize}
\end{theorem}

\begin{proof}
  It follows from Proposition~\ref{prop:N-M-list} and
  Theorem~\ref{thm:minima=triple-moduli}.  For (3) we also use that
  the moduli spaces of semistable bundles $M(p,a)$ and $M(q,b)$ are
  connected, and that the moduli spaces of stable bundles $M^s(p,a)$
  and $M^s(q,b)$ are connected (in fact these moduli spaces are
  irreducible.)
\end{proof}

In view of Theorem~\ref{thm:minima=triple-moduli}, it is important to 
understand where  
$2g-2$ lies in the range (given by Proposition \ref{prop:alpha-range})
for the stability parameter $\alpha$. 

\begin{proposition}\label{prop:MW}
Let $(n_1,n_2,d_1,d_2)$ be the type of the triples arising from 
$\mathcal{N}(a,b)$, as in Theorem~\ref{thm:minima=triple-moduli}.
Let $\tau$ be the Toledo invariant as in \eqref{toledo-invariant}. 
\begin{itemize}
\item[$(1)$] $ 2g-2\geq \alpha_m$, with equality  
if and only if $\tau=0$.
\item[$(2)$] If $p\neq q$ then the 
Milnor--Wood inequality $|\tau|\leq \min\{p,q\}(2g-2)$ is equivalent 
to the condition $2g-2\leq \alpha_M$. Moreover, 
$2g-2=\alpha_M$ if and only  if $|\tau|= \min\{p,q\}(2g-2)$. 
\item[$(3)$] If $p=q$, then the 
Milnor--Wood inequality is equivalent to the condition $\alpha_m\geq 
0$. Moreover $\alpha_m=0$ if and only if $|\tau|=p(2g-2)$. 
\end{itemize}
\end{proposition}
\begin{proof}
Let $(E,\Phi)\in\mathcal{N}(a,b)$. We assume  that 
$a/p\geq b/q$, i.e.\ $\tau\geq 0$ 
(the case $a/p\leq b/q$ is analogous).
Then $\mu_1=b/q+ 2g-2$, $\mu_2=a/p$, and hence
$\alpha_m=\mu_1-\mu_2=b/q- a/p +2g-2$. We thus have
$\alpha_m\leq 2g-2$.
Hence, the Toledo invariant  \eqref{toledo-invariant} is
\begin{equation}
\label{toledo-alpha_m}
\tau= \frac{2pq}{p+q}(\frac{a}{p}-\frac{b}{q})=\frac{2pq}{p+q}(2g-2 -\alpha_m),
\end{equation}
from which  we see that $\alpha_m\leq 2g-2$ is equivalent to 
$\tau\geq 0$ and, in particular, $\alpha_m= 2g-2$ if and only if
$\tau= 0$. This proves (1).
{}From (\ref{toledo-alpha_m}) we see that $\alpha_m\geq 0$ is
equivalent to 
\begin{equation}
\label{weak-toledo}
\tau\leq \frac{2pq}{p+q}(2g-2).
\end{equation}
If $p=q$ then (\ref{weak-toledo}) is equivalent to $\tau\leq p(2g-2)$, which
is the Milnor--Wood inequality, proving (2). If $p\neq q$, then 
(\ref{weak-toledo}) is weaker than the Milnor--Wood inequality. 
But in this case $\alpha_M$ is an upper bound for $\alpha$ and in order
to have existence we must assume
\begin{equation}
\label{existence-Higgs-triples}
2g-2\leq \alpha_M.
\end{equation}
Suppose that $p> q$ (the case $p<q$ is analogous), then
$\alpha_M=\frac{2p}{p-q}\alpha_m$ and hence
 (\ref{toledo-alpha_m}) 
can be written as
\begin{equation}
\label{toledo-alpha_M}
\tau=\frac{2pq}{p+q}(2g-2 -\frac{p-q}{2p}\alpha_M).
\end{equation}
Hence, from (\ref{toledo-alpha_M}) we see that
 (\ref{existence-Higgs-triples}) is equivalent to  $\tau\leq q(2g-2)$,
and in particular we have $2g-2=\alpha_M$  if and only if 
$\tau= q(2g-2)$. Which conludes the proof of (3).
\end{proof}
\begin{remark}\label{remark:MWbound}
The above proposition gives another explanation for the Milnor--Wood 
inequality in Corollary \ref{cor:toledo}. Using the fact that the 
non-emptiness of 
$\mathcal{M}(a,b)$ is equivalent to the non-emptiness of either
$\mathcal{N}_{2g-2}(p,q,a+p(2g-2),b)$ or 
$\mathcal{N}_{2g-2}(q,p,b+q(2g-2),a)$, we see that the Milnor--Wood
inequality is equivalent to the condition that 
$2g-2$ lies within the range where  $\alpha$-polystable triples of the 
given kind exist. {}From this, we see that  in order to study 
$\mathcal{N}(a,b)$ for different values of the Toledo invariant, we 
have to study the moduli space of triples for $\alpha$ that may be 
lying anywhere in  the $\alpha$-range, including in the extreme 
values $\alpha_m$ and $\alpha_M$. 

The rest of the paper is devoted to this question. The strategy we 
adopt is the following: When 
$p\neq q$, and assuming that $2g-2$ is not an extreme, we give a 
explicit birational  description of  what we call the {\em large} 
moduli space. This is the moduli space for 
$\alpha_M^-=\alpha_M-\epsilon$, for 
$\epsilon>0$ small enough so that there is no any  critical value in
 $[\alpha_M^-,\alpha_M)$.  Once we have identified this,
we need to study the changes that happen when we cross a critical
value. We do this with all the critical values  until we get to 
$\alpha=2g-2$. When $p=q$ there is no upper bound for $\alpha$ but it 
turns out that the moduli space remains unchanged after a certain 
value of $\alpha$ (Sec. \ref{sec:stabilize}). We can thus  apply the 
same principle as before taking this as the {\em large} moduli space. 
\end{remark}
%%%%%%%%%%%%%%%%%%%%%%%%%%%%%%%%%%%%%%%%%%%%%%%%%%%
\subsection{Extensions and deformations of triples}
\label{sec:extensions-of-triples}
%%%%%%%%%%%%%%%%%%%%%%%%%%%%%%%%%%%%%%%%%%%%%%%%%%%
Let $T'=(E'_1,E'_2,\phi')$ and $T''=(E''_1,E''_2,\phi'')$ be two triples
and, as usual, let
\begin{align*}
  (\mathbf{n}',\mathbf{d}') &= (n_{1}',n_{2}',d_{1}',d_{2}'), \\
  (\mathbf{n}'',\mathbf{d}'') &= (n_{1}'',n_{2}'',d_{1}'',d_{2}''),
\end{align*}
where $n'_i = \deg(E'_i)$, $n''_i = \deg(E''_i)$, $d'_i = \deg(E'_i)$
and $d''_i = \deg(E''_i)$.  Let $\Hom(T'',T')$ denote the linear  
space of homomorphisms from 
$T''$ to
$T'$, and let $\Ext^1(T'',T')$  denotes the linear space of extensions 
of the form
\begin{displaymath}
  0 \lto T' \lto T \lto T'' \lto 0.
\end{displaymath}
In order to analyse $\Ext^1(T'',T')$ one considers the
complex of sheaves 
\begin{equation}
  \label{eq:extension-complex}
    C^{\bullet}(T'',T') \colon {E_{1}''}^{*} \otimes E_{1}' \oplus
  {E_{2}''}^{*} \otimes E_{2}'
  \overset{c}{\lto}
  {E_{2}''}^{*} \otimes E_{1}',
\end{equation} 
where the map $c$ is defined by
\begin{displaymath}
  c(\psi_{1},\psi_{2}) = \phi'\psi_{2} - \psi_{1}\phi''.
\end{displaymath}
\begin{proposition}
  \label{prop:hyper-equals-hom}
  There are natural isomorphisms
  \begin{align*}
    \Hom(T'',T') &\cong \HH^{0}(C^{\bullet}(T'',T')), \\
    \Ext^{1}(T'',T') &\cong \HH^{1}(C^{\bullet}(T'',T')),
  \end{align*}
and a long exact sequence associated to the complex
$C^{\bullet}(T'',T')$:
\begin{equation}
  \label{eq:long-exact-extension-complex}
\begin{array}{ccccccc}
  0 &\lto \mathbb{H}^0(C^{\bullet}(T'',T')) &
  \lto & H^0({E_{1}''}^{*} \otimes E_{1}' \oplus {E_{2}''}^{*} \otimes E_{2}')
  & \lto &  H^0({E_{2}''}^{*} \otimes E_{1}') \\
    &  \lto \mathbb{H}^1(C^{\bullet}(T'',T')) &
  \lto &  H^1({E_{1}''}^{*} \otimes E_{1}' \oplus {E_{2}''}^{*} \otimes E_{2}')
&  \lto & H^1({E_{2}''}^{*} \otimes E_{1}') \\
&   \lto \mathbb{H}^2(C^{\bullet}(T'',T')) & \lto & 0. & & 
\end{array}
\end{equation}
\end{proposition}
\begin{proof}
A proof can be given by using cocycle representatives of the hypercohomology 
 classes, in a manner similar to what is done in the  study of  deformations
(\cite{biswas-ramanan:1994}). In fact the result  is a  special case of a much 
more general result proved in  \cite{gothen-king:2002}. 
\end{proof}
We introduce the following notation:
\begin{align*}
  h^{i}(T'',T') &= \dim\HH^{i}(C^{\bullet}(T'',T')), \\
  \chi(T'',T') &= h^0(T'',T') - h^1(T'',T') + h^2(T'',T').
\end{align*}
\begin{proposition}
  \label{prop:chi(T'',T')}
  For any holomorphic triples $T'$ and $T''$ we have
  \begin{align*}
    \chi(T'',T') &= \chi({E_{1}''}^{*} \otimes E_{1}')
    + \chi({E_{2}''}^{*} \otimes E_{2}')
    - \chi({E_{2}''}^{*} \otimes E_{1}')  \\
    &= (1-g)(n''_1 n'_1 + n''_2 n'_2 - n''_2 n'_1) \\
    & \quad + n''_1 d'_1 - n'_1 d''_1
    + n''_2 d'_2 - n'_2 d''_2
    - n''_2 d'_1 + n'_1 d''_2.
  \end{align*}
\end{proposition}
\begin{proof}
  Immediate from the long exact sequence
  \eqref{eq:long-exact-extension-complex} and the Riemann--Roch formula.
\end{proof}
\begin{corollary}
  \label{cor:chi-relation}
  For any extension $0 \to T' \to T \to T'' \to 0$ of triples,
  \begin{displaymath}
    \chi(T,T) = \chi(T',T') + \chi(T'',T'') + \chi(T'',T') +
    \chi(T',T''). 
  \end{displaymath}
  \qed
\end{corollary}
\begin{remark} Proposition \ref{prop:chi(T'',T')} shows that 
$\chi(T'',T')$ depends only  on the topological invariants 
$(\mathbf{n}',\mathbf{d}')$ and 
$(\mathbf{n}'',\mathbf{d}'')$ of $T'$ and $T''$.  Whenever convenient
we shall therefore use the notation
\begin{displaymath}
  \chi(\mathbf{n}'',\mathbf{d}'',\mathbf{n}',\mathbf{d}')
  = \chi(T'',T').
\end{displaymath}
\end{remark}
\noindent The following vanishing results for $\HH^{0}$ and $\HH^{2}$  
play a central role in our problem.
 \begin{proposition}
   \label{prop:h0-vanishing}
   Suppose that $T'$ and $T''$ are $\alpha$-semistable. 
   
   \begin{enumerate}
   \item If $\mu_\alpha(T')<\mu_\alpha (T'')$ then  
 $\HH^{0}(C^{\bullet}(T'',T')) \cong 0$
  \item If $\mu_\alpha(T')=\mu_\alpha (T'')$ and $T''$ is
   $\alpha$-stable, then
  $$
     \HH^{0}(C^{\bullet}(T'',T')) \cong 
     \begin{cases}
       \C \quad &\text{if $T' \cong T''$} \\
       0 \quad &\text{if $T' \not\cong T''$}. 
     \end{cases}
 $$
 \end{enumerate}
   \end{proposition}
 \begin{proof} By Proposition \ref{prop:hyper-equals-hom} we can 
 identify $\HH^{0}(C^{\bullet}(T'',T'))$ with $\Hom(T'',T')$. The 
 statements (1) and (2) are thus the direct analogs for triples of the 
 same results for semistable bundles. The proof is identical. Suppose 
 that we can find a $h:T''\rightarrow T'$ is a non-trivial homorphism 
 of triples $h:T''\rightarrow T'$. If 
 $T'=(E_1',E_2',\Phi')$ and $T''=(E_1'',E_2'',\Phi'')$ then $h$ is 
 given by a pair of bundle maps $u_i:E_i''\rightarrow E_i''$ for 
 $i=1,2$ such that $\Phi'\circ u_2=u_1\circ\Phi''$. We can thus define
 subtriples of $T''$ and $T'$ respectively by 
 $T_N=(\ker(u_1),\ker(u_2),\Phi'')$ and $T_I=(\im(u_1),\im(u_2),\Phi')$,  where in $T_I$, it is in general necessary to take the saturations of 
 the image $\im(u_1)$ and $\im(u_2)$. By the semistability conditions, 
 we get
 $$
 \mu_\alpha(T_N)\le\mu_{\alpha}(T'')\le\mu_{\alpha}(T_I)\le\mu_\alpha(T')\ .
 $$ 
 The conclusions follow directly from this.
 \end{proof}
\begin{proposition}
  \label{prop:h2-vanishing}
Suppose that triples $T'$ and $T''$ are $\alpha$-semistable with 
$\mu_\alpha(T')=\mu_\alpha (T'')$. Then 
\begin{enumerate}
\item $\HH^{2}(C^{\bullet}(T'',T')) = 0$ whenever $\alpha > 2g-2$.
\item If one of $T'$, $T''$ is   $\alpha+\epsilon$-stable for some
$\epsilon\ge 0$, then $\HH^{2}(C^{\bullet}(T'',T')) = 0$ whenever
$\alpha \geq2g-2$.
\end{enumerate}
\end{proposition}
\begin{proof}
{}From (\ref{eq:long-exact-extension-complex}) it is clear that the 
vanishing of 
$\HH^{2}(C^{\bullet}(T'',T'))$ is equivalent to the surjectivity of the map
$$
  H^1({E_{1}''}^{*} \otimes E_{1}' \oplus {E_{2}''}^{*} \otimes E_{2}')
  \lto  H^1({E_{2}''}^{*} \otimes E_{1}')\ .
$$
By Serre duality this is equivalent to the injectivity of the map
\begin{equation}
  \label{eq:petri}
\begin{array}{ccc}
H^0({E_1'}^{*} \otimes E_{2}''\otimes K) &   \overset{P}{\lto} &
H^0({E_{1}'}^{*} \otimes E_{1}''\otimes K) \oplus H^0({E_{2}'}^{*}
\otimes E_{2}''\otimes K)\\
\psi &\longmapsto & ((\phi''\otimes \Id)\circ \psi, \psi\circ \phi').
\end{array}
\end{equation}
\noindent {\it Proof of (1):} Suppose that $P$ is not injective. Then there is 
a non-trivial homomorphism $\psi:E_1' \to E_2''\otimes K$ in 
$\ker P$. Let $I=\im\psi$ and $N=\ker\psi$. Since 
$(\phi''\otimes \Id_K)\circ\psi=0$, $I\subset\ker\phi''$ and hence 
$T_I''=(0, I\otimes K^*, 0)$ is a proper subtriple of  $T''$. Similarly,
the fact that $\psi\circ \phi'=0$ implies that $\im\phi'\subset N$ and thus
 $T_N'=(\ker\psi,E_2',\phi')$ is a proper subtriple of $T'$.  Let $k=\rk(N)$ and 
$l=\deg(N)$. Then, from the exact sequence
\begin{displaymath}
  0 \lto N \lto E_1' \lto I \lto 0.
\end{displaymath}
we see that $\rk(I)=n_1'-k$ and $\deg(I)=d_1'-l$. Hence
\begin{align*}
 \mu_\alpha(T_N') &= \frac{l+d_2'}{k+n_2'}+\alpha \frac{n_2'}{k+n_2'}, \\
\mu_\alpha(T_I'') &= \frac{d_1'-l}{n_1'-k}+ 2-2g +\alpha.
\end{align*}
Adding these two expressions, and clearing denominators we see that
\begin{displaymath}
d_1'+d_2'+(n_1'-k)(2-2g)+\alpha(n_1'+n_2'-k)= 
(k+n_2')\mu_\alpha(T_N')+ (n_1'-k)\mu_\alpha (T_I'').
\end{displaymath}
But  $\mu_\alpha(T_N')\leq \mu_\alpha(T')$, $\mu_\alpha(T_I'')\leq
 \mu_\alpha(T'')$ and $\mu_\alpha(T')=\mu_\alpha(T'')$.
{}From this we obtain that 
\begin{equation}\label{eqn:alpha-est}
d_1'+d_2'+(n_1'-k)(2-2g)+\alpha(n_1'+n_2'-k)\leq d_1'+d_2'+\alpha 
n_2'\ , 
\end{equation}
and hence
\begin{displaymath}
\alpha(n_1'-k)\leq(n_1'-k)(2g-2)\ .
\end{displaymath}
Since  $n_1'-k>0$ we get  that  $\alpha\leq 2g-2$. Hence $P$ must be 
injective if the hypotheses of the part (1) of the proposition are 
satisfied. 

\noindent {\it Proof of (2):} Suppose that $T''$ is $\alpha+\epsilon$-stable 
for some $\epsilon\ge 0$. It follows that 
$\mu_{\alpha+\epsilon}(T_I'')< 
\mu_{\alpha+\epsilon}(T'')$, i.e
$$
\mu_{\alpha}(T_I'')-
\mu_{\alpha}(T'')<\epsilon(\frac{n_2''}{n_1''+n_2''}-1)\le 0\ .
$$
Thus, following exactly the same argument as in the proof of (1), we 
get a strict inequality in (\ref{eqn:alpha-est}).  We conclude that 
that if P is not injective then $\alpha< 2g-2$, i.e.\ if $\alpha\ge 
2g-2$ then $P$ must be injective.
If  $T'$ is $\alpha+\epsilon$-stable for some $\epsilon\ge 0$ then we 
get that
$$
\mu_{\alpha}(T_N')-
\mu_{\alpha}(T')<\epsilon(\frac{n_2'}{n_1'+n_2'}-\frac{n_2'}{k+n_2'})\le 0\ .
$$
The rest of the argument is the same as in the case that $T''$ is 
$\alpha+\epsilon$-stable.
\end{proof}
\begin{corollary}\label{dimension-ext1}
Let $T'$ and $T''$ be two holomorphic triples.
\begin{itemize}
\item[(1)]
$\dim \Ext^1(T'',T')=h^0(T'',T')+h^2(T'',T')- \chi(T'',T')$.
\item[(2)] If $T'$ and $T''$ are  $\alpha$-semistable,  
 $\mu_\alpha(T')=\mu_\alpha (T'')$,
and $\alpha> 2g-2$, then
 $$
\dim \Ext^1(T'',T')=    h^0(T'',T')  - \chi(T'',T').
$$
The same holds for $\alpha\geq 2g-2$ if in addition $T'$ or $T''$ is 
$\alpha+\epsilon$-stable for some $\epsilon\ge 0$.
\end{itemize}
\end{corollary}
Since the  space of infinitesimal deformations of $T$ is
isomorphic to $\HH^{1}(C^{\bullet}(T,T))$, these
considerations also apply to studying deformations of a holomorphic
triple $T$.  To be precise, one has the following.
\begin{theorem}\label{thm:smoothdim}
Let $T=(E_1,E_2,\phi)$ be an $\alpha$-stable triple of type 
$(n_1,n_2,d_1,d_2)$.
\begin{itemize}
\item[(1)] The Zariski tangent space at the point defined by $T$ is isomorphic
to
$$
\HH^{1}(C^{\bullet}(T,T)).
$$
\item[(2)]
If $\HH^{2}(C^{\bullet}(T,T))= 0$, then the moduli space of $\alpha$-stable
triples is smooth in  a neighbourhood of the point defined by $T$.
\item[(3)]
$\HH^{2}(C^{\bullet}(T,T))= 0$ if and only if the homomorphism
$$
  H^1(E_1^* \otimes E_1 \oplus E_2^* \otimes E_2)
  \lto  H^1(E_2^* \otimes E_1)
$$
in the corresponding long exact sequence is surjective.
\item[(4)] At a smooth point $T\in \mathcal{N}^s_\alpha(n_1,n_2,d_1,d_2)$
the dimension of the moduli space of $\alpha$-stable triples is
\begin{align}
  \dim \mathcal{N}^s_\alpha(n_1,n_2,d_1,d_2)
  &= h^{1}(T,T) = 1 - \chi(T,T) \notag \\
  &= (g-1)(n_1^2 + n_2^2 - n_1 n_2) - n_1 d_2 + n_2 d_1 + 1.
    \label{eq:dim-triples}
\end{align}
\item[(5)] If $\phi$ is injective or surjective then $T=(E_1,E_2,\phi)$
defines a smooth point in the moduli space.
\item[(6)] If $\alpha\geq 2g-2$, then $T$ defines a smooth point
in the moduli space, and hence $\mathcal{N}^s_\alpha(n_1,n_2,d_1,d_2)$ is 
smooth.
\end{itemize}
\end{theorem}
\begin{proof}
Statements (1) and (2) are proved in \cite{bradlow-garcia-prada:1996} 
(see also \cite{biswas-ramanan:1994}).
(3) follows from the long exact sequence (\ref{eq:long-exact-extension-complex}) 
with $T=T'=T''$. (4) follows from (1), (2) and 
Propositions~\ref{prop:chi(T'',T')} and \ref{dimension-ext1}.
(5) is proved in \cite[Prop. 6.3]{bradlow-garcia-prada:1996}. (6) is a 
consequence of Proposition~\ref{prop:h2-vanishing}.
\end{proof}
%%%%%%%%%%%%%%%%%%%%%%%%%%%%%%%%%%%%%%
\section{Crossing critical values}
\label{crossing-critical-values}
%%%%%%%%%%%%%%%%%%%%%%%%%%%%%%%%%%%%%%
In this section we study the differences between the stable loci 
$\mathcal{N}^s_\alpha(\mathbf{n},\mathbf{d})$ in the moduli spaces 
$\mathcal{N}_\alpha(\mathbf{n},\mathbf{d})$, for fixed
values of $\mathbf{n}=(n_1,n_2)$ and 
$\mathbf{d} = (d_1,d_2)$ but different values of $\alpha$. Since
in this section $\mathbf{n}$ and $\mathbf{d}$ are fixed, we use the 
abbreviated notation 
$$\mathcal{N}^s_{\alpha}=\mathcal{N}^s_\alpha(\mathbf{n},\mathbf{d})
\ \quad\ \mathrm{and}\ \quad\ 
\mathcal{N}_{\alpha}=\mathcal{N}_\alpha(\mathbf{n},\mathbf{d})\ .$$
\noindent Our main result is that for all $\alpha\ge 2g-2$ 
any differences between the $\mathcal{N}^s_{\alpha}$ are confined to 
subvarieties of positive codimension. In particular, the number of 
components and the irreducibility properties of the spaces 
$\mathcal{N}^s_{\alpha}$ are the same for all $\alpha$ satisfying 
$\alpha\ge 2g-2$ and $\alpha_m<\alpha <\alpha_M$\footnotemark .
\footnotetext{When $n_1\ne n_2$ the bounds $\alpha_m$ and $\alpha_M$ 
are as in (\ref{alpha-bounds-m}) and (\ref{alpha-bounds-m}). When 
$n_1=n_2$ we adpot the convention that $\alpha_M=\infty$}If 
the coprimality condition $\GCD(n_2,n_1+n_2,d_1+d_2)=1$ is satisfied, 
then $\mathcal{N}^s_{\alpha}=\mathcal{N}_{\alpha}$ at all 
non-critical vales of $\alpha$, so the results apply to   
$\mathcal{N}_{\alpha}$ for all non-critical $\alpha\ge 2g-2$.

We begin with a set theoretic description of the differences between 
two spaces $\mathcal{N}^s_{\alpha_1}$ and  $\mathcal{N}^s_{\alpha_2}$ 
when $\alpha_1$ and $\alpha_2$ are separated by a critical value (as 
defined in section \ref{sec:triples-definitions}). For the rest of 
this section we adopt the following notation: Let $\alpha_c$ be a 
critical value such that 
\begin{equation}\label{eqn:alphac-range}
\alpha_m <\alpha_c <\alpha_M\ .
\end{equation}
\noindent Set 
\begin{equation}\label{eqtn:alpha-c-pm}
\alpha_c^+ = \alpha_c + \epsilon\quad ,\quad 
\alpha_c^- = \alpha_c - \epsilon\ ,
\end{equation}
\noindent where $\epsilon > 0$ is small enough 
so that $\alpha_c$ is the only critical value in the interval 
$(\alpha_c^-,\alpha_c^+)$. 
%%%%%%%%%%%%%%%%%%%%%
\subsection{Flip Loci}\label{subs:fliploci}
%%%%%%%%%%%%%%%%%%%%%%
\begin{definition}Let $\alpha_c\in (\alpha_m,\alpha_M)$ 
be a critical value for triples of type $(\mathbf{n},\mathbf{d})$.  
We define {\it flip loci} 
$\mathcal{S}_{\alpha_c^{\pm}}\subset\mathcal{N}^s_{\alpha_c^{\pm}}$ 
by the conditions that the points in $\mathcal{S}_{\alpha_c^+}$ represent 
triples which are $\alpha_c^+$-stable  but $\alpha_c^-$-unstable, 
while the points in $\mathcal{S}_{\alpha_c^-}$ represent triples 
which are 
$\alpha_c^-$-stable  but $\alpha_c^+$-unstable.
\end{definition}

\begin{remark}The definition of $\mathcal{S}_{\alpha_c^+}$ can be 
extended to the extreme case $\alpha_c=\alpha_m$. However, since all 
$\alpha_m^+$-stable triples must be $\alpha_m^-$-unstable, we see that
$\mathcal{S}_{\alpha_m^+}=\mathcal{N}^s_{\alpha_m^+}$. Similarly, when 
$n_1\ne n_2$ we get $\mathcal{S}_{\alpha_M^-}=\mathcal{N}^s_{\alpha_M^-}$.
 The only interesting cases are thus those those for which 
 $\alpha_m<\alpha_c <\alpha_M$. 
\end{remark}

\begin{lemma}\label{lemma:fliploci}In the above notation:
\begin{equation}\label{eqn:Nsalpha}
\mathcal{N}^s_{\alpha_c^+}-\mathcal{S}_{\alpha_c^+}=
\mathcal{N}^s_{\alpha_c}=
\mathcal{N}^s_{\alpha_c^-}-\mathcal{S}_{\alpha_c^-}\ .
\end{equation}
\end{lemma}

\begin{proof} By definition we can identify 
$\mathcal{N}^s_{\alpha_c^+}-\mathcal{S}_{\alpha_c^+}=
\mathcal{N}^s_{\alpha_c^-}-\mathcal{S}_{\alpha_c^-}$.

Suppose now that $t$ is a point in 
$\mathcal{N}^s_{\alpha_c^+}-\mathcal{S}_{\alpha_c^+}=
\mathcal{N}^s_{\alpha_c^-}-\mathcal{S}_{\alpha_c^-}$, but that $t$ is not
in $\mathcal{N}^s_{\alpha_c}$. Let $T$ be a triple representing $t$. 
Then $T$ has a subtriple $T'\subseteq T$ for which 
$\mu_{\alpha_c}(T')\ge\mu_{\alpha_c}(T)$, and also 
$\mu_{\alpha_c^{\pm}}(T')<\mu_{\alpha_c^{\pm}}(T)$.
This is not possible, and hence $t\in\mathcal{N}^s_{\alpha_c}$. 

Finally, suppose that $t\in\mathcal{N}^s_{\alpha_c}$ and let $T$ be a 
triple representing $t$. Then $\mu_{\alpha_c}(T')<\mu_{\alpha_c}(T)$ 
for all subtriples $T'\subset T$. But since the set of possible 
values for 
$\mu_{\alpha_c}(T')$ is a discrete subset of $\R$, we can find a 
$\delta>0$ such that $\mu_{\alpha_c}(T')-\mu_{\alpha_c}(T)\le -\delta$ 
for all subtriples $T'\subset T$. Thus 
$\mu_{\alpha_c^{\pm}}(T')-\mu_{\alpha_c^{\pm}}(T)<0$. That is, $t$ 
is in $\mathcal{N}^s_{\alpha^{\pm}}$, and hence 
$\mathcal{N}^s_{\alpha} \subseteq
\mathcal{N}^s_{\alpha^{\pm}}-\mathcal{S}_{\alpha^{\pm}}$. 
\end{proof}

\noindent Our goal is to show that the flip loci 
$\mathcal{S}_{\alpha_c^{\pm}}$ are contained in subvarieties of 
positive codimension in $\mathcal{N}^s_{\alpha_c^{\pm}}$ 
respectively. 

\begin{proposition}\label{prop:vicente}
Let $\alpha_c\in (\alpha_m,\alpha_M)$  be a critical value for 
triples of type 
$(\mathbf{n},\mathbf{d})=(n_1,n_2,d_1,d_2)$. 
Let $T=(E_1,E_2,\phi)$ be a triple of this type. 
\begin{enumerate}
\item Suppose that $T$ represents a point in 
$\mathcal{S}_{\alpha_c^+}$, i.e.\ suppose that $T$ is
$\alpha_c^+$-stable but $\alpha_c^-$-unstable.  Then $T$ has a 
description as the middle term in an extension 
 \begin{equation}\label{destab}
 0\to T'\to T\to T'' \to 0
 \end{equation}
in which 
 \begin{enumerate}
 \item  $T'$ and $T''$ are both $\alpha_c^+$-stable, with 
 $\mu_{\alpha_c^+}(T')<\mu_{\alpha_c^+}(T)$,
 \item $T'$ and $T''$ are both $\alpha_c$-semistable with
 $\mu_{\alpha_c}(T')=\mu_{\alpha_c}(T)$,
\item $\frac{n_2'}{n_1'+n_2'}$ is a maximum among all proper subtriples
 $T'\subset T$  which satisfy {\rm (b)},
 \item $n_1'+n_2'$ is a minimum among all subtriples which satisfy {\rm (c)}.
 \end{enumerate}
 
\item Similarly, if $T$ represents a point in 
$\mathcal{S}_{\alpha_c^-}$, i.e.\ if $T$ is $\alpha_c^-$-stable but
$\alpha_c^+$-unstable,  then $T$ has a 
description as the middle term in an extension 
$\mathrm{(\ref{destab})}$
in which 
 \begin{enumerate}
 \item  $T'$ and $T''$ are both $\alpha_c^-$-stable with 
 $\mu_{\alpha_c^-}(T')<\mu_{\alpha_c^-}(T)$,
 \item $T'$ and $T''$ are both $\alpha_c$-semistable with
 $\mu_{\alpha_c}(T')=\mu_{\alpha_c}(T)$,
 \item $\frac{n_2'}{n_1'+n_2'}$ is a minimum among all proper subtriples
 $T'\subset T$  which satisfy {\rm (b)},
 \item $n_1'+n_2'$ is a minimum among all subtriples which satisfy {\rm (c)}.
 \end{enumerate}
\end{enumerate}
\end{proposition}
{\em Proof.\/} In both cases (i.e.\ (1) and (2)), since its stability 
property changes at 
$\alpha_c$, the triple $T$ must be strictly $\alpha_c$-semistable, i.e.\ 
it must have a proper subtriple $T'$ with 
$\mu_{\alpha_c}(T')=\mu_{\alpha_c}(T).$  We can thus consider the (non-empty) set
 $$
 {\mathcal F}_1=\{T'\subsetneq T\ |\
 \mu_{\alpha_c}(T')= \mu_{\alpha_c}(T)\ \}.
 $$
\noindent ({\em Proof of (1)})  
Suppose first that $T$ is 
$\alpha_c^+$-stable but $\alpha_c^-$-unstable. We observe that if
$T'\in {\mathcal F}_1$, then $\frac{n_2'}{n_1'+n'_2}< \frac{n_2}{n_1+n_2}$, since
otherwise $T$ could not be $\alpha_c^+$-stable. But the allowed 
values for 
$\frac{n'_2}{n_1'+n'_2}$ are limited by the constraints 
$0\le n_1'\le n_1$, $0\le n'_2\le n_2$ and $n_1'+n'_2\ne 0$. We can thus define
 $$
  \lambda_0=\max\left\{\frac{n'_2}{n_1'+n'_2}\ \biggm|\ T'\in {\mathcal F}_1\ \right\}\
 $$
and set 
 $$
 {\mathcal F}_2=\left\{T_1\subset {\mathcal F}_1\ \biggm|\\
 \frac{n'_2}{n_1'+n'_2}=\lambda_0\ \right\}.
 $$
\noindent Now let $T'$ be any triple in $\mathcal F_2$. Since $T'$ has maximal 
$\alpha_c$-slope, we can assume that $T''=T/T'$ is a locally free 
triple, i.e.\ if 
$T''=(E_2'',E_1'',\Phi)$ then $E_2''$ and $E_1''$ are both locally 
free.  Furthermore, since $T$ is 
$\alpha_c$-semistable and $\mu_{\alpha_c}(T')=
\mu_{\alpha_c}(T)=\mu_{\alpha_c}(T'')$, it follows that both
$T'$ and $T''$ are $\alpha_c$-semistable and of the same $\alpha_c$-slope.
We now show that $T''$ is $\alpha_c^+$-stable. Suppose not. Then 
there is a proper subtriple $\tilde{T}''\subset T''$ with 
$\mu_{\alpha_c^+}(\tilde{T}'')\ge\mu_{\alpha_c^+}(T'')$. 
However, since we can assume that $\alpha_c^+$ is not a critical 
value for triples of type $(\tilde{T}'')$, we must have 
$$\mu_{\alpha_c^+}(\tilde{T}'')>\mu_{\alpha_c^+}(T'')\ .$$  
\noindent Thus, since $(T'')$ is $\alpha_c$-semistable, we must have
\begin{itemize}
\item $\mu_{\alpha_c}(\tilde{T}'')\le\mu_{\alpha_c}(T'')$ 
\item $\frac{\tilde{n}_2''}{\tilde{n}_1''+\tilde{n}_2''}>\frac{n''_2}{n_1''+n''_2}$
\end{itemize}
\noindent  If $\mu_{\alpha_c}(\tilde{T}'')<\mu_{\alpha_c}(T'')$, say 
$\mu_{\alpha_c}(\tilde{T}'')=\mu_{\alpha_c}(T'')-\delta$, then in order to have
$\mu_{\alpha_c^+}(\tilde{T}'')>\mu_{\alpha_c^+}(T'')$ we must have
$$\frac{\tilde{n}_2''}{\tilde{n}_1''+\tilde{n}_2''}>
\frac{n''_2}{n_1''+n''_2}+\frac{\delta}{\epsilon}\ .$$
\noindent Letting $\epsilon$ approach zero, we see that 
$\frac{\tilde{n}_2''}{\tilde{n}_1''+\tilde{n}_2''}$ must be arbitrarily large. 
This cannot be if $0\le \tilde{n}_1''\le n_1''$ and 
$0\le \tilde{n}_2''\le n''_2$ (and $\tilde{n}_1''+\tilde{n}_2''>0$).  
We may thus assume that 
$\mu_{\alpha_c}(\tilde{T}'')=\mu_{\alpha_c}(T'').$ 
Consider now the subtriple $\tilde{T}'\subset T$ defined by the 
pull-back diagram 
 $$
 0\to T'\to \tilde{T}'\to \tilde{T}'' \to 0 \ .
 $$
\noindent This has $\mu_{\alpha_c}(\tilde{T}')=\mu_{\alpha_c}(T'')=\mu_{\alpha_c}(T)$ and 
thus 
$$\frac{\tilde{n}_2'}{\tilde{n}_1'+\tilde{n}_2'}\le\lambda_0=
\frac{n'_2}{n_1'+n'_2}\ .$$
\noindent It follows from this and the above extension that
$$\frac{\tilde{n}_2''}{\tilde{n}_1''+\tilde{n}_2''}\le\lambda_0=
\frac{n'_2}{n_1'+n'_2}\ .$$
\noindent However, since $\mu_{\alpha_c}(T')=\mu_{\alpha_c}(T)$ but 
$\mu_{\alpha_c^+}(T')<\mu_{\alpha_c^+}(T)$, we have that
$$\frac{n'_2}{n_1'+n'_2}<\frac{n''_2}{n_1''+n''_2}\ .$$
\noindent Combining the previous two inequalities we get
$$\frac{\tilde{n}_2''}{\tilde{n}_1''+\tilde{n}_2''}<\frac{n''_2}{n_1''+n''_2}\ $$
\noindent which is a contradiction. 
Now take $T'\in\mathcal F_2$ with \it minimum rank (i.e.\ minimum 
$n_1'+n'_2$) in\rm\ $\mathcal F_2$. We claim that $T'$ is 
$\alpha_c^+$-stable. If not, then as before
it has a proper subtriple $\tilde{T}'$ with 
\begin{itemize}
\item $\mu_{\alpha_c}(\tilde{T}')\le\mu_{\alpha_c}(T')$
\item $\frac{\tilde{n}_2'}{\tilde{n}_1'+\tilde{n}_2'}>\frac{n'_2}{n_1'+n'_2}$.
\end{itemize}
\noindent But then $\tilde{n}_1'+\tilde{n}_2'<n_1'+n'_2$, which contradicts 
the minimality of $n_1'+n'_2$. Thus $T'$ is $\alpha_c^+$-stable. 
Moreover, since $T$ is $\alpha_c^+$-stable it follows that 
$\mu_{\alpha_c^+}(T')<\mu_{\alpha_c^+}(T)$.
Thus taking $T'\in\mathcal F_2$ with minimum rank, and $T''=T/T'$, we 
get a description of $T$ as an extension in which (a)-(d) are 
satisfied. 
\bigskip

\noindent {\em Proof of (2).}
If 
$T$ is $\alpha_c^-$-stable but 
$\alpha_c^+$-unstable, then $\frac{n'_2}{n_1'+n'_2}> \frac{n_2}{n_1+n_2}$ for all
$T'\in\mathcal F_1$. The proof of (i) must thus be modified as
follows. With 
 $$
 \lambda_0=\min\left\{\frac{n'_2}{n_1'+n'_2}\ \biggm|\ T'\in {\mathcal F}_1\ \right\}\
 $$
\noindent we can define 
 $$
 {\mathcal F}_2=\left\{T'\subset {\mathcal F}_1\ \biggm|\\
 \frac{n'_2}{n_1'+n'_2}=\lambda_0\ \right\}
 $$
\noindent and select $T'\in\mathcal F_2$ such that $T'$ has minimal 
rank in 
$\mathcal F_2$. It follows in a similar fashion to that above that
$T$ has a description as
 $$
 0\to T'\to T\to T'' \to 0
 $$
in which all the requirements of the proposition are satisfied. 
\qed
\begin{remark}Unlike for Jordan-Holder filtrations for semistable 
objects, the filtrations produced by the above proposition are always 
of length two, i.e.\ always yield a description of the semistable 
object as an extension of stable objects. This is achieved by 
exploiting the extra `degree of freedom' provided by the parameter 
$\alpha_c$. The true advantage of never having to consider extensions 
of length greater than two is that it removes the need for inductive 
procedures in the analysis of the flip loci. 
\end{remark}
\begin{definition}
\label{def:S-plusminus} 
Let $\alpha_c\in (\alpha_m,\alpha_M)$ 
be a critical value for triples of type 
$(\mathbf{n},\mathbf{d})$.  Let
$(\mathbf{n}',\mathbf{d}') = (n'_1,n'_2,d'_1,d'_2)$ and
$(\mathbf{n}'',\mathbf{d}'') = (n''_1,n''_2,d''_1,d''_2)$ be such that
\begin{equation}\label{eqn:C1}
(\mathbf{n},\mathbf{d})=(\mathbf{n}',\mathbf{d}')+ 
(\mathbf{n}'',\mathbf{d}'')\ , 
\end{equation}
(i.e.\ $n_1 = n'_1 + n''_1$, $n_2 = n'_2 + n''_2$, $d_1 = d'_1 + 
d''_1$, and $d_2 = d'_2 + d''_2$), and also 
\begin{equation}\label{eqn:C2}
 \frac{d'_1+d'_2}{n'_1+n'_2}+\alpha_c\frac{n'_2}{n'_1+n'_2}=
\frac{d''_1+d''_2}{n''_1+n''_2}+\alpha_c\frac{n''_2}{n''_1+n''_2}.
\end{equation}
\begin{enumerate}
\item Define 
$\tilde{\mathcal{S}}_{\alpha_c^+}(\mathbf{n}'',\mathbf{d}'',\mathbf{n}',
\mathbf{d}')$ to be the set of all isomorphism classes of extensions 
\begin{displaymath}
  0 \lto T' \lto T \lto T'' \lto 0,
\end{displaymath}
where $T'$ and $T''$ are $\alpha_c^+$-stable triples with topological 
invariants $(\mathbf{n}',\mathbf{d}')$ and 
$(\mathbf{n}'',\mathbf{d}'')$ respectively, and the isomorphism is
on the triple $T$. 
\item Define 
$\tilde{\mathcal{S}}^{0}_{\alpha_c^+}(\mathbf{n}'',\mathbf{d}''
,\mathbf{n}',\mathbf{d}') 
\subset \tilde{\mathcal{S}}_{\alpha_c^+}(\mathbf{n}'',\mathbf{d}'',
\mathbf{n}',\mathbf{d}')$
to be the set of all extensions for which moreover  $T$ is 
$\alpha_c^+$-stable.  In an analogous manner, define 
$\tilde{\mathcal{S}}_{\alpha_c^-}(\mathbf{n}'',\mathbf{d}'', 
\mathbf{n}',\mathbf{d}')$ and
$\tilde{\mathcal{S}}^{0}_{\alpha_c^-}(\mathbf{n}'',\mathbf{d}'' 
,\mathbf{n}',\mathbf{d}') 
\subset \tilde{\mathcal{S}}_{\alpha_c^+}(\mathbf{n}'',\mathbf{d}'',
\mathbf{n}',\mathbf{d}')$.
\item  Define
\begin{align*}
\tilde{\mathcal{S}}_{\alpha_c^+}= \bigcup
\tilde{\mathcal{S}}_{\alpha_c^+}(\mathbf{n}'',\mathbf{d}'',
\mathbf{n}',\mathbf{d}')\ &,\ 
\tilde{\mathcal{S}}^0_{\alpha_c^+} = \bigcup
\tilde{\mathcal{S}}^0_{\alpha_c^+}(\mathbf{n}'',\mathbf{d}'',
\mathbf{n}',\mathbf{d}')\\
\end{align*}
\noindent where the union is over all $(n'_1,n'_2,d'_1,d'_2)$ and
$(n''_1,n''_2,d''_1,d''_2)$ such that the above conditions apply, and 
also $\frac{n'_2}{n'_1+n'_2}<\frac{n''_2}{n''_1+n''_2}$. 
\item  Similarly, define 
\begin{align*}
\tilde{\mathcal{S}}_{\alpha_c^-}= \bigcup
\tilde{\mathcal{S}}_{\alpha_c^-}(\mathbf{n}'',\mathbf{d}'',\mathbf{n}',\mathbf{d}')\ 
&,\ 
\tilde{\mathcal{S}}^0_{\alpha_c^-}= \bigcup
\tilde{\mathcal{S}}^0_{\alpha_c^-}(\mathbf{n}'',\mathbf{d}'',\mathbf{n}',\mathbf{d}')
\end{align*}
\noindent where the union is over all $(n'_1,n'_2,d'_1,d'_2)$ and
$(n''_1,n''_2,d''_1,d''_2)$ such that the above conditions apply, and 
also $\frac{n'_2}{n'_1+n'_2}>\frac{n''_2}{n''_1+n''_2}$. 
\end{enumerate}
\end{definition}
\begin{remark}It can happen that $\tilde{\mathcal{S}}^0_{\alpha_c^-}$ 
or $\tilde{\mathcal{S}}^0_{\alpha_c^-}$ is empty. For instance there 
may be no possible choices of  $(n'_1,n'_2,d'_1,d'_2)$ and 
$(n''_1,n''_2,d''_1,d''_2)$ which satisfy all the required 
conditions. In this case, the implication of the next lemma is that 
one or both of the flip loci $\mathcal{S}_{\alpha_c^{\pm}}$ is empty. 
\end{remark}
\begin{lemma}\label{lemma:vmaps}
There are maps, say 
$v^{\pm}:\tilde{\mathcal{S}}^0_{\alpha_c^{\pm}}
\longrightarrow\mathcal{N}^s_{\alpha_c^{\pm}}$, 
which map triples to their equivalence classes. The images contain 
the flip loci $\mathcal{S}_{\alpha_c^{\pm}}$. 
\end{lemma}
\begin{proof}The existence of the maps is clear. The second statement,
about the images of the maps, follows by Proposition 
\ref{prop:vicente}. Indeed, suppose that $T$ represents a point in 
$\mathcal{S}_{\alpha_c^{+}}$ and that
 $$
 0\to T'\to T\to T'' \to 0
 $$
is an extension of the type described in proposition 
\ref{prop:vicente}, with 
$T'$ a triple of type $(\mathbf{n}',\mathbf{d}')$  and $T''$
a triple of type $(\mathbf{n}'',\mathbf{d}'')$. Then  
$(\mathbf{n}',\mathbf{d}')$ and $(\mathbf{n}',\mathbf{d}')$ satisfy
conditions (\ref{eqn:C1}) and (\ref{eqn:C2}). Furthermore, since 
$\mu_{\alpha^+_c}(T')< \mu_{\alpha^+_c}(T'')$, we must have
$\frac{n'_2}{n'_1+n'_2}<\frac{n''_2}{n''_1+n''_2}$. Thus $T$ 
is contained in  $v^+(\tilde{\mathcal{S}}^0_{\alpha_c^+})$. A similar 
argument shows that $\mathcal{S}_{\alpha_c^{+}}$ is contained in  
$v^-(\tilde{\mathcal{S}}^0_{\alpha_c^-})$
\end{proof}
%%%%%%%%%%%%%%%%%%%%%%%
\subsection{Codimension Estimates}\label{subs:codim}
%%%%%%%%%%%%%%%%%%%%%%%%
Let $\alpha_c\in (\alpha_m,\alpha_M)$  be a critical value for 
triples of type $(\mathbf{n},\mathbf{d})$.  Fix 
$(\mathbf{n}',\mathbf{d}') = (n'_1,n'_2,d'_1,d'_2)$ and 
$(\mathbf{n}'',\mathbf{d}'') = (n''_1,n''_2,d''_1,d''_2)$ as in 
Definition \ref{def:S-plusminus}. For simplicity we shall denote the 
moduli spaces of $\alpha_c^{\pm}$-semistable triples of topological 
type 
$(\mathbf{n}',\mathbf{d}')$, respectively 
$(\mathbf{n}'',\mathbf{d}'')$, by 
\begin{displaymath}
  \mathcal{N}_{\alpha_c^{\pm}}' =
  \mathcal{N}_{\alpha_c^{\pm}}(\mathbf{n}',\mathbf{d}')
  \quad\text{and}\quad
  \mathcal{N}_{\alpha_c^{\pm}}'' =
  \mathcal{N}_{\alpha_c^{\pm}}( \mathbf{n}'',\mathbf{d}'').
\end{displaymath}
\begin{proposition}\label{prop:Sprojbundle} If $\alpha_c\ge 2g-2$ then the space
$\tilde{\mathcal{S}}_{\alpha_c^\pm}(\mathbf{n}'',\mathbf{d}'',\mathbf{n}',\mathbf{d}')$
is a locally trivial fibration over $\mathcal{N}_{\alpha_c^{\pm}}' 
\times\mathcal{N}_{\alpha_c^{\pm}}''$, with projective fibers of dimension 
$$-\chi(\mathbf{n}'',\mathbf{d}'',\mathbf{n}',\mathbf{d}')-1\ .
$$
In particular, 
$\tilde{\mathcal{S}}_{\alpha_c^\pm}
(\mathbf{n}'',\mathbf{d}'',\mathbf{n}',\mathbf{d}')$ has dimension 
$$
1-\chi(\mathbf{n}',\mathbf{d}',\mathbf{n}',\mathbf{d}') 
 -\chi(\mathbf{n}'',\mathbf{d}'',\mathbf{n}'',\mathbf{d}'')
 -\chi(\mathbf{n}'',\mathbf{d}'',\mathbf{n}',\mathbf{d}')\ ,
$$
\noindent where 
$\chi(\mathbf{n}',\mathbf{d}',\mathbf{n}',\mathbf{d}')$ etc. are 
as in section \ref{sec:extensions-of-triples}. 
\end{proposition}
\begin{proof} 
{}From the defining properties of $\tilde{\mathcal{S}}_{\alpha_c^\pm} 
(\mathbf{n}'',\mathbf{d}'',\mathbf{n}',\mathbf{d}')$ there is map 
\begin{equation}\label{eqn:fibration}
\tilde{\mathcal{S}}_{\alpha_c^\pm}
(\mathbf{n}'',\mathbf{d}'',\mathbf{n}',\mathbf{d}')\longrightarrow 
 \mathcal{N}_{\alpha_c^{\pm}}' \times\mathcal{N}_{\alpha_c^{\pm}}''
\end{equation}
\noindent which sends an extension 
  $$
 0\to T'\to T\to T'' \to 0
 $$
 \noindent to the pair $([T'],[T''])$, where $[T']$ denotes the class 
 represented by $T'$ and similarly for $[T'']$. We first examine the 
 fibers of this map.  
 
 Notice that $T'$ and $T''$
 satisfy the hypothesis of Proposition \ref{prop:h2-vanishing}
 and therefore of Corollary \ref{dimension-ext1}. Notice moreover that, since 
 $\mu_{\alpha_c^{\pm}}(T')<\mu_{\alpha_c^{\pm}}(T'')$, it is not possible
 to have $T'\cong T''$. Thus  (cf.\ Corollary \ref{dimension-ext1} and 
 Proposition \ref{prop:h0-vanishing}(2)) we 
 get
 \begin{align}\label{eq:fiberdim}
 \mathrm{dim}\mathbb{P}(\mathrm{Ext}^1(T'',T'))&=
 \mathrm{dimExt}^1(T'',T')-1 \notag\\
 &=-\chi(T'',T')-1 \notag\\
 &= -\chi(\mathbf{n}'',\mathbf{d}'',\mathbf{n}',\mathbf{d}')-1\ ,
 \end{align}
 \noindent which is independent of $T'$ and $T''$.
It remains to establish that the fibration (\ref{eqn:fibration}) is 
locally trivial.  If the coprime conditions 
$\GCD(n_1',n_2',d_1'+d_2')=1=\GCD(n_1'',n_2'',d_1''+d_2'')$ hold then the moduli
spaces $ \mathcal{N}_{\alpha_c^{\pm}}'$ and 
$\mathcal{N}_{\alpha_c^{\pm}}''$ are fine moduli spaces (cf.\ 
\cite{Schmitt}). That is, there are universal objects, say $\mathcal{U}'$ and 
$\mathcal{U}''$, defined over
$\mathcal{N}_{\alpha_c^{\pm}}'\times X$ and 
$\mathcal{N}_{\alpha_c^{\pm}}''\times X$. These can be viewed as
coherent sheaves of algebras (cf.\ \cite{AG2}), or more precisely as 
examples of the $Q$-bundles considered in \cite{gothen-king:2002}. 
Pulling these back to 
$\mathcal{N}_{\alpha_c^{\pm}}'\times \mathcal{N}_{\alpha_c^{\pm}}''\times X$
we can construct $Hom(\mathcal{U}'',\mathcal{U}')$ (where we have 
abused notation for the sake of clarity). Taking the projection from 
$\mathcal{N}_{\alpha_c^{\pm}}'\times 
\mathcal{N}_{\alpha_c^{\pm}}''\times X$ onto
$\mathcal{N}_{\alpha_c^{\pm}}'\times \mathcal{N}_{\alpha_c^{\pm}}''$, 
we can then construct the first direct image sheaf. By the results in 
\cite{gothen-king:2002}, we can identify the fibers as hypercohomology groups 
which, in this case, parameterize extensions of triples. We thus 
obtain 
$\tilde{\mathcal{S}}_{\alpha_c^\pm}$ as the projectivization of the 
first direct image of 
$Hom(\mathcal{U}'',\mathcal{U}')$. 
If the coprime conditions fail, then the universal objects do not 
exist globally. However they still exist locally over (analytic) open 
sets in the stable locus in the base 
$\mathcal{N}_{\alpha_c^{\pm}}'\times 
\mathcal{N}_{\alpha_c^{\pm}}''$. This is sufficient for our purpose since
by construction the image of the map in \eqref{eqn:fibration} lies in 
the stable locus. The result now follows from (\ref{eq:fiberdim}) and 
formula (\ref{eq:dim-triples}) (in Theorem 
\ref{thm:smoothdim}) as applied to 
 $\mathcal{N}_{\alpha_c^{\pm}}'$ and 
 $\mathcal{N}_{\alpha_c^{\pm}}''$.
\end{proof}
\begin{corollary}\label{cor:codim-est} If $\alpha_c> 2g-2$ then the loci 
$\mathcal{S}_{\alpha_c^{\pm}}\subset\mathcal{N}^s_{\alpha_c^{\pm}}$ are 
locally contained in subvarieties of  codimension bounded below by 
$$
\mathrm{min}\{-\chi(\mathbf{n}'',\mathbf{d}'',\mathbf{n}',\mathbf{d}')\}\ ,
$$ 
\noindent where the minimum is over all $(\mathbf{n}',\mathbf{d}')$ and 
$(\mathbf{n}'',\mathbf{d}'')$ which satisfy
(\ref{eqn:C1}) and (\ref{eqn:C2}) and also 
$\frac{n'_2}{n'_1+n'_2}<\frac{n''_2}{n''_1+n''_2}$ (in the case of
$\mathcal{S}_{\alpha_c^{+}}$) or
$\frac{n'_2}{n'_1+n'_2}>\frac{n''_2}{n''_1+n''_2}$
(in the case of 
$\mathcal{S}_{\alpha_c^{-}}$). The same is true for 
$\mathcal{S}_{\alpha_c^{+}}$ when $\alpha_c= 2g-2$.
\end{corollary}
\begin{proof} If $\alpha_c > 2g-2$ then we can assume $\alpha_c^{\pm}\ge 2g-2$. 
Clearly also, $\alpha_c^{+}\ge 2g-2$ when $\alpha_c=2g-2$.  Thus by 
Theorem \ref{thm:smoothdim} the moduli spaces 
$\mathcal{N}^s_{\alpha^{\pm}}$ are smooth and have dimension 
$1-\chi(n_1,n_2,d_1,d_2)$.  By Corollary \ref{cor:chi-relation} and 
Proposition \ref{prop:Sprojbundle} we get 
\begin{align*}
\mathrm{dim}\mathcal{N}^s_{\alpha^{\pm}}=1 &-\chi(n'_1,d'_1,n'_2,d'_2)
-\chi(n''_1,d''_1,n''_2,d''_2)\\ 
&-\chi(n'_1,d'_1,n''_2,d''_2)
-\chi(n''_1,d''_1,n'_2,d'_2)\\ =&\mathrm{dim} 
\tilde{\mathcal{S}}_{\alpha_c^\pm}(\mathbf{n}'',\mathbf{d}'',
\mathbf{n}',\mathbf{d}')-\chi(n''_1,d''_1,n'_2,d'_2).
\end{align*}
\end{proof}
In order to complete the estimate on the size of the flip loci we 
thus need to estimate the expressions $\chi(T'',T')= 
\chi(\mathbf{n}'',\mathbf{d}'',\mathbf{n}',\mathbf{d}')$. 
The basic idea is to identify $\chi(T',T'')$ as a hypercohomology 
Euler characteristic for the complex 
$C^{\bullet}(T'',T')$ defined in \ref{eq:extension-complex}
and to notice that the complex is itself a holomorphic triple. As 
such it ought to satisfy a stability condition induced from the 
stability condition of $T'$ and 
$T''$.  The right way to obtain the stability condition for
$C^{\bullet}(T'',T')$ should be via the Hitchin--Kobayashi
correspondence, in a way analogous to the case of vector bundles: to 
show that the tensor product of two stable vector bundles is 
semistable, one constructs a flat unitary connection on the tensor 
product from flat unitary connections on each of the factors.  In the 
case of holomorphic triples one considers solutions to the so-called 
coupled vortex equations instead of flat connections.  However, it 
turns out that one cannot construct (at least in a simple way) a 
solution to the coupled vortex equations for $C^{\bullet}(T'',T')$ 
from solutions on $T'$ and $T''$.  Instead one needs to consider a 
so-called \emph{holomorphic chain}.  These objects were studied in 
\cite{alvarez-garcia-prada:2001} and generalize holomorphic triples.
Next we recall the relevant results and definitions. 
%%%%%%%%%%%%%%%%%%%%%%%%%%
\subsection{Holomorphic chains}\label{sec:holomorphic-chains}
%%%%%%%%%%%%%%%%%%%%%%%%%%
A \emph{holomorphic chain} is a diagram 
\begin{displaymath}
  \mathcal{C} \colon E_{m} \overset{\phi_m}{\lto} E_{m-1}
  \overset{\phi_{m-1}}{\lto} \cdots
  \overset{\phi_{1}}{\lto} E_{0},
\end{displaymath}
where each $E_i$ is a holomorphic bundle and $\phi_{i} \colon E_{i} \to 
E_{i-1}$ is a holomorphic map.  Let 
\begin{align*}
  \mu(\mathcal{C}) &= \mu(E_{0} \oplus \cdots \oplus E_{m}), \\
  \lambda_{i}(\mathcal{C}) &=   \frac{\rk(E_{i})}
       {\sum_{i=0}^{m}\rk(E_{i})}, \quad i=0,\ldots, m.
\end{align*}
For $\balpha = (\alpha_1,\ldots,\alpha_m) \in \R^{m}$, the 
\emph{$\balpha$-slope} of $\mathcal{C}$ is defined to be
\begin{displaymath}
  \mu_{\balpha}(\mathcal{C}) = \mu(\mathcal{C}) +
  \sum_{i=1}^{m}\alpha_{i}\lambda_{i}(\mathcal{C}).
\end{displaymath}
The notions of $\balpha$-(semi)stability and $\balpha$-polystability 
are then defined via the standard $\balpha$-slope condition on 
sub-chains. Note that holomorphic chains generalize holomorphic 
triples: a holomorphic triple is a holomorphic chain of length 2, and 
the stability notions coincide, taking $\balpha=(\alpha)$. 
There are natural gauge theoretic equations for holomorphic chains, 
which we now describe.  Define $\btau = (\tau_{0}, \ldots, \tau_{m}) 
\in \R^{m+1}$ by
\begin{equation}
  \label{eq:tau-alpha}
  \tau_{i} = \mu_{\alpha}(\mathcal{C}) - \alpha_{i}, \quad i=0,\ldots, m,
\end{equation}
where we make the convention $\alpha_0 = 0$.  Then $\balpha$ can be 
recovered from $\btau$ by 
\begin{equation}
  \label{eq:alpha-tau}
  \alpha_{i} = \tau_{0} - \tau_{i}, \quad i=0,\ldots, m.
\end{equation}
The 
\emph{$\btau$-vortex equations} or \emph{chain vortex equations}
\begin{displaymath}
  \sqrt{-1}\Lambda F(E_{i}) + \phi_{i+1} \phi_{i+1}^* - \phi_{i}^*\phi_{i}
  = \tau_{i} \Id_{E_{i}}, \quad i=0,\ldots, m,
\end{displaymath}
are equations for Hermitian metrics on $E_{0}, \ldots, E_{m}$.  Here 
$F(E_{i})$ is the curvature of the Hermitian connection on $E_{i}$,
$\Lambda$ is contraction with the K\"ahler form and $\vol(X) = 2\pi$.
By convention $\phi_{i} = 0$ for $i \neq 1, \ldots, m$. 
With these preliminaries we can now state the Hitchin--Kobayashi 
correspondence for holomorphic chains. 
\begin{theorem}[\cite{alvarez-garcia-prada:2001}]
  \label{thm:chain-hitchin-kobayashi}
  A holomorphic chain $\mathcal{C}$ is $\balpha$-polystable if and
  only if the $\btau$-vortex equations have a solution, where
  $\balpha$ and $\btau$ are related by \eqref{eq:tau-alpha}.
\end{theorem}
\begin{remark}
  To be consistent with the previous notation for holomorphic
  triples we use here a slightly different notation from that of
  \cite{alvarez-garcia-prada:2001}: in that paper a parameter vector
  $\tilde{\balpha} =
  (\tilde{\alpha}_0,\tilde{\alpha}_1,\ldots,\tilde{\alpha}_m)$ is used; in
  these terms our $\balpha$ is given by $\balpha = (\tilde{\alpha}_0 -
  \tilde{\alpha}_1,\ldots,\tilde{\alpha}_0 - \tilde{\alpha}_m)$.
\end{remark}
%\subsection{A solution to the chain vortex equations}
%\label{sec:chain-vortex-solution}
Consider the length $3$ holomorphic chain 
\begin{equation}
  \label{eq:H-chain}
  \widetilde{C^{\bullet}}(T'',T') \colon
  {E_{1}''}^{*} \otimes E_{2}'
  \xrightarrow{a_2}
  {E_{1}''}^{*} \otimes E_{1}' \oplus {E_{2}''}^{*} \otimes E_{2}'
  \xrightarrow{a_1}
  {E_{2}''}^{*} \otimes E_{1}' ,
\end{equation}
where 
\begin{align*}
  a_{2}(\psi) &= (\phi'\psi, -\psi\phi''), \\
  a_{1}(\psi_{1},\psi_{2}) &= \phi'\psi_{2} - \psi_{1}\phi''.
\end{align*}
We shall sometimes write this chain briefly as 
\begin{displaymath}
  \widetilde{C^{\bullet}}(T'',T') \colon
  C_{2} \xrightarrow{a_{2}} C_{1} \xrightarrow{a_{1}} C_{0}.
\end{displaymath}
Note that the last two terms of $\widetilde{C^{\bullet}}(T'',T')$ is 
just the complex $C^{\bullet}(T'',T')$.  Note also that 
$\widetilde{C^{\bullet}}(T'',T')$ is not in general a complex.  Our goal
in this section is to prove, using the Hitchin--Kobayashi 
correspondence for chains, that if $T'$ and $T''$ are 
$\alpha$-polystable then $\widetilde{C^{\bullet}}(T'',T')$ is
$\balpha$-polystable for a suitable choice of $\balpha$.
\begin{lemma}
  \label{lemma:chain-vortex-solution}
  Let $T'$ and $T''$ be holomorphic triples and suppose we have
  solutions to the $(\tau_{1}',\tau_{2}')$-vortex equations on $T'$
  and the $(\tau_{1}'',\tau_{2}'')$-vortex equations on $T''$, such
  that $\tau_{1}' - \tau_{1}''= \tau_{2}' - \tau_{2}''$.
  Then the induced Hermitian metric on $\widetilde{C^{\bullet}}(T'',T')$
  satisfies the chain vortex equations
  \begin{align}
    \label{eq:chain-vortex-1}
    \sqrt{-1}\Lambda F(C_{0}) + a_{1}{a_{1}}^{*}
      &= \tilde{\tau}_{0} \Id_{C_{0}}, \\
    \label{eq:chain-vortex-2}
    \sqrt{-1}\Lambda F(C_{1}) + a_{2}{a_{2}}^{*} - {a_{1}}^{*}a_{1}
      &= \tilde{\tau}_{1} \Id_{C_{1}}, \\
    \label{eq:chain-vortex-3}
    \sqrt{-1}\Lambda F(C_{2}) - {a_{2}}^{*}a_{2}
      &= \tilde\tau_{2} \Id_{C_{2}}.
  \end{align}
  for $\btau = (\tilde{\tau}_{0},\tilde{\tau}_{1},\tilde{\tau}_{2})$
  given by
  \begin{align*}
    \tilde{\tau}_{0} &= \tau_{1}' - \tau_{2}'', \\
    \tilde{\tau}_{1} &= \tau_{1}' - \tau_{1}''
                      = \tau_{2}' - \tau_{2}'', \\
    \tilde{\tau}_{2} &= \tau_{2}' - \tau_{1}''. \\
  \end{align*}
\end{lemma}
\begin{proof}
  We shall only show that the induced Hermitian metric satisfies
  \eqref{eq:chain-vortex-2}, since the proofs that it satisfies the
  two remaining equations are similar (but simpler).
  
  The coupled vortex equations for $T'$ and $T''$ are
  \begin{align*}
    \sqrt{-1}\Lambda F(E_{1}')
      + \phi'{\phi'}^{*}&=\tau_{1}' \Id_{E_{1}'}, &
      \sqrt{-1}\Lambda F(E_{1}'')
      + \phi''{\phi''}^{*}&=\tau_{1}'' \Id_{E_{1}''}, \\
      \sqrt{-1}\Lambda F(E_{2}')
      - {\phi'}^{*}\phi'&=\tau_{2}' \Id_{E_{2}'}, &
      \sqrt{-1}\Lambda F(E_{2}'')
      - {\phi''}^{*}\phi''&= \tau_{2}'' \Id_{E_{2}''}. \\
  \end{align*}
  We shall write the left hand side of \eqref{eq:chain-vortex-2} in
  terms of these known data of the triples $T'$ and $T''$.
  First, we note that
  \begin{displaymath}
    F({E_{i}'}^{*}) = -F(E_{i}')^{t}, \quad i=1,2,
  \end{displaymath}
  and similarly for $F({E_{i}''}^{*})$.  Hence
  \begin{align}
    F(C_{1}) &= F({E_{1}''}^{*} \otimes E_{1}'
      \oplus {E_{2}''}^{*} \otimes E_{2}') \notag \\
    &= \bigl(
      F({E_{1}''}^{*}) \otimes \Id + \Id \otimes F(E_{1}'),\
      F({E_{2}''}^{*}) \otimes \Id + \Id \otimes F(E_{2}')
      \bigr) \notag \\
    &= \bigl(
      -F({E_{1}''})^{t} \otimes \Id + \Id \otimes F(E_{1}'),\
      -F({E_{2}''})^{t} \otimes \Id + \Id \otimes F(E_{2}')
      \bigr). \label{eq:F(C_1)}
  \end{align}
  Next we calculate $a_{1}^{*}$: note that for $\xi \otimes x \in
  C_{0}$ and $(\eta_1 \otimes y_1, \eta_2 \otimes y_2) \in C_{1}$ we
  have
  \begin{align*}
    &\bigl\langle a_{1}^{*}(\xi \otimes x),\
    (\eta_1 \otimes y_1, \eta_2 \otimes y_2) \bigr\rangle_{C_{1}} \\
    =&\bigl\langle \xi \otimes x,\
    a_{1}(\eta_1 \otimes y_1,\eta_2 \otimes y_2) \bigr\rangle_{C_{0}} \\
    =&\bigl\langle \xi \otimes x,\ - \eta_1\phi'' \otimes y_1
    + \eta_2 \otimes \phi'(y_2) \bigr\rangle_{C_{0}} \\
    =&\bigl\langle \xi \otimes x,\ - {\phi''}^{t}(\eta_1)
    \otimes y_1 + \eta_2 \otimes \phi'(y_2) \bigr\rangle_{C_{0}} \\
    =& -\bigl\langle \xi, {\phi''}^{t}(\eta_1) \bigr\rangle_{{E_{2}''}^{*}}
    \ \bigl\langle x, y_1 \bigr\rangle_{E_{1}'} 
    \ +\ \bigl\langle \xi,\eta_2 \bigr\rangle_{{E_{2}''}^{*}}
    \ \bigl\langle x,\phi'(y_2) \bigr\rangle_{E_{1}'} \\
    =& -\bigl\langle {\phi''}^{t*}(\xi), \eta_1 \bigr\rangle_{{E_{1}''}^{*}}
    \ \bigl\langle x, y_1 \bigr\rangle_{E_{1}'} 
    \ +\ \bigl\langle \xi,\eta_2 \bigr\rangle_{{E_{2}''}^{*}}
    \ \bigl\langle {\phi'}^{*}(x),y_2 \bigr\rangle_{E_{2}'} \\
    =&\bigl\langle (-{\phi''}^{t*}(\xi) \otimes x,
    \xi \otimes {\phi'}^{*}(x)),\
    (\eta_1 \otimes y_1, \eta_2 \otimes y_2) \bigr\rangle_{C_{1}}.
  \end{align*}
  Hence,
  \begin{equation}
    \label{eq:a_1^*}
    a_{1}^{*}(\xi \otimes x) = \bigl(-{\phi''}^{t*}(\xi) \otimes x,\
    \xi \otimes {\phi'}^{*}(x)\bigr).
  \end{equation}
  Similarly, to calculate $a_{2}^{*}$ consider $\xi \otimes x \in
  C_{2}$ and $(\eta_1 \otimes y_1, \eta_2 \otimes y_2) \in C_{1}$,
  then
  \begin{align*}
    &\bigl\langle a_{2}^{*}(\eta_1 \otimes y_1, \eta_2 \otimes y_2),\
    \xi \otimes x \bigr\rangle_{C_{2}} \\
    =& \bigl\langle (\eta_1 \otimes y_1, \eta_2 \otimes y_2),\
    a_{2}(\xi \otimes x) \bigr\rangle_{C_{1}} \\
    =& \bigl\langle (\eta_1 \otimes y_1, \eta_2 \otimes y_2),\
    (\xi \otimes \phi'(x),-\xi\phi'' \otimes x) \bigr\rangle_{C_{1}} \\
    =& \bigl\langle \eta_1 \otimes y_1, \xi \otimes \phi'(x)
    \bigr\rangle_{{E_{1}''}^{*} \otimes E_{1}'}
    \ +\ \bigl\langle \eta_2 \otimes y_2, -{\phi''}^{t}(\xi) \otimes x
    \bigr\rangle_{{E_{2}''}^{*} \otimes E_{2}'} \\
    =& \bigl\langle \eta_1 \otimes {\phi'}^{*}(y_1), \xi \otimes x
    \bigr\rangle_{{E_{1}''}^{*} \otimes E_{2}'}
    \ +\ \bigl\langle -{\phi''}^{t*}(\eta_2) \otimes y_2, \xi \otimes x
    \bigr\rangle_{{E_{1}''}^{*} \otimes E_{2}'}.
  \end{align*}
  Thus,
  \begin{equation}
    \label{eq:a_2^*}
    a_{2}^{*}(\eta_1 \otimes y_1, \eta_2 \otimes y_2)
    = \eta_1 \otimes {\phi'}^{*}(y_1)
    -{\phi''}^{t*}(\eta_2) \otimes y_2.
  \end{equation}
  Using \eqref{eq:a_2^*} and \eqref{eq:a_1^*} we can now calculate for
  $(\eta_1 \otimes y_1, \eta_2 \otimes y_2) \in C_{1}$:
  \begin{multline}
    \label{eq:a_2a_2^*}
    a_{2}a_{2}^{*}(\eta_1 \otimes y_1, \eta_2 \otimes y_2) \\
    = \bigl(\eta_1 \otimes \phi'{\phi'}^{*}(y_1)
    -{\phi''}^{t*}(\eta_2) \otimes \phi'(y_2),\
    -{\phi''}^{t}(\eta_1) \otimes {\phi'}^{*}(y_1)
    +{\phi''}^{t}{\phi''}^{t*}(\eta_2) \otimes y_2 \bigr),
  \end{multline}
  and
  \begin{multline}
    \label{eq:a_1^*a_1}
    a_{1}^{*}a_{1}(\eta_1 \otimes y_1, \eta_2 \otimes y_2) \\
    =\bigl({\phi''}^{t*}{\phi''}^{t}(\eta_1) \otimes y_1
    -{\phi''}^{t*}(\eta_2) \otimes \phi'(y_2),\
    -{\phi''}^{t}(\eta_1) \otimes {\phi'}^{*}(y_1)
    +\eta_2 \otimes {\phi'}^{*}\phi'(y_2) \bigr).
  \end{multline}
  Putting together \eqref{eq:F(C_1)}, \eqref{eq:a_2a_2^*} and
  \eqref{eq:a_1^*a_1} we finally obtain
  \begin{multline}
    \label{eq:chain-1-lhs}
    \bigl(\sqrt{-1}\Lambda F(C_{1}) + a_{2}{a_{2}}^{*} -
    {a_{1}}^{*}a_{1}\bigr)
    (\eta_1 \otimes y_1, \eta_2 \otimes y_2) \\
    = \Bigl( \eta_1 \otimes \bigl(\sqrt{-1}\Lambda F(E_{1}') +
    \phi'{\phi'}^{*} \bigr)(y_1)
    + \bigl(-\sqrt{-1}\Lambda F(E_{1}'')^{t} +
    {\phi''}^{t*}{\phi''}^{t} \bigr)(\eta_1) \otimes y_1, \\
    \eta_2 \otimes \bigl(\sqrt{-1}\Lambda F(E_{2}')
    - {\phi'}^{*}\phi'\bigr)(y_2)
    - \bigl(\sqrt{-1}\Lambda F(E_{2}'')^{t} -
    {\phi''}^{t}{\phi''}^{t*} \bigr)(\eta_2) \otimes y_2 \Bigr).
  \end{multline}
  Notice that the unpleasant mixed term
  $\bigl(-{\phi''}^{t*}(\eta_2) \otimes \phi'(y_2),
  -{\phi''}^{t}(\eta_1) \otimes {\phi'}^{*}(y_1) \bigr)$ appears both
  in $a_{1}^{*}a_{1}$ and $a_{2}a_{2}^{*}$ and therefore cancels.
  This would not have been the case if we had considered the 
  vortex equations on the triple $C^{\bullet}(T'',T')$ and is the
  reason why we must consider the chain $\widetilde{C^{\bullet}}(T'',T')$.
  Combining \eqref{eq:chain-1-lhs} with the coupled vortex equations
  (or their transposes) for the triples $T'$ and $T''$ we get
  \begin{multline}
    \bigl(\sqrt{-1}\Lambda F(C_{1}) + a_{2}{a_{2}}^{*} -
    {a_{1}}^{*}a_{1}\bigr)
    (\eta_1 \otimes y_1, \eta_2 \otimes y_2) \\
    = \bigl((\tau_{1}'-\tau_{1}'')\eta_1 \otimes y_1,
    (\tau_{2}'-\tau_{2}'')\eta_2 \otimes y_2 \bigr).
  \end{multline}
  Since $\tau_{1}'-\tau_{1}'' = \tau_{2}'-\tau_{2}''$ this
  concludes the proof.
\end{proof}
\begin{proposition}
  \label{prop:H-alpha-stable}
  Let $T'$ and $T''$ be $\alpha$-polystable triples.  Then the
  holomorphic chain $\widetilde{C^{\bullet}}(T'',T')$ is
  $\balpha$-polystable for
  \begin{math}
    \balpha = (\alpha_1,\alpha_2) = (\alpha,2\alpha).
  \end{math}
\end{proposition}
\begin{proof}
  Since $T'$ and $T''$ are $\alpha$-polystable, it follows
  from the Hitchin--Ko\-ba\-yashi correspondence for triples that they
  support solutions to the $(\tau_{1}',\tau_{2}')$- and
  $(\tau_{1}',\tau_{2}')$-vortex equations, respectively, where
  $\alpha = \tau_{1}' - \tau_{2}'= \tau_{1}'' - \tau_{2}''$.  Notice
  that $\tau_{1}' - \tau_{1}''= \tau_{2}' - \tau_{2}''$.  Thus Lemma
  \ref{lemma:chain-vortex-solution} shows that
  $\widetilde{C^{\bullet}}(T'',T')$ supports a solutions to the chain
  vortex equations for $\btau = (\tau_{1}' - \tau_{2}'',\tau_{1}' -
  \tau_{1}'',\tau_{2}' - \tau_{1}'')$.  
  Now the Hitchin--Kobayashi correspondence for chains
  (Theorem~\ref{thm:chain-hitchin-kobayashi}) and
  \eqref{eq:alpha-tau} imply that $\widetilde{C^{\bullet}}(T'',T')$ is
  $\balpha$-polystable for 
  \begin{align*}
    \alpha_1 &= \tau_{1}' - \tau_{2}'' - \tau_{2}' + \tau_{2}'' =
    \alpha, \\
    \alpha_2 &= \tau_{1}' - \tau_{2}'' - \tau_{2}' + \tau_{1}'' =
    2\alpha. 
  \end{align*}
\end{proof}
The following Lemmas will be needed in the next section: 
\begin{lemma}\label{lemma:stability-for-H_*}
  Let $T'=(E'_1,E'_2,\phi')$ and $T''=(E''_1,E''_2,\phi'')$ be triples
  for which the chain $\widetilde{C^{\bullet}}(T'',T')$ is 
  $\balpha=(\alpha,2\alpha)$-poly-stable.  Let
  \begin{align*}
    C_{1} &= {E_{1}''}^{*} \otimes E_{1}' \oplus
    {E_{2}''}^{*} \otimes E_{2}', \\
    C_{0} &= {E_{2}''}^{*} \otimes E_{1}',
  \end{align*}
  and $a_{1} \colon C_{1} \to C_{0}$ be defined as in 
  (\ref{eq:H-chain}).
  Then the following inequalities hold.
  \begin{align}
    \label{eq:ker-a_1}
    \deg(\ker(a_{1})) &\leq \rk(\ker(a_{1}))(\mu_{\alpha}(T') -
    \mu_{\alpha}(T'')), \\
    \label{eq:im-a_1}
    \deg(\im(a_{1})) &\leq \bigl(\rk(C_{0}) - \rk(\im(a_{1}))\bigr)
    (\mu_{\alpha}(T'') - \mu_{\alpha}(T') - \alpha) + \deg(C_{0}).
  \end{align}
\end{lemma}
\begin{proof} If $\rk(\ker(a_{1}))=0$ then (\ref{eq:ker-a_1}) is obvious.
Assume therefore that $\rk(\ker(a_{1})>0$.  Using $\ker(a_{1})$, we 
can then define a quotient of the chain 
  $\widetilde{C^{\bullet}}(T'',T')$ by 
  \begin{displaymath}
      \mathcal{K}\colon  0 \to \ker(a_{1}) \to 0.
  \end{displaymath}
  Thus, since
  $\mu_{\balpha}(\mathcal{K}) = \mu(\ker(a_{1})) + \alpha$, it
  follows that
  \begin{align*}
    \mu(\ker(a_{1})) + \alpha &\leq
    \mu_{\balpha}(\widetilde{C^{\bullet}}(T'',T'))
    = \mu_{\alpha}(T') - \mu_{\alpha}(T'') + \alpha.
  \end{align*}
  We therefore have
  \begin{displaymath}
    \mu(\ker(a_{1})) \leq \mu_{\alpha}(T') - \mu_{\alpha}(T''),
  \end{displaymath}
  which is equivalent to \eqref{eq:ker-a_1}.
The second inequality, i.e.\ (\ref{eq:im-a_1}), is obvious if 
$\rk(\im(a_{1})=\rk(C_0)$. We thus assume $\rk(\im(a_{1})<\rk(C_0)$.
Using the cokernel $\coker(a_{1})$ (or its saturation if it is not 
torsion free), we can define a sub-chain of the chain 
$\widetilde{C^{\bullet}}(T'',T')$ by 
  \begin{displaymath}
      \mathcal{Q}\colon  0 \to 0 \to \coker(a_{1}).
  \end{displaymath}
By the $\balpha$-poly-stability of    
$\widetilde{C^{\bullet}}(T'',T')$ we have  
$\mu_{\balpha}(\mathcal{Q})\ge 
\mu_{\balpha}(\widetilde{C^{\bullet}}(T'',T'))$. This, together with
the fact that 
$$\mu(\coker(a_1)\le\frac{\deg(C_{0})-\deg(\im(a_1))}{\rk(C_0)-\rk(\im(a_1))}\ 
,$$ 
\noindent leads directly to \ref{eq:im-a_1}. 
\end{proof}
\begin{lemma} \label{lem:linear-algebra-triples}
  Let 
  \begin{displaymath}
    0 \lto T' \lto T \lto T'' \lto 0 
  \end{displaymath}
  be an extension of $\alpha$-polystable triples, where $T'$ and $T''$
  are non-zero and $\mu_{\alpha}(T')=\mu_{\alpha}(T'')=\mu_{\alpha}(T)$.  
  Let $\alpha_m$ and $\alpha_M$ be the extreme $\alpha$ values for the 
  triple $T$, as defined in (\ref{alpha-bounds-m}) and 
  (\ref{alpha-bounds-M}), with the convention that $\alpha_M=\infty$ 
  if $n_1=n_2$ in $T$.
  Then the map
  \begin{displaymath}
    a_{1}\colon {E_{1}''}^{*} \otimes E_{1}' \oplus {E_{2}''}^{*}
    \otimes E_{2}' \to {E_{2}''}^{*} \otimes E_{1}'
  \end{displaymath}
  is not an isomorphism if $\alpha_m<\alpha<\alpha_M$.
  \end{lemma}
\begin{proof} If $a_1$ is an isomorphism then, applying 
 Lemma~\ref{lem:linear-algebra} fibrewise, it follows that one
  of the following alternatives must occur:
  \begin{itemize}
  \item[(a)] $E'_1 = E''_2 =0$ and $\phi'=\phi''=0$.
  \item[(b)] $E''_1=0$, $E'_1, E'_2, E''_2 \neq 0$ and $\phi'\colon E'_2
    \overset{\cong}{\longrightarrow} E'_1$.
  \item[(c)] $E'_2=0$, $E'_1, E''_1, E''_2 \neq 0$ and $\phi''\colon E''_2
    \overset{\cong}{\longrightarrow} E''_1$.
  \end{itemize}
  We shall consider each case in turn.
  \emph{Case (a).}  In this case $T'=(0,E_2,0), \ T''=(E_1,0,0),\ 
  T=(E_1,E_2,0)$, and the extension $0 \to T' \to T \to T'' \to 0$ is 
  trivial. It follows from $\mu_{\alpha}(T')=\mu_{\alpha}(T)$ that
  $\alpha=\mu(E_1)-\mu(E_2)=\alpha_m$.
  \emph{Case (b).}  In this case we have 
  $n_1 = n'_1$ and $n_2 = n'_2 + n''_2 = n'_1 + n''_2$.  Hence $n_2 > 
  n_1$. Furthermore, from $\mu_{\alpha}(T')=\mu_{\alpha}(T)$ we get
  $\mu(E_1)+\frac{\alpha}{2}=\mu(E_1\oplus E_2)+\frac{n_2}{n_1+n_2}$, 
  i.e.\ $\alpha=\frac{2n_2}{n_2-n_1}\alpha_m=\alpha_M$.
  \emph{Case (c).}  In this case we have 
 $n_2 = n''_2 $ and  $n_1 = n'_1+n_1'' = n'_1 + n''_2$.  Hence $n_1 > 
  n_2$. Furthermore, from $\mu_{\alpha}(T'')=\mu_{\alpha}(T)$ we get
$\alpha=\frac{2n_1}{n_1-n_2}\alpha_m=\alpha_M$.
If $n_1=n_2$ then case (a) is the only possibility, so  
$\alpha=\alpha_m$. If  $n_1\ne n_2$, then (a) or exactly one of (b) and (c) 
are the only possibilities, depending on whether $n_1<n_2$ or 
$n_1>n_2$. In both cases we see that $\alpha=\alpha_m$ or $\alpha=\alpha_M$. 
\end{proof}
%%%%%%%%%%%%%%%%%%%%%
%\subsection{Proof of Proposition \ref{prop:hyper-1-vanishing-criterion}}
\subsection{Estimate of $\chi(\mathbf{n}'',\mathbf{d}'',\mathbf{n}',\mathbf{d}')$ and
comparison of moduli spaces} 
\label{sec:chi-estimate}
%%%%%%%%%%%%%%%%%%%%%%%%%
Let $\alpha_c\in (\alpha_m,\alpha_M)$ be a critical value with 
$\alpha_c\ge 2g-2$, and let $\alpha_c^{\pm}$ be as in 
(\ref{eqtn:alpha-c-pm}). Fix $(\mathbf{n}',\mathbf{d}') = 
(n'_1,n'_2,d'_1,d'_2)$ and $(\mathbf{n}'',\mathbf{d}'') = 
(n''_1,n''_2,d''_1,d''_2)$ such that equations (\ref{eqn:C1}) and 
(\ref{eqn:C2}) are satisfied, and such that 
$\frac{n'_2}{n'_1+n'_2}\ne\frac{n''_2}{n''_1+n''_2}$.
Let 
$T'$ and $T''$ be triples of type 
$(\mathbf{n}',\mathbf{d}') = (n'_1,n'_2,d'_1,d'_2)$ and 
$(\mathbf{n}'',\mathbf{d}'') = (n''_1,n''_2,d''_1,d''_2)$ 
respectively. Suppose that $T'$ and $T''$ satisfy:
\begin{itemize}
\item both are $\alpha_c$-semistable with 
$\mu_{\alpha_c}(T')=\mu_{\alpha_c}(T'')$,
\item if $\frac{n'_2}{n'_1+n'_2}<\frac{n''_2}{n''_1+n''_2}$ then
$T'$ and $T''$ are both $\alpha^+_c$-stable,
\item if $\frac{n'_2}{n'_1+n'_2}>\frac{n''_2}{n''_1+n''_2}$ then
$T'$ and $T''$ are both $\alpha^-_c$-stable.
\end{itemize}
It follows by Proposition \ref{prop:H-alpha-stable} that  the chain 
$\widetilde{C^{\bullet}}(T'',T')$ is 
$(\alpha_c^+,2\alpha_c^+)$-polystable (if
$\frac{n'_2}{n'_1+n'_2}<\frac{n''_2}{n''_1+n''_2}$) or  it is
$(\alpha_c^-,2\alpha_c^-)$-polystable ( if  
$\frac{n'_2}{n'_1+n'_2}<\frac{n''_2}{n''_1+n''_2}$).
\begin{lemma}\label{lemma:chain-poly-2} The chain $\widetilde{C^{\bullet}}(T'',T')$ 
is $(\alpha_c,2\alpha_c)$ polystable. 
\end{lemma}
\begin{proof} The critical values for the chain form a discrete set of points in 
the $(\alpha_1,\alpha_2)$ plane. We can thus pick 
$\epsilon>0$ so that, with $\alpha_c^{\pm}=\alpha_c\pm \epsilon$, the point
$(\alpha_c^{\pm},2\alpha_c^{\pm})$ is not a critical point. We can in fact assume 
that there are no critical points in 
$B_{\epsilon}^{\circ}(\alpha_c,2\alpha_c)$, i.e.\ in the punctured ball of radius 
$\epsilon$ centered at $(\alpha_c,2\alpha_c)$. Thus $(\alpha_c^{\pm},2\alpha_c^{\pm})$ 
polystability is equivalent to $(\alpha_c,2\alpha_c)$ polystability. 
\end{proof}
\begin{proposition}\label{prop:hyper-1-vanishing-criterion}
  Let $T'=(E'_1,E'_2,\phi')$ and $T''=(E''_1,E''_2,\phi'')$ be triples
  as above.  In particular, suppose that both are $\alpha_c$-semstable and 
 that $\mu_{\alpha_c}(T')=\mu_{\alpha_c}(T'')$. Then
  \begin{displaymath}
    \chi(T'',T')\le 1-g < 0
  \end{displaymath}
\noindent if one of the following holds:
\begin{enumerate}
\item $\alpha_c > 2g-2$ and $T'$, $T''$ are both are 
$\alpha_c^-$-stable,
\item $\alpha_c \ge 2g-2$ and $T'$, $T''$ are both are 
$\alpha_c^+$-stable.
\end{enumerate}
\end{proposition}
\begin{proof}
{}From the long exact sequence \eqref{eq:long-exact-extension-complex}
  and the Rie\-mann-Roch formula we obtain
  \begin{equation}
    \label{eq:chi-T''-T'}
    \chi(T'',T')
    = (1-g)\bigl(\rk(C_{1}) - \rk(C_{0})\bigr) + \deg(C_{1}) -
    \deg(C_{0}). 
  \end{equation}
  
 \noindent where $C_1$ and $C_0$ are as in (\ref{eq:H-chain}).  
 {}From Proposition \ref{prop:H-alpha-stable}, we know that
  $\widetilde{C^{\bullet}}(T'',T')$ is
  $(\alpha_c^+,2\alpha_c^+)$-poly-stable or  
  $(\alpha_c^-,2\alpha_c^-)$-poly-stable.
 By the Lemma \ref{lemma:chain-poly-2} it 
  follows that $\widetilde{C^{\bullet}}(T'',T')$ is also
  $\balpha=(\alpha_c,2\alpha_c)$-poly-stable. Thus Lemma 
  \ref{lemma:stability-for-H_*}
  applies, and we can use the estimates (\ref{eq:ker-a_1}) and 
  (\ref{eq:im-a_1}). Together with
  \begin{align}
    \deg(C_{1}) &= \deg(\ker(a_{1})) + \deg(\im(a_{1}))\ , 
    \label{eq:deg(H_1)}\\
    \rk(C_{1})&=\rk(\ker(a_{1})) + \rk(\im(a_{1})) \ ,\label{eq:rk(H_1)}
  \end{align}
  
  \noindent these yield
    \begin{multline*}
      \deg(C_{1}) \leq (\mu_{\alpha_c}(T') - \mu_{\alpha_c}(T''))
    \bigl(\rk(C_{1}) - \rk(C_{0})\bigr) \\
    - \alpha_c\bigl(\rk(C_{0}) - \rk(\im(a_{1}))\bigr) + \deg(C_{0}).
    \end{multline*}
  
  \noindent Using that $\mu_{\alpha_c}(T') = \mu_{\alpha_c}(T'')$, we can then 
  deduce that
  \begin{displaymath}
      \deg(C_{1}) - \deg(C_{0})
      \leq  -\alpha_c\bigl(\rk(C_{0}) - \rk(\im(a_{1}))\bigr).
  \end{displaymath}
  Combining this with \eqref{eq:chi-T''-T'} we get
  \begin{equation}\label{eqn:chi-est}
    \chi(T'',T')
    \leq (1-g)\bigl(\rk(C_{1}) - \rk(C_{0})\bigr)
      - \alpha_c\bigl(\rk(C_{0}) - \rk(\im(a_{1}))\bigr)\ .
  \end{equation}
  If $\alpha_c\ge 2g-2$ then we get 
 $$
 \chi(T'',T')\le (1-g)\bigl(\rk(C_{0})+\rk(C_{1}) - 2\rk(\im(a_{1}))\bigr),
 $$
with equality if and only if $\alpha_c=2g-2$. Furthermore $\rk(a_{1}) 
\leq \rk(C_{0})$ and $\rk(a_{1}) \leq \rk(C_{1})$, with equality in both if and only if $a_1$ is an isomorphism. Thus in all cases we get 
$\chi(T'',T')\le 0$, with equality if and only if $\alpha_c=2g-2$ and
$a_1$ is an isomorphism. But by Lemma \ref{lem:linear-algebra-triples} if 
$\alpha_m<\alpha_c<\alpha_M$ then $a_1$ cannot be an isomorphism. 
Thus in all cases we get 
$\rk(C_{0})+\rk(C_{1})-2\rk(\im(a_{1})\ge 1)$ and hence $\chi(T'',T')\le 1-g$. 
 \end{proof}
 
 \noindent Combining Proposition \ref{prop:hyper-1-vanishing-criterion}
 with Corollary \ref{cor:codim-est}, we obtain
 
 \begin{theorem}\label{thm:codim}
Let $\alpha_c\in (\alpha_m,\alpha_M)$  be a critical value for 
triples of type $(\mathbf{n},\mathbf{d})$.  If $\alpha_c> 2g-2$ then 
the loci $\mathcal{S}_{\alpha_c^{\pm}} 
\subset\mathcal{N}^s_{\alpha_c^{\pm}}$ are contained in subvarieties of 
codimension at least $g-1$. In particular, they are contained in 
subvarieties of strictly positive codimension if $g\ge 2$. If 
$\alpha_c = 2g-2$ then the same is true for 
$\mathcal{S}_{\alpha_c^+}$.
 \end{theorem}
\begin{corollary}\label{cor:birationality}Let $\alpha_1$ and 
$\alpha_2$ be any two values in $(\alpha_m,\alpha_M)$ 
such that $\alpha_m <\alpha_1 <\alpha_2 <\alpha_M$ and $2g-2 
\le\alpha_1$. Then 
\begin{itemize}
\item The moduli spaces $\mathcal{N}^s_{\alpha_1}$ and 
$\mathcal{N}^s_{\alpha_2}$ have the same number of connected 
components, and 
\item The moduli space $\mathcal{N}^s_{\alpha_1}$ is irreducible if and only if
$\mathcal{N}^s_{\alpha_2}$ is.
\end{itemize} 
\end{corollary}

\begin{proof} This follows immediately from Theorem \ref{thm:codim} 
if $\alpha_1$ and $\alpha_2$ are non-critical, and from Theorem 
\ref{thm:codim} together with Lemma \ref{lemma:fliploci}
if either of them is critical.
\end{proof}

%%%%%%%%%%%%%%%%%%%%%%%%%%%%%%%%%%%%
\section{Special values of $\alpha$}\label{sect:special-alpha}
%%%%%%%%%%%%%%%%%%%%%%%%%%%%%%%%%%%%

Throughout this section we will assume that the triple 
$(E_1,E_2,\phi)$ has type $(n_1,n_2,d_1,d_2)$, with $n_1\geq n_2$. 
The case $n_1<n_2$ can be dealt with  via duality of  triples. In 
this section we identify some critical values in the range 
$(\alpha_m,\alpha_M)$ and describe their significance for the structure of 
$\alpha$-stable triples.

%%%%%%%%%%%%%%%%%%%%%%%%%%%%%%%%%%%%%%%%%%%%%%%%%%%%%%%%%%%%%%%
\subsection{The kernel of $\phi$ and the  parameter $\alpha$}
%%%%%%%%%%%%%%%%%%%%%%%%%%%%%%%%%%%%%%%%%%%%%%%%%%%%%%%%%%%%%%%

\begin{definition}
  For each integer $0\leq j < n_2$\ set
  \begin{equation}
    \label{eq:1.30}
    \alpha_j=\frac{2n_1n_2}{n_2(n_1-n_2)+
      (j+1)(n_1+n_2)}(\mu_1-\mu_2)\ .
  \end{equation}
\end{definition}
\begin{proposition}
  \label{prop:kernel-rank}
  Let $T=(E_1,E_2,\phi)$\ be a triple in which $n_1\geq n_2$. Let
  $N\subset E_2$\ be the kernel of $\phi: E_2\longrightarrow
  E_1$. Suppose that $T$\ is $\alpha$-semistable for some
  $\alpha>\alpha_j$.  Then $N$\ has rank at most $j$. In
  particular, if $T$\ is $\alpha$-semistable for some
  $\alpha>\alpha_0$\ then $N = 0$, i.e.\ $\phi$\ is
  injective.  
\end{proposition}
\begin{proof} Suppose that
\begin{equation}
\rk(N)=k>0\ .\label{eq:1.40}
\end{equation}
We consider the subtriples $T_N=(0,N,0)$\ and
$T_I=(I,E_2,\phi)$, where $I$\ denotes the
image sheaf $\im(\phi)$. If $N \neq 0$, then the triple
$T_N$\ is a proper subtriple, and so is $T_I$\ since $n_1\geq n_2$.
The $\alpha$-semistability condition applied to $T_N$\ yields
$$
\mu_N+\alpha\leq\mu+\frac{\alpha n_2}{n_1+n_2}\ ,
$$
where $\mu_N$\ denotes the slope of $N$\ and $\mu$\ is the slope of 
$E_1\oplus E_2$. Rearranging, we get
\begin{equation}
\alpha n_1\leq(n_1+n_2)(\mu -\mu_N)\ .\label{eq:1.41}
\end{equation}
The $\alpha$-semistability condition applied to $T_I$\ yields
\begin{equation}
\mu(E_2\oplus I)+\alpha\frac{n_2}{i+n_2}\leq
\mu+\alpha\frac{n_2}{n_1+n_2}\ ,\label{eq:1.42}
\end{equation}
where $i=\rk(I)$.  Furthermore, from the exact sequence
\begin{equation}
0 \longrightarrow N \longrightarrow E_2 \to I \longrightarrow 0\ ,
\label{eq:1.43}
\end{equation}
we get
\begin{align}
k+i&=n_2\ ,\label{eq:1.44} \\
k\mu_N+i\mu_I&=n_2\mu_2\ .\label{eq:1.45}
\end{align}
Using \eqref{eq:1.45} we can write
\begin{equation}
\mu(E_2\oplus I)=
\frac{2n_2\mu_2-k\mu_N}{2n_2-k} \label{eq:1.46}
\end{equation}
and hence \eqref{eq:1.42} yields
\begin{equation}
-k(n_1+n_2)\mu_N\leq(n_1+n_2)((2n_2-k)\mu-2n_2\mu_2)
+\alpha n_2(n_2-k-n_1)\ .\label{eq:1.47}
\end{equation}
Combining $k$ times \eqref{eq:1.41} and \eqref{eq:1.47} yields
\begin{align}
\alpha &\leq\frac{2n_2(n_1+n_2)}{n_2(n_1-n_2)+k(n_1+n_2)}(\mu-\mu_2)
\notag \\
&=\frac{2n_1n_2}{n_2(n_1-n_2)+k(n_1+n_2)}(\mu_1-\mu_2). \label{eq:1.48}
\end{align}
We have thus shown that if $\rk(N)=k$\ and the triple is
$\alpha$-stable, then
$$
\alpha\leq\alpha_{k-1}
$$
where $\alpha_{k-1}$\ is given by \eqref{eq:1.30} with $j=k-1$. 
Equivalently, if the triple is $\alpha$-semistable for some 
$\alpha>\alpha_{k-1}$, then 
$\rk(N)\neq k$. But
$$
\alpha_{k-1}>\alpha_{k}>\dots > \alpha_{n_2-1}\ .
$$
We can thus conclude that if the triple is $\alpha$-semistable with
$\alpha>\alpha_{k-1}$, then the rank of $N$\ is strictly less than $k$.
In particular, if $\alpha>\alpha_0$, where
\begin{equation}\label{alpha0}
\alpha_0 =\frac{2n_1n_2}{n_2(n_1-n_2)+
      (n_1+n_2)}(\mu_1-\mu_2), 
\end{equation}
then $T$ is injective.
\end{proof}
We thus have the following.
\begin{corollary}\label{injective-triples}
Let $\alpha>\alpha_0$, where $\alpha_0$ is given by (\ref{alpha0}).
\begin{itemize}
\item[(1)]
An $\alpha$-semistable  triple $(E_1,E_2,\phi)$  defines a sequence  of
the form
\begin{equation}
0 \lto E_2 \overset{\phi}{\lto} E_1 \lto F\oplus S \lto 0\ ,
\label{eq:3.20}
\end{equation}
where $F$\ is locally free and $S$\ is a torsion sheaf.
\item[(2)] If $n_1=n_2$ then
an $\alpha$-semistable  triple $(E_1,E_2,\phi)$  defines a sequence  of
the form
\begin{equation}
0 \lto E_2 \overset{\phi}{\lto} E_1 \lto S \lto 0,
\end{equation}
 where $S$ is a torsion sheaf of degree  $d_1 -d_2$.
\end{itemize}
\end{corollary}

\begin{lemma}\label{lemma:alpha-m-form} Let $\alpha_0$ be given by (\ref{alpha0}).
\begin{itemize}
\item[(1)] If $n_1>n_2$ then
\begin{equation}
\alpha_0=\frac{n_2(n_1-n_2)}{n_2(n_1-n_2)+n_1+n_2}\alpha_M
=\frac{2n_1n_2}{n_2(n_1-n_2)+n_1+n_2}\alpha_m ,
\end{equation}
where $\alpha_m$ and $\alpha_M$ are given by (\ref{alpha-bounds-m})
and (\ref{alpha-bounds-M}), respectively.
 
\item[(2)] 
If $n_1=n_2=n$ then 
\begin{equation}
\label{alpha0-equal-ranks}
\alpha_0=n\alpha_m=n(\mu_1-\mu_2)=d_1-d_2.
\end{equation}
\item[(3)] If $n_1\ge n_2$ then $\alpha_0\ge\alpha_m$,
with equality if and only if $\alpha_m=0$ or $n_2=1$. 
\end{itemize}
\end{lemma}
\begin{proof} Parts (1) and (2) are immediate. Using (1)
we compute 
$$\alpha_0-\alpha_m 
=\frac{n_1+n_2}{n_2(n_1-n_2)+n_1+n_2}(n_2-1)\alpha_m\ ,$$ from which 
(3) follows. 
\end{proof}

%%%%%%%%%%%%%%%%%%%%%%%%%%%%%%%%%%%%%%%%%%%%%%%%%%%%%%%%%%%%%%%
\subsection{The cokernel of $\phi$ and the  parameter $\alpha$}
%%%%%%%%%%%%%%%%%%%%%%%%%%%%%%%%%%%%%%%%%%%%%%%%%%%%%%%%%%%%%%%

In this section we will assume that $n_1>n_2$. The range for $\alpha$
is then $[\alpha_m,\alpha_M]$,  where $\alpha_m$ and $\alpha_M$ are given 
by (\ref{alpha-bounds-m})
and (\ref{alpha-bounds-m}). Let us define 
\begin{equation}
\alpha_t:=\alpha_M-\frac{n_1+n_2}{n_2(n_1-n_2)}.\label{alpha-torsion}
\end{equation}
\begin{proposition}
  \label{prop:torsion-degree-bound}
  Suppose that a triple $T=(E_1,E_2,\phi)$  of the 
form (\ref{eq:3.20}) with $n_1>n_2$ 
  is $\alpha$-semistable for some $\alpha>\alpha_m$.  Then
  $$
  s\leq\frac{n_2(n_1-n_2)}{(n_1+n_2)}(\alpha_M-\alpha)\ ,
  $$
  where $s$\ is the degree of $S$.  In particular, if
  $ \alpha >\alpha_t$, then $S=0$.
\end{proposition}
\begin{proof} If $T=(E_1,E_2,\phi)$\ is of the form
  in \eqref{eq:3.20}, with $S \neq 0$, then we can find a proper
  subtriple $T'=(E'_1,E_2,\phi)$\ of the form
\begin{equation}
0 \lto E_2 \overset{\phi}{\lto} E'_1 \lto S   \lto 0\ ,
\label{eq:3.30}
\end{equation}
Indeed, $E'_1$\ is the kernel of the sheaf map
$E_1\longrightarrow F\oplus S \longrightarrow S$. Notice that
$n'_1=n_2$\ and $d'_1=d_2+s$, where $n'_1\ ,\ d'_1$\ denote the rank
and degree of $E'_1$,  etc. We compute
\begin{equation}
\Delta_{\alpha}(T')=
\frac{n_1}{n_1+n_2}(\mu_2-\mu_1)+
\frac{\alpha}{2}\left(\frac{n_1-n_2}{n_1+n_2}\right)+
\frac{s}{2n_2}\ .\label{eq:3.31}
\end{equation}
But
$$
\frac{n_1}{n_1+n_2}(\mu_1-\mu_2)=
\frac{\alpha_M}{2}\left(\frac{n_1-n_2}{n_1+n_2}\right)
$$
and hence
\begin{equation}
\Delta_{\alpha}(T')=
\frac{n_1-n_2}{2(n_1+n_2)}\left( \alpha-\alpha_M +
\frac{n_1+n_2}{n_2(n_1-n_2)}s \right). \label{eq:3.32}
\end{equation}
If the triple is $\alpha$-semistable then
$\Delta_{\alpha}(T')\leq 0$\ and the
result follows.
\end{proof}

Let us define
 \begin{equation}\label{alpha_e}
  \alpha_e=
\max \{\alpha_m, \alpha_0, \alpha_t\}.
 \end{equation}
The following is then immediate.

\begin{proposition} \label{prop:alpha_e}
Let $\alpha>\alpha_e$. An $\alpha$-semistable triple $(E_1,E_2,\phi)$
defines an extension
\begin{equation}\label{free-extension}
0 \lto E_2 \overset{\phi}{\lto} E_1 \lto F \lto 0,
\end{equation}
with $F$ locally free.
\end{proposition}
It turns out that for  extension like (\ref{free-extension}),
arising from semistable triples the dimension of $H^1(E_2\otimes F^*)$
does not depend on the given triple. More precisely:
\begin{proposition} \label{vanishing-morphisms}
Let  $(E_1,E_2,\phi)$ be an $\alpha$-semistable triple defining an
extension like (\ref{free-extension}). Then 
$H^0(E_2\otimes F^*)=0$ and hence 
\begin{equation}\label{dimension-extensions}
\dim H^1(E_2\otimes F^*)=n_2d_1-n_1d_2+n_1(n_1-n_2)(g-1).
\end{equation}
\end{proposition}
\begin{proof}
{}From \cite[Lemma 4.5]{bradlow-garcia-prada:1996} we 
have that the 
$\alpha$-semistability of $(E_1,E_2,\phi)$  for arbitrary $\alpha$  
implies that  $H^0(E_1\otimes E_2^*)=0$.
{}From  (\ref{free-extension}), we have an injective homomorphism 
$F^*\to E_1^*$, which after  tensoring with $E_2$
gives that   $H^0(E_2\otimes F^*)$ injects in 
$H^0(E_1\otimes E_2^*)$, and hence the desired vanishing.
By Riemann--Roch we obtain (\ref{dimension-extensions}).
\end{proof}

 %%%%%%%%%%%%%%%%%%%%%%%%%%%%%%%%%%%%%%%%%%%%%%%%%%%%%%%
\section{Moduli space of triples  with  $n_1\neq n_2$}
\label{sec:n_1-neq-n_2}
%%%%%%%%%%%%%%%%%%%%%%%%%%%%%%%%%%%%%%%%%%%%%%%%%%%%%%

Throughout this section we assume that  $n_1>n_2$.
As usual, the case $n_1<n_2$ can be dealt with by triples duality. 
Recall that the allowed range for the stability parameter is
$\alpha_m\leq\alpha\leq\alpha_M$,
where $\alpha_m=\mu_1-\mu_2$ and 
$\alpha_M=\frac{2n_1}{n_1-n_2}\alpha_m$, and we assume  that 
$\mu_1-\mu_2>0$.  

We are now  ready to  jump all the way up to the large extreme value of the 
range, and 
describe the moduli space for $\alpha_M$.  

%%%%%%%%%%%%%%%%%%%%%%%%%%%%%%%%%%%%%%%%%%%%%%%%%%
\subsection{Moduli space  for   $\alpha=\alpha_M$}
%%%%%%%%%%%%%%%%%%%%%%%%%%%%%%%%%%%%%%%%%%%%%%%%%

\begin{proposition}\label{moduli-alpha-M}
Let $T=(E_1,E_2,\phi)$ be an $\alpha_M$-polystable triple then  
$E_1=\im\phi\oplus F$, and $T$ decomposes as the direct sum of two 
$\alpha_M$-polystable triples of the same $\alpha_M$-slope, $T'$ and 
$T''$, where $T'=(\im\phi,E_2,\phi)$, and  $T''=(F,0,0)$.
In particular, $T$ is never $\alpha_M$-stable. Moreover,
$E_2\cong\im \phi$ and $E_2$  and  $F$ are polystable.
\end{proposition}
\begin{proof}
 By Proposition, \ref{prop:alpha_e} $T$ defines an extension
\begin{equation}
0 \lto E_2 \overset{\phi}{\lto} E_1 \lto F \lto 0,
\end{equation}
with $F$ locally free. Let $T'=(\im\phi,E_2,\phi)$. Of course
$\phi:E_2\to\im\phi$ is an isomorphism, and
$$
\mu_{\alpha_M}(T')=\mu(E_2)+\frac{\alpha_M}{2},
$$
but this is equal to $\mu_{\alpha_M}(T)$ and hence
$T$ cannot be $\alpha_M$-stable and must decompose
as $T'\oplus T''$, where $T''=(F,0,0)$. It is clear from
the polystability of $T$  that
$T'$ and $T''$ are $\alpha_M$-polystable with the same 
$\alpha_M$-slope. Applying the $\alpha_M$-semistability 
condition to the
subtriples $(E_2',\phi(E_2'),\phi)\subset T'$ and 
$(F',0,0)\subset T''$,  we obtain that $\mu(E_2')\leq \mu(E_2)$
and $\mu(F')\leq \mu(F)$, and hence $E_2$ and $F$ are semistable.
In fact the polystability of $T'$ and $T''$ imply
the polystability of $E_2$ and $F$.
\end{proof}

Using Proposition 
\ref{same-rank-same-slope} and Corollary \ref{cor:largealphamoduli}
from Section \ref{subs:d_1=d_2}, we obtain the following corollary of 
Proposition \ref{moduli-alpha-M}:

\begin{corollary}
\label{cor:N-alpha-M}
Suppose that $n_1>n_2$ and $\mu_1-\mu_2>0$. Then 
\begin{equation}
  \label{eqn:alpha-Mmoduli}
  \begin{split}
\mathcal{N}_{\alpha_M}(n_1,n_2,d_1,d_2)&\cong
\mathcal{N}_{\alpha_M}(n_2,n_2,d_2,d_2)
\times M(n_1-n_2,d_1-d_2)\\
&\cong
M(n_2,d_2)
\times M(n_1-n_2,d_1-d_2)
\end{split}
\end{equation}
where $M(n,d)$ denotes the moduli space of polystable bundles of rank 
$n$ and degree $d$. In particular, 
$\mathcal{N}_{\alpha_M}(n_1,n_2,d_1,d_2)$ is irreducible. 
\end{corollary}

%%%%%%%%%%%%%%%%%%%%%%%%%%%%%%%%%%%%%%%%%%%%%%%%%%%%%%%
\subsection{Moduli space for large $\alpha$}
%%%%%%%%%%%%%%%%%%%%%%%%%%%%%%%%%%%%%%%%%%%%%%%%%%%%%%% 
%%%%%%%%%%%%%%%%%%%%%%%%%%%%%%%%%%%%%%%%%%%%%%%%%%%%%%
%\subsection{Large $\alpha$ and semistability of $E_2$}
%%%%%%%%%%%%%%%%%%%%%%%%%%%%%%%%%%%%%%%%%%%%%%%%%%%%%%

Let $\alpha_L$ be the largest critical value in 
 $(\alpha_m,\alpha_M)$, and let $\mathcal{N}_L$
(respectively $\mathcal{N}^s_L$) denote the moduli space of 
$\alpha$-polystable (respectively 
$\alpha$-stable) triples for
$\alpha_L<\alpha<\alpha_M$. We refer to $\mathcal{N}_L$ as the `large
$\alpha$' moduli space.

If a triple  $T=(E_1,E_2,\phi)$ defines an extension of the form
(\ref{free-extension}), then  $I=\im\phi$\ is a subbundle with 
torsion free 
  quotient in $E_1$, and  $\phi:E_2\lto I$ is an
  isomorphism. Thus we get a subtriple $T_I=(I,E_2,\phi)$ in
  which the bundles have the same rank and degree, and $\phi$ is an
  isomorphism.

\begin{proposition}  \label{prop:E2-stability}
Let $T=(E_1,E_2,\phi)$ represent a point in $\mathcal{N}_L$,
i.e.\ suppose that the triple is $\alpha$-semistable for some
$\alpha$\ in the range $\alpha_L<\alpha<\alpha_M$. Then
\begin{itemize}
\item[$(1)$] the triple $T_I=(I,E_2,\phi)$\ is $\alpha_M$-semistable
\item[$(2)$] the bundle $E_2$\ is semistable
\end{itemize}
\end{proposition}
\begin{proof}
(1). Let $T'=(E_1',E_2',\phi')$ be any subtriple of
$T_I$.  Since $T'$ is also a subtriple of $T$, we get
\begin{equation}
\mu_{\alpha}(T')\le
\mu_{\alpha}(T)\ .\label{eq:3.80}
\end{equation}

\noindent A direct computation shows that
\begin{align}
\mu_{\alpha}(T)&=
\mu_{\alpha}(T_I)+\frac{n_1-n_2}{n_1+n_2}(\mu(F)-\mu_2-\frac{\alpha}{2})
\notag \\
&=\mu_{\alpha}(T_I)+
\frac{n_1-n_2}{2(n_1+n_2)}(\alpha_M-\alpha), \label{eq:3.81}
\end{align}

\noindent where in the last line we have used the fact that $n_1>n_2$ in $T$
and hence
$\alpha_M=\frac{2n_1}{n_1-n_2}(\mu_1-\mu_2)=2(\mu(F)-\mu_2)$.
Thus for all $\alpha_L<\alpha<\alpha_M$ we have

$$\mu_{\alpha}(T')-\mu_{\alpha}(T_I)\le
\frac{n_1-n_2}{2(n_1+n_2)}(\alpha_M-\alpha)\ .$$

\noindent Taking the limit $\alpha\rightarrow\alpha_M$, we get

$$\mu_{\alpha_M}(T')-\mu_{\alpha_M}(T_I)\le 0\ ,$$

\noindent i.e.\ $T_I$ is $\alpha_M$-semistable.

(2). Let $E'_2\subset E_2$\ be any proper subsheaf.  Then
$(\phi(E'_2),E'_2,\phi)$\ is a subtriple of
$T_I$. Since $\phi:E_2\lto\phi(E_2)$ is
an isomorphism, this subtriple has $\mu(\phi(E'_2)=\mu(E'_2)$\ and
$n'_2=n'_1$. The $\alpha_M$-semistability condition of
$T_I$  thus gives
$$\mu(E'_2)+\frac{\alpha_M}{2}\leq\mu_2+\frac{\alpha_M}{2}\ ,$$ \noindent
(where we have made use of the fact that
$\mu(\phi(E_2)=\mu(E_2)=\mu_2$).  It follows from this that
$\mu(E'_2)\leq\mu_2$, i.e.\ that $E_2$\ is semistable.
\end{proof}

%%%%%%%%%%%%%%%%%%%%%%%%%%%%%%%%%%%%%%%%%%%%%%%%%%%%%%
%\subsection{Large $\alpha$ and semistability of $F$}
%%%%%%%%%%%%%%%%%%%%%%%%%%%%%%%%%%%%%%%%%%%%%%%%%%%%%%

\begin{proposition}
  \label{prop:triple-stable-implies-quotient-semistable}
  Suppose that the triple $T=(E_1,E_2,\phi)$\ is of the form in 
(\ref{free-extension}), i.e. 
  $$
   0 \lto E_2 \overset{\phi}{\lto} E_1 \lto F \lto 0,
  $$
with $F$\ locally free. Then there is an $\epsilon> 0$ such that $F$
is semistable if the triple is $\alpha$-semistable for any 
$\alpha >\alpha_M-\epsilon$. Indeed the
conclusion holds for any
\begin{equation}
0<\epsilon < \frac{2}{m(m-1)^2}\ ,\label{eq:3.60}
\end{equation}
where $m=n_1-n_2=\rk(F)$.
\end{proposition}
\begin{proof} Let $F'\subset F$\ be any proper
  subsheaf. 
Denote the rank and slope of $F$\ (resp.
  $F'$) by $m$\ and $\mu_F$\ (resp.\ $m'$\ and $\mu_{F'}$). We
  can always find $E'_1\subset E_1$\ such that
  $F'=E'_1/E_2$, i.e.\ such that we have
$$
0 \lto E_2 \overset{\phi}{\lto} E'_1 \lto F' \lto 0.
$$
Let $T'=(E_1',E_2,\phi)$. For convenience, define
\begin{equation}
\Delta_{\alpha}\equiv\Delta_{\alpha}(T')=
\mu_{\alpha}(T')-\mu_{\alpha}(T)\ .
\label{eq:3.61} 
\end{equation}
Using
\begin{align}
n_1&=n_2+m\ , \notag \\
n'_1&=n_2+m'\ , \notag \\
n_1\mu_1 &= n_2\mu_2+m\mu_F\ ,\label{eq:3.62}\\
n'_1\mu'_1 &= n_2\mu_2+m'\mu_{F'}\ , \notag
\end{align}
we get
\begin{equation}
\mu_{F'}-\mu_F=\frac{(2n_2+m)(2n_2+m')}{2n_2m'}\Delta_{\alpha}-
\left(\frac{m-m'}{2m'}\right)(\alpha -2(\mu_F-\mu_2))\ .\label{eq:3.63}
\end{equation}
But $2(\mu_F-\mu_2)=\alpha_M$. Thus, setting
\begin{equation}
\alpha=\alpha_M-\epsilon\ ,\label{eq:3.64}
\end{equation}
we get
\begin{equation}
\mu_{F'}-\mu_F=\frac{(2n_2+m)(2n_2+m')}{2n_2m'}\Delta_{\alpha}+
\left(\frac{m-m'}{m'}\right)\frac{\epsilon}{2}\ .\label{eq:3.65}
\end{equation}
If now we take
$$
\frac{\epsilon}{2}<\frac{1}{m(m-1)^2}\ ,
$$
then for all $0<m'<m$ we get
\begin{equation}
\left(\frac{m-m'}{m'}\right)\frac{\epsilon}{2}<
\frac{1}{m(m-1)}\ .\label{eq:3.66}
\end{equation}
Hence, if the triple is $\alpha$-semistable, so that 
$\Delta_{\alpha}\le 0$, then we get 
\begin{equation}
\mu_{F'}-\mu_F<\frac{1}{m(m-1)}\ .\label{eq:3.67}
\end{equation}
Since $\mu_F$ and $\mu_{F'}$ are rational numbers, the first with 
denominator $m$, and the second with denominator $m'\le (m-1)$, 
equation (\ref{eq:3.67}) equivalent to the condition 
$\mu_{F'} - \mu_F\leq 0$. 
\end{proof}

We can combine   
Propositions~\ref{prop:alpha_e},
\ref{prop:triple-stable-implies-quotient-semistable} and 
\ref{prop:E2-stability} to obtain  the following.
\begin{proposition}
\label{prop:triple-stable-implies-bundles-semistable} 
Let $T=(E_1,E_2,\phi)$ be an  $\alpha$-semistable triple for
some $\alpha$\ in the range $\alpha_L<\alpha<\alpha_M$. Then
$T$ is of the form
$$
0 \lto E_2 \overset{\phi}{\lto} E_1 \lto F \lto 0, 
$$
with  $F$ locally free, and $E_2$ and $F$ are semistable.
\end{proposition}
In the converse direction we have:
\begin{proposition}
Let $T=(E_1,E_2,\phi)$ be a triple 
of  the form
$$
0 \lto E_2 \overset{\phi}{\lto} E_1 \lto F \lto 0, 
$$
with   $F$  locally free.
If  $E_2$ is semistable  and $F$ is stable then $T$ 
 is $\alpha$-stable for  $\alpha=\alpha_M-\epsilon$ in the range 
$\alpha_L<\alpha<\alpha_M$. 
\end{proposition}
\begin{proof}
Any subtriple $T'=(E_1',E_2',\phi')$ defines a commutative diagram
  $$
  \begin{CD}
  0@>>>E_2@>\phi>>E_1@>>>F@>>>0\\
  @.@AAA@AAA@AAA\\
  0@>>>E'_2@>\phi'>>E'_1@>>>F'@>>>0,
  \end{CD}
$$
where $F'\subset F$. Then
\begin{align}
\Delta_{\alpha}\equiv\Delta_{\alpha}(T') & = 
\mu_{\alpha}(T')-\mu_{\alpha}(T) \notag\\
&=\mu(E_1'\oplus E_2')-\mu(E_1\oplus E_2)+
\alpha(\frac{n_2'}{n_1'+n_2'}- \frac{n_2}{n_1+n_2}).
\end{align}
 Denote the rank and slope of $F$\ (resp.
  $F'$) by $m$\ and $\mu_F$\ (resp.\ $m'$\ and $\mu_{F'}$).
Using
\begin{align}
n_1&=n_2+m\ , \notag \\
n'_1&=n'_2+m'\ , \notag \\
n_1\mu_1 &= n_2\mu_2+m\mu_F\ ,\notag \\
n'_1\mu'_1 &= n'_2\mu_2+m'\mu_{F'}\ , \notag
\end{align}
and the fact that $\alpha_M=2(\mu_F-\mu_2)$, and setting
 $\alpha=\alpha_M-\epsilon$,
we obtain 
\begin{align}
\Delta_{\alpha} & = 
\frac{2n_2'+\mu_2'+m\mu_{F'}}{2n_2'+m'}-\frac{2n_2+\mu_2+m\mu_F}{2n_2+m}
+ 2(\mu_F-\mu_2)(\frac{n_2'}{2n_2'+m'}- \frac{n_2}{2n_2+m}) \notag \\
& -\epsilon(\frac{n_2'}{2n_2'+m'}- \frac{n_2}{2n_2+m})\notag \\
&= \frac{2n_2'}{2n_2'+m'}(\mu_2'-\mu_2) + \frac{m'}{2n_2'+m'}(\mu_{F'}-\mu_F)
-\epsilon(\frac{n_2'm-n_2m'}{(2n_2'+m')(2n_2+m)}).  \label{alpha-epsilon}
\end{align}
Clearing denominators in (\ref{alpha-epsilon}) we obtain 
\begin{align}
\hat\Delta_{\alpha}=n_2'(2\Delta_2-\frac{m}{2n_2+m}\epsilon) +
                    m'(\Delta_F+\frac{n_2}{2n_2+m}\epsilon),
\end{align}
where 
\begin{equation}
\hat\Delta_{\alpha}= (2n_2'+m')\Delta_{\alpha}, \quad\
\Delta_2 =\mu_2'-\mu_2, \quad\ \mathrm{and}\ \quad\ 
\Delta_F=\mu_{F'}-\mu_F. 
\end{equation}
Now suppose that  $E_2$ is semistable and $F$ is stable. The 
semistability of $E_2$ implies that 
$$
2\Delta_2-\frac{m}{2n_2+m}\epsilon<0,
$$
while the  stability of $F$ implies there exists $\delta>0$ such that
$$\Delta_F\leq -\delta <0.$$
Thus by taking 
$$
\epsilon <\frac{2n_2 +m}{n_2}\delta,
$$
we have 
$$
\Delta_F+\frac{n_2}{2n_2+m}\epsilon<0,
$$
and hence $\Delta_\alpha <0$, completing the proof. 
\end{proof}
\begin{theorem}\label{thm:largealpha}
Let $n_1>n_2$ and $d_1/n_1> d_2/n_2$. The moduli space 
$\mathcal{N}^s_L(n_1,n_2,d_1,d_2)$ is  smooth 
of dimension 
$$
(g-1)(n_1^2 + n_2^2 - n_1 n_2) - n_1 d_2 + n_2 d_1 + 1,
$$
and is birationally equivalent to a $\mathbb{P}^N$-fibration over
$M^s(n_1-n_2,d_1-d_2) \times M^s(n_2,d_2)$, where 
$$N=n_2d_1-n_1d_2+n_1(n_1-n_2)(g-1)-1.$$
In particular, $\mathcal{N}_L^s(n_1,n_2,d_1,d_2)$ is non-empty and
irreducible. 

If $\GCD(n_1-n_2,d_1-d_2)=1$ and $\GCD(n_2,d_2)=1$, the birational 
equivalence is an isomorphism. 

Moreover, $\mathcal{N}_L(n_1,n_2,d_1,d_2)$ is irreducible and hence 
birationally equivalent to $\mathcal{N}_L^s(n_1,n_2,d_1,d_2)$. 
\end{theorem}
\begin{proof}
For every triple $T=(E_1,E_2,\phi)$  in $\mathcal{N}^s_L(n_1,n_2,d_1,d_2)$,
the homomorphism  $\phi$ is injective and hence,  by (5) in Proposition 
\ref{thm:smoothdim}, $T$ defines a smooth point in the moduli space, whose dimension is 
then given by (4) in Proposition \ref{triples-critical-range}. 

Given 
$F\in M^s(n_1-n_2,d_1-d_2)$ and $E_2\in M^s(n_2,d_2)$, we know
from Proposition \ref{triples-critical-range} that
every extension 
$$
0 \lto E_2 \overset{\phi}{\lto} E_1 \lto F \lto 0, 
$$
determines  a triple $T=(E_1,E_2,\phi)$ in  $\mathcal{N}^s_L(n_1,n_2,d_1,d_2)$.
These extensions are  classified by
  $H^1(E_2\otimes F^*)$. In fact two classes defining the same element
in the projectivatization $\mathbb{P}H^1(E_2\otimes F^*)$ define equivalent
extensions and therefore equivalent triples.
Now, 
$$
\deg(E_2\otimes F^*)=(n_1-n_2)d_2-n_2(d_1-d_2)= n_1n_2(\mu_2-\mu_1)<0
$$
and, since $E_2\otimes F^*$\ is semistable, then
$H^0(E_2\otimes F^*)=0$.  Hence, by the Riemann--Roch theorem
$$h^1(E_2\otimes F^*)=n_2d_1-n_1d_2+n_1(n_1-n_2)(g-1).$$
In particular this dimension is constant as $F$ and $E_2$ vary in the
corresponding moduli spaces.

We  can describe this globally in terms of Picard sheaves. To do that  
we consider first the case in which $\GCD(n_1-n_2,d_1-d_2)=1$ and  
$\GCD(n_2,d_2)=1$. In this situation there exist universal 
 bundles $\mathbb{F}$  and  $\mathbb{E}_2$ over $X \times M(n_1-n_2,d_1-d_2)$
and $X\times M(n_2,d_2)$, respectively.
Consider the canonical projections 
$\pi:X M(n_1-n_2,d_1-d_2)\times M(n_2,d_2)\times \to
M(n_1-n_2,d_1-d_2) \times M(n_2,d_2)$
$\nu:X \times M(n_1-n_2,d_1-d_2)\times M(n_2,d_2)\to 
X \times  M(n_1-n_2,d_1-d_2)$, and 
 $\pi_2:X \times M(n_1-n_2,d_1-d_2)\times M(n_2,d_2)\to 
X\times M(n_2,d_2)$. 
The Picard sheaf
$$
\mathcal{S}:=R^1\pi_*(\pi_2^*\mathbb{E}_2\otimes \nu^* \mathbb{F^*}),
$$
is then locally free (a Picard bundle) and we can identify 
$\mathcal{N}_L(n_1,n_2,d_1,d_2)=\mathcal{N}_L^s(n_1,n_2,d_1,d_2)$
with  $\mathcal{P}=\mathbb{P}(\mathcal{S})$. This is indeed a 
$\mathbb{P}^N$ fibration
with $N=n_2d_1-n_1d_2+n_1(n_1-n_2)(g-1)-1$, which in particular
is non-empty since $M(n_1-n_2,d_1-d_2)$ and  $M(n_2,d_2)$ are 
non-empty and $N>0$.

If $\GCD(n_1-n_2,d_1-d_2)\neq 1$ and  $\GCD(n_2,d_2)\neq 1$, the 
universal bundles and hence the Picard bundle do not exist but the 
projectivization over $M^s(n_1-n_2,d_1-d_2) \times M^s(n_2,d_2)$ 
does. One way to show this is to work in the open set $R$ of the Quot 
scheme corresponding to stable bundles. The point is that an 
appropriate linear group 
$\GL$ acts on $R$, with the centre acting trivially and such that 
$\PGL$ acts freely with the quotient being $M^s(n_1-n_2,d_1-d_2) 
\times M^s(n_2,d_2)$. For the action on the projective bundle 
associated to the universal bundle over $R$, the centre of $\GL$ 
still acts trivially, and one can use standard descent arguments to 
obtain a $\mathbb{P}^N$ fibration 
$\mathcal{P}$ over $M^s(n_1-n_2,d_1-d_2) \times M^s(n_2,d_2)$. 

W now show that $\mathcal{N}_L^s(n_1,n_2,d_1,d_2) - \mathcal{P}$  has 
strict positive codimension in  $\mathcal{N}_L^s(n_1,n_2,d_1,d_2)$. 
This follows from two facts. The first one is that  any family of 
semistable bundles of rank $n_1-n_2$ and degree $d_1-d_2$ depends on 
a number of parameters strictly less than the dimension of 
$M^s(n_1-n_2,d_1-d_2)$ (similarly for any family of semistable 
bundles of rank $n_2$ and degree $d_2$). The second fact is that the 
dimension of $H^1(E_2\otimes F^*)$ is fixed by the Riemann--Roch 
theorem (we use here that $E_2$ and $F$ are semistable).

To prove the last statement, i.e.\ to extend the results to 
$\mathcal{N}_L(n_1,n_2,d_1,d_2)$, we consider the family 
$\tilde{\mathcal{P}}$ of equivalence classes of extensions 
$$
0 \lto E_2 \overset{\phi}{\lto} E_1 \lto F \lto 0, 
$$
where $F$ and $E_2$ are semistable. Clearly, $\tilde{\mathcal{P}}$ 
contains the family $\mathcal{P}$. The family $\tilde{\mathcal{P}}$ 
is irreducible. This is because  since $F$ and $E_2$ are semistable  
they vary (for fixed ranks and degrees) in irreducible families 
$\mathcal{F}$ and $\mathcal{E}_2$, respectively, and  as shown above
$H^0(E_2\otimes F^*)=0$. Hence  
$\tilde{\mathcal{P}}$ is a projective bundle over $\mathcal{F}\times 
\mathcal{E}_2$.  {}From Proposition 
\ref{prop:triple-stable-implies-bundles-semistable}, 
we know that $\mathcal{N}_L(n_1,n_2,d_1,d_2)\subset 
\tilde{\mathcal{P}}$, and since 
$\alpha$-semistability is an open condition  we have that
$\mathcal{N}_L(n_1,n_2,d_1,d_2)$ is irreducible.
\end{proof}

\begin{remark}
If $n_1<n_2$, we have an analogous theorem for 
$\mathcal{N}^s_L(n_1,n_2,d_1,d_2)$ via the isomorphism
$$
\mathcal{N}^s_\alpha(n_1,n_2,d_1,d_2) = \mathcal{N}^s_\alpha(n_2,n_1,-d_2,-d_1)
$$
given by duality.
\end{remark}

%%%%%%%%%%%%%%%%%%%%%%%%%%%%%%%%%%%%%%%%%%%%%%%%%%%%%%%%
\subsection{Moduli space for $2g-2\leq \alpha<\alpha_M$}
%%%%%%%%%%%%%%%%%%%%%%%%%%%%%%%%%%%%%%%%%%%%%%%%%%%%%%%%

\begin{theorem}\label{thm:irreducibility-moduli-stable-triples}
Let $\alpha$ be any value in the range $2g-2\leq\alpha< 
\alpha_M$. Then 
$\mathcal{N}^s_\alpha$ is birationally equivalent to 
$\mathcal{N}^s_L$.  In particular it is non-empty and irreducible. 
\end{theorem}
\begin{proof}
This follows from Corollary \ref{cor:birationality} and Theorem 
\ref{thm:largealpha}. 
\end{proof}

\begin{corollary}\label{cor:gcd=1}
Let $(n_1,n_2,d_1,d_2)$ be such that $\GCD(n_2,n_1+n_2,d_1+d_2)=1$.
If  $\alpha$ is a  generic value satisfying  
$2g-2\leq\alpha<\alpha_M$, then  $\mathcal{N}_\alpha$ is birationally 
equivalent to $\mathcal{N}_L$, and in particular it is irreducible. 
\end{corollary}
\begin{proof}
{}From (4) in Proposition \ref{triples-critical-range} one has that 
$\mathcal{N}_\alpha=\mathcal{N}_\alpha^s$ if
$\GCD(n_2,n_1+n_2,d_1+d_2)=1$ and $\alpha$ is generic.
In particular, $\mathcal{N}_L=\mathcal{N}_L^s$, and hence the result
follows from Theorem \ref{thm:irreducibility-moduli-stable-triples}.
\end{proof}

%%%%%%%%%%%%%%%%%%%%%%%%%%%%%%%%%%%%%%%%%%%%%%%%%%%%%%%%%
\section {Moduli space of triples with $n_1 = n_2$}
\label{sec:n_1=n_2}
%%%%%%%%%%%%%%%%%%%%%%%%%%%%%%%%%%%%%%%%%%%%%%%%%%%%%%%%%
Throughout this section we will assume that $n_1=n_2=n$
and $d_1\geq d_2$.

%%%%%%%%%%%%%%%%%%%%%%%%%%%%%%%%%%%%%%%%%%%%%%%%%%%%%%%%%%%%
\subsection{Moduli space for $d_1=d_2$}\label{subs:d_1=d_2}
%%%%%%%%%%%%%%%%%%%%%%%%%%%%%%%%%%%%%%%%%%%%%%%%%%%%%%%%%%%%

\begin{proposition}\label{same-rank-same-slope}
Suppose  that $n_1=n_2=n$ and $d_1=d_2=d$. Let $T=(E_1,E_2,\phi)$ be 
a triple of type $(n_1,n_2,d_1,d_2)$, and let $\alpha>0$. Then $T$ is 
$\alpha$-(poly)stable if and only if $E_1$ and $E_2$ are (poly)stable and 
$\phi$ is an isomorphism. 
\end{proposition}
\begin{proof}
In this case $\alpha_0=\alpha_m=0$  and hence for  every 
$\alpha$-semistable triple $T=(E_1,E_2,\phi)$ with $\alpha>0$, $\phi$ 
must be injective and therefore an isomorphism. The polystability of 
$E_1$ and $E_2$ is now  straightforward to see. To show the converse, 
suppose that $E_1$ and $E_2$ are both polystable and let 
$T'=(E_1',E_2',\phi')$ be any subtriple of $T$. 

\begin{align*}
\mu_\alpha(T')&= \mu(E_1'\oplus E_2')+\alpha\frac{n_2'}{n_1'+n_2'}\\
              &\leq  \mu(E_1\oplus E_2)+\alpha\frac{n_2'}{n_1'+n_2'}\\
              &\leq  \mu_\alpha(T)+
                \alpha(\frac{n_2'}{n_1'+n_2'}-\frac{1}{2})\\
                            &\leq  \mu_\alpha(T),
\end{align*}
since $n_1'\geq n_2'$ for $\phi$ is injective. 

\end{proof}

\begin{corollary}\label{cor:largealphamoduli}
There is a surjective map from $\mathcal{N}_\alpha(n,n,d,d)$ to 
$M(n,d)$ which defines an isomorphism between the two moduli spaces. 
In particular $\mathcal{N}_\alpha(n,n,d,d)$ is non-empty and 
irreducible. 

\end{corollary}
\begin{proof}
{}From Proposition \ref{same-rank-same-slope} it is clear that 
we have a surjective map, say 
$$\pi:\mathcal{N}_\alpha(n,n,d,d)\to M(n,d)\ .$$
Suppose that 
$\pi([T])=\pi([T'])$, where $[T]$ and $[T']$ are points 
in $\mathcal{N}_\alpha(n,n,d,d)$ represented by triples 
$T=(E,E,\phi)$ and $T'=(E',E',\phi')$ respectively. We may assume 
that $T$ and $T'$ are polystable triples, and hence that $E$ and 
$E'$ are polystable bundles.  Thus, since 
$\pi([T])=\pi([T'])$, we can find an isomorphism 
$h_1:E\mapsto E'$. Set $h_2=\phi'\circ h_1\circ\phi^{-1}$  
(remember that $\phi$ and $\phi'$ are bundle isomorphisms!). 
Then $(h_1,h_2)$ defines an isomorphism form $T$ to $T'$. 
Thus $\pi$ is injective.
\end{proof}

%%%%%%%%%%%%%%%%%%%%%%%%%%%%%%%%%%%%%%%%%%%%%%%%%%%%%%%%%%%%%%%
\subsection{Bounds on $E_1$ and $E_2$ for  $\alpha>\alpha_0$}
%%%%%%%%%%%%%%%%%%%%%%%%%%%%%%%%%%%%%%%%%%%%%%%%%%%%%%%%%%%%%%%

\begin{lemma}  \label{lemma:technical}
  Let $(E_1,E_2,\phi)$ be a triple with
  $\ker\phi=0$. Let $(E'_1,E'_2,\phi')$ be a subtriple with
  $n'_1=n'_2=n'$. Thus we get the following diagram, in which $S$ and
  $S'$ are torsion sheaves:
  $$
  \begin{CD}
  0@>>>E_2@>\phi>>E_1@>>>S@>>>0\\
  @.@AAA@AAA@AAA\\
  0@>>>E'_2@>\phi'>>E'_1@>>>S'@>>>0.
  \end{CD}
  $$
  Then
  \begin{align*}
  \Delta_{\alpha}(T')\equiv
  \mu_{\alpha}(T')-\mu_{\alpha}(T)
  &=(\mu(E'_2)-\mu_2)+\frac{1}{2}\left(\frac{s'}{n'}-\frac{s}{n}\right)\\  &=(\mu(E'_1)-\mu_1)-\frac{1}{2}\left(\frac{s'}{n'}-\frac{s}{n}\right).
  \end{align*}
  Here $s$\ and $s'$\ are the degrees of $S$\ and $S'$\ respectively.
\end{lemma}
\begin{proof} {}From the above diagram we get
\begin{align*}
 n\mu_2+s&=n\mu_1\ ,\\
 n\mu(E'_2)+s'&=n\mu(E'_1)\ .
\end{align*}
Thus
\begin{align*}
\mu_{\alpha}(T')=
\frac{1}{2}(\mu(E'_1)+\mu(E'_2))+\frac{\alpha}{2}
&=\frac{1}{2}\left(2\mu(E'_2) + \frac{s'}{n'}\right)+\frac{\alpha}{2}\\
&=\frac{1}{2}\left(2\mu(E'_1)- \frac{s'}{n'}\right)+\frac{\alpha}{2},
\end{align*}
and similarly for $\mu_{\alpha}(T)$.
\end{proof}
\begin{proposition}
  \label{prop:slope-bound-on-subbdls}
  Let $(E_1,E_2,\phi)$\ be an $\alpha$-semistable triple with 
$\ker\phi=0$. Then
\begin{itemize}
\item[$(1)$] For any subsheaf $E'_1\subset E_1$
  $$
  \mu(E'_1)\leq\mu_1+\frac{1}{2}(n-1)(\mu_1-\mu_2)\ .
  $$
\item[$(2)$] For any subsheaf $E'_2\subset E_2$
  $$
  \mu(E'_2)\leq\mu_2+\frac{1}{2}(\mu_1-\mu_2)\ .
  $$
\end{itemize}
\end{proposition}
\begin{proof} Since $\ker\phi=0$\ the results of Lemma
  \ref{lemma:technical} apply. Furthermore, any subsheaf $E'_1\subset
  E_1$\ is part of a subtriple $(E'_1, E'_2,\phi')$\ with
  $n'_1=n'_2=n'$. Likewise, given any subsheaf $E'_2\subset E_2$, we
  can take $E'_1=\phi(E'_2)$. Thus we can use the results of Lemma
  \ref{lemma:technical}, plus the fact that $\alpha$-stability implies
  $\Delta_{\alpha}(T')<0$\ for all subtriples, to conclude
  $$
  \mu(E'_1)-\mu_1-\frac{1}{2}\left(\frac{s'}{n'}-
  \frac{s}{n}\right)<0\
  $$
  for all $E'_1\subset E_1$. Similarly
  $$
  \mu(E'_2)-\mu_2+\frac{1}{2}\left(\frac{s'}{n'}-
    \frac{s}{n}\right) < 0
  $$
  for all $E'_2\subset E_2$.  The results now follow using the fact
  that $0\leq s'\leq s$\ and $1\leq n'< n$.
\end{proof}

%%%%%%%%%%%%%%%%%%%%%%%%%%%%%%%%%%%%
\subsection{Stabilization of moduli}
\label{sec:stabilize}
%%%%%%%%%%%%%%%%%%%%%%%%%%%%%%%%%%%%

\begin{theorem}[Stabilization Theorem]
  \label{thm:stabilization}
Let $\alpha_0$ be as in  (\ref{alpha0-equal-ranks}).  
\begin{itemize}
\item[$(1)$] If  $\alpha_1,\ \alpha_2$\ be any real numbers
  such that $\alpha_0<\alpha_1\leq\alpha_2$,
$$
\mathcal{N}_{\alpha_1}(n,n,d_1,d_2)\subseteq\mathcal{N}_{\alpha_2}(n,n,d_1,d_2). 
$$
\item[$(2)$] There is a real number
  $\alpha_{L}\geq\alpha_0$\ such that
$$
\mathcal{N}_{\alpha_1}(n,n,d_1,d_2)
=\mathcal{N}_{\alpha_2}(n,n,d_1,d_2)
$$
for all $\alpha_1\geq \alpha_2>\alpha_L$.
\end{itemize}
\end{theorem}
\begin{proof}
  (1). Recall from \ref{injective-triples} that if $\alpha>\alpha_0$\ then any triple, 
  say $T=(E_1,E_2,\phi)$,
  in $\mathcal{N}_\alpha(n,n,d_1,d_2)$\ has $\rk(\phi)=n$. It follows
  that in any subtriple, say $T'=(E'_1,E'_2,\phi')$, the rank of $E'_1$\ 
  is at least as big as the rank of $E'_2$, i.e.\ $n'_1\geq n'_2$. We
  treat the cases $n'_1> n'_2$\ and $n'_1= n'_2$\ separately. In both
  cases we must show that
$$
\Delta_{\alpha_1}(T')\leq0\ \Rightarrow\
\Delta_{\alpha_2}(T')\leq0\ 
$$
if $\alpha_1\leq\alpha_2$.
If $n'_1=n'_2$\ then for any $\alpha$
\begin{equation}
\Delta_{\alpha}(T')=
\mu(E'_1\oplus E'_2)-\mu(E_1\oplus E_2)\ .\label{eq:2.60}
\end{equation}
In particular,  $\Delta_{\alpha}(T')$\ is independent of
$\alpha$\ and hence
$\Delta_{\alpha_1}(T')=
\Delta_{\alpha_2}(T')$.
If $n'_1>n'_2$, then for any $\alpha$
\begin{align}
\Delta_{\alpha}(T')&=
\mu(E'_1\oplus E'_2)-\mu(E_1\oplus E_2)+
\left(\frac{n'_2}{n'_1+n'_2}-\frac{1}{2}\right)\alpha. \label{eq:2.61}
\end{align}
For each subtriple, $\Delta_{\alpha}(T')$\ is thus a linear
function of $\alpha$, with slope
\begin{equation}
\lambda(T')\ =\ \left(\frac{n'_2}{n'_1+n'_2}-\frac{1}{2}\right)\ =\
\frac{n'_2-n'_1}{2(n'_1+n'_2)}\ \label{eq:2.62}
\end{equation}
and constant term
\begin{equation}
M(T')= \mu(E'_1\oplus E'_2)-\mu(E_1\oplus E_2)\ .\label{eq:2.63}
\end{equation}
We see that if $n'_1>n'_2$\ then $\lambda(T')<0$. It follows
from this that
$$
\Delta_{\alpha_1}(T')\leq 0\ \Longrightarrow\
\Delta_{\alpha_2}(T')\leq 0\ 
$$
if $\alpha_1\leq\alpha_2$.
(2). Consider any $\alpha_1,\alpha_2$\ such that
$\alpha_0<\alpha_1\leq\alpha_2$. By Part (1), the difference (if any)
between $\mathcal{N}_{\alpha_1}$\ and $\mathcal{N}_{\alpha_2}$\ is due
entirely to triples which are $\alpha_2$-stable but not
$\alpha_1$-stable. Any such triple must have a subobject, say
$T'=(E'_1,E'_2,\phi')$, for which
\begin{equation}
\Delta_{\alpha_2}(T')\leq 0 <
\Delta_{\alpha_1}(T')\ .\label{eq:2.64}
\end{equation}
As in (1), we need only consider subobjects for which the rank of
$E'_1$\ is at least as big as the rank of $E'_2$, i.e.\ $n'_1\geq
n'_2$. Clearly \eqref{eq:2.64} is not possible for a subobject with
$n'_1=n'_2$\ (since in that case
$\Delta_{\alpha_1}(T')=\Delta_{\alpha_2}(T')$). Suppose
then that $n'_1>n'_2$. By \eqref{eq:2.61} and the fact that for such a
subobject $\lambda(T')<0$, we get that
\begin{equation}
\Delta_{\alpha}(T') \geq0\  \iff
\alpha \leq \frac{M(T')}{-\lambda(T')}\ .\label{eq:2.64a}
\end{equation}
We claim that there is a bound, $\alpha_L$, depending
  only on the degrees and ranks of $E_1$\ and $E_2$, such that
  \begin{equation}
    \frac{M(T')}{-\lambda(T')}\leq \alpha_L
    \label{eq:2.65}
  \end{equation}
  for all possible subtriples with $n'_1>n'_2$.
For a triple  $T=(E_1,E_2,\phi)$  in
$\mathcal{N}_{\alpha_2}$ 
Proposition~\ref{prop:slope-bound-on-subbdls} applies, giving
upper bounds on slopes of subsheaves of both $E_1$\ and $E_2$. Using
these bounds we compute
\begin{equation}
M(T')\leq
\frac{nn'_1}{2(n'_1+n'_2)}(\mu_1-\mu_2)\ .\label{eq:2.66}
\end{equation}
Combined with \eqref{eq:2.62}, this gives
\begin{align*}
\frac{M(T')}{-\lambda(T')}&\leq
\frac{nn'_1}{(n'_1-n'_2)}(\mu_1-\mu_2)\\
&\leq n(n-1)(\mu_1-\mu_2).
\end{align*}
We can thus take
\begin{equation}
\alpha_L=n(n-1)(\mu_1-\mu_2)\ .\label{eq:2.67}
\end{equation}
We can now complete the proof of Part (2): if
$\alpha_1>\alpha_L$\ then no triple in
$\mathcal{N}_{\alpha_2}$\ can have a subobject satisfying
\eqref{eq:2.64a}.  Hence
$\mathcal{N}_{\alpha_2}=\mathcal{N}_{\alpha_1}$.
\end{proof}
\begin{remark}
If $n=2$ then  $\alpha_L=\alpha_0=d_1-d_2$, i.e.\ the stabilization
parameter coincides with the injectivity  parameter.
\end{remark}
It is clear from (\ref{eq:2.67}) that $\alpha_L=0$ correspond to the 
following especial cases. 
\begin{proposition} The condition $\alpha_L=0$ holds if and only if
$n=1$ or $\alpha_m=0$. Hence if $\epsilon$ is any positive real 
number:
\begin{itemize}
\item[(1)]
If $n=1$, the moduli space for every $\alpha\in (\alpha_m,\infty)$ is 
isomorphic to $\mathcal{N}_{\alpha_m+\epsilon}(1,1,d_1,d_2)$ ,
\item[(2)] If $\alpha_m=0$, the moduli space 
for  every $\alpha\in (0,\infty)$ is isomorphic to 
$\mathcal{N}_{\epsilon}(n,n,d_1,d_2)$.
\end{itemize}
\end{proposition}

%%%%%%%%%%%%%%%%%%%%%%%%%%%%%%%%%%%%%%%%%%%%%%%%%%%%%%%%%%%%%%%%
\subsection{Moduli for large $\alpha$ and $\alpha\geq 2g-2$}
\label{sec:alpha0}
%%%%%%%%%%%%%%%%%%%%%%%%%%%%%%%%%%%%%%%%%%%%%%%%%%%%%%%%%%%%%%%%%

Let $\alpha>\alpha_0$. By Corollary
\ref{injective-triples}, we know that all triples in
$\mathcal{N}_\alpha(n,n,d_1,d_2)$ are of the form
\begin{equation}\label{ext-tor}
0 \lto E_2 \overset{\phi}{\lto} E_1 \lto S \lto 0,
\end{equation}
where $S$ is a torsion sheaf of degree $d=d_1-d_2$.
\begin{theorem}[Markman-Xia \cite{markman-xia:2001}]
\label{thm:irreducible-alpha0}
There is an irreducible  family $\mathcal{S}$ parameterizing
quotients $E_1\lto S\lto 0$, where $E_1$ is a rank $n$ and degree
$d_1$ locally free coherent sheaf varying on a bounded family, and
$S$ is a torsion sheaf of degree $d>0$.
\end{theorem}

\begin{theorem}\label{thm:irreducibility-alpha0}
Let $\alpha>\alpha_0$, then $\mathcal{N}_\alpha(n,n,d_1,d_2)$ is
irreducible.
\end{theorem}
\begin{proof}
Since $\alpha>\alpha_0$, an $\alpha$-semistable triple
$T=(E_1,E_2,\phi)$ defines a sequence as in (\ref{ext-tor}) and hence a
quotient $ E_1 \lto S \lto 0$ in $\mathcal{S}$. That $E_1$ varies on
a bounded family is a consequence of (1) in Proposition
\ref{prop:slope-bound-on-subbdls}. Indeed, let $E$ be a vector bundle
of degree $d$ and rank $n$ satisfying
\begin{equation}
\label{B-bound}
\mu(E')\le B
\end{equation}
for all subbundles $E'\subset E$, and fixed $B$. But then (see, for
instance, the proof of Theorem 5.6.1 in \cite{LeP}) we can find a
line bundle
$L$ of sufficiently high degree such that $H^1(E\otimes L)=0$ for all
$E$ which satisfy (\ref{B-bound}). The irreducibility of
$\mathcal{N}_\alpha(n,n,d_1,d_2)$ follows now from
the irreducibility of $\mathcal{S}$ and the fact that
  $\alpha$-semistability is an open condition.
\end{proof}

By analogy with the $n_1\neq n_2$, let us denote by 
$\mathcal{N}_L(n,n,d_1,d_2)$ the `large $\alpha$' moduli space, i.e.\ 
the moduli space of $\alpha$-semistable triples for any 
$\alpha\in (\alpha_L,\infty)$. Since $\alpha_L\geq \alpha_0$ we have 
 that all triples in
$\mathcal{N}_L(n,d_1,d_2)$ are of the form
$$
0 \lto E_2 \overset{\phi}{\lto} E_1 \lto S \lto 0\ ,
$$
and that $E_1$\ and $E_2$\ are bounded by the constraints in
Proposition \ref{prop:slope-bound-on-subbdls}.  

In the converse direction we have
\begin{proposition}
  \label{prop:bundles-stable-implies-triple-stable}
Let $T=(E_1,E_2,\phi)$ be a triple such that $\ker\phi=0$
If $E_1$ and $E_2$ are semistable, then $T$ is  $\alpha$-semistable
for large enough  $\alpha$,  i.e.
  $T \in \mathcal{N}(n,n,d_1,d_2)$.  If
  either $E_1$\ or $E_2$\ is stable, then $T$\ is
  $\alpha$-stable.  
\end{proposition}
\begin{proof} 
  Since $\ker\phi=0$, it follows (as in the
  proof of Theorem \ref{thm:stabilization}) that in any subtriple, say
  $T'=(E'_1,E'_2,\phi')$, the rank of $E'_1$\ is at least as big as the
  rank of $E'_2$, i.e.\  $n'_1\geq n'_2$.  If $n'_1>n'_2$, then 
  \eqref{eq:2.61}, \eqref{eq:2.62} and  \eqref{eq:2.63} apply, with
  $\lambda(T')<0$\ and
  $\frac{M(T')}{-\lambda(T')}\leq \alpha_L$.  Thus $\Delta_{\alpha}<0$\ 
  whenever $\alpha>\alpha_L$. For subtriples with $n'_1=n'_2$, 
  equation \eqref{eq:2.60} applies, i.e.
  $$
  \Delta_{\alpha}(T')= \mu(E'_1\oplus E'_2)-\mu(E_1\oplus
  E_2)\ 
  $$
for any $\alpha$. For such subtriples, and for any $\alpha$, it thus 
follows that 
\begin{itemize}
\item[$(1)$] $\Delta_{\alpha}(T')\leq 0$\ if both $E_1$\ and $E_2$\ are
  semistable, and
\item[$(2)$] $\Delta_{\alpha}(T') < 0$\ if at least one of the bundles is 
stable.
\end{itemize}
\end{proof}

\begin{theorem} \label{thm:existence-n1=n2}
The moduli space $\mathcal{N}^s_L(n,n,d_1,d_2)$ is non-empty.
\end{theorem}
\begin{proof}
Our strategy is to  show that  there exist rank $n$ 
 stable bundles $E_1$ and $E_2$ of degree $d_1$ and $d_2$, respectively,
and a torsion sheaf $S$ of degree $d_1-d_2$,  fitting  in an exact sequence
$$
0\lto E_2\lto  E_1\lto S \lto  0.
$$ 
The result will then  follow from Proposition \ref{prop:bundles-stable-implies-triple-stable}.

To prove this, let $E$ be a vector bundle, and let $\Quot^d(E)$ be 
the Quot scheme of quotients $E\lto S$ where $S$ is a torsion
sheaf of degree $d$. The basic fact we will need  is the  following result 
of  R. Hern\'andez (\cite{hernandez}):
Let $\psi:\mathcal{O}^n\lto S$ be an element in
 $\Quot^d(\mathcal{O}^n)$, then for a  generic  $S$, the vector bundle
$E=\ker\psi$ is stable.
Notice that $\deg E=-d$. By tensoring with a line bundle $L$ of 
big enough  degree (depending on $n$ and $d$), we can extend the result
to $\Quot^d(L^n)$, so that the kernel of any element in
$\Quot^d(L^n)$ has a given (fixed) degree.

Let $L$ be a line bundle of degree $m$ and  $d''>0$ such that
$d_1=nm-d''$.
By Hern\'andez result, if  $\psi:L^n\lto S''\in \Quot^{d''}(L^n)$
is generic, then $E_1=\ker\psi$ is a stable bundle of rank $n$ and degree
$d_1$. Let $d=d_1-d_2$ and consider a generic element 
$\eta: E_1\lto S\in \Quot^d(E_1)$. Let $E_2=\ker\eta$, and let $S'$ the
cokernel of the natural inclusion $E_2\lto L^n$. We have the
following commutative diagram:

$$\begin{array}{ccccccccc}
&&0&&0&& 0 &&\\
&&\downarrow&&\downarrow& & \downarrow&  & \\
0&\rightarrow&E_2&\rightarrow&E_1 &\rightarrow& S & \rightarrow & 0 \\
&&\parallel &&\downarrow& &\downarrow&& \\
0&\rightarrow& E_2&\rightarrow& L^n  &\rightarrow& S'  &\rightarrow& 0\\
&&\downarrow&&\downarrow& &\downarrow&& \\
& &0&\rightarrow& S''& = & S'' &\rightarrow& 0\\
&& &&\downarrow& &\downarrow&& \\
 & & & & 0  & & 0.  & &  \\
\end{array}
$$

We see from the diagram that  $E_2$ coincides with the kernel
of $L^n\lto S'$. If $S'$ is general enough we can again apply
the basic result of Hern\'andez  and conclude that $E_2$ is stable.
To show that this is indeed the case, we observe that the diagram
defines a map
$$
\Quot_0^d(E_1) \times  \Quot_0^{d''}(L^n) \lto \Quot_0^{d+d''}(L^n),   
$$
where $\Quot_0$ denotes an open non-empty subscheme of $\Quot$,  which
is surjective and finite.
\end{proof}

\begin{proposition}\label{prop:nnNL}
The moduli space $\mathcal{N}_L(n,n,d_1,d_2)$ is 
birationally  equivalent to a 
$\mathbb{P}^N$-fibration $\mathcal{P}$ over
$M^s(n,d_2)\times \Div^d(X)$, where $N=n(d_1-d_2)-1$.
\end{proposition}

\begin{proof}
Let $E_2$ be a 
rank 
$n$ and degree $d_2$ vector bundle  and let 
$S$ be  a torsion sheaf of degree $d>0$. We construct $E_1$ as an 
extension 
\begin{equation}\label{torsion-extension}
0 \lto E_2 \lto E_1 \lto S \lto 0.
\end{equation}
Such extensions are parameterized by $\Ext^1(S,E_2)$. Suppose
that   $S$ 
is of the form $S=\mathcal{O}_D$, where $D$ is a divisor
in $\Div^d(X)$. Let $L$ be a 
line bundle. Consider  the  short exact sequence 
$$
0 \lto L^*(-D) \lto L^* \lto \mathcal{O}_D \lto 0,
$$
and apply to it the functor $\Hom(\cdot,E_2)$, to  obtain the long
exact sequence
\begin{equation}
  \label{eq:long-exact-torsion}
\begin{array}{ccccccc}
  0 &\lto &H^0(E_2\otimes L) &\lto& H^0(E_2\otimes L(D)) &\lto&  \\ 
  \Ext^1(\mathcal{O}_D,E_2) & \lto & H^1(E_2\otimes L) & \lto& 
H^1(E_2\otimes L(D)) &\lto &  0.
\end{array}
\end{equation}
We thus have
$$
\dim \Ext^1(\mathcal{O}_D,E_2)=\chi(E_2\otimes L)-\chi(E_2\otimes L(D))=nd,
$$
which is positive since we are assuming $d>0$.
Taking $L$ so that $\deg(L)>>0$, we have that 
$H^1(E_2\otimes L)=0$. If $E_2$ is semistable (or more generally
it moves in a  bounded family)  we can take the  same $L$ 
for every $E_2$. Then
$$
\Ext^1(\mathcal{O}_D,E_2)=H^0(E_2\otimes L(D))/H^0(E_2\otimes L).
$$
Let $\mathcal{P}$ be the set of equivalence classes of
extensions \eqref{torsion-extension}, where $E_2$ is stable
then $\mathcal{P}$ is a  $\mathbb{P}^N$-fibration over
$M^s(n,d_2)\times \Div^d(X)$, where  
$N=nd-1=\dim \mathbb{P}(\Ext^1(\mathcal{O}_D,E_2))$. 
Setting $d=d_1-d_2$, a simple 
computation shows that
$$
\dim\mathcal{P}=
(g-1)(n_1^2 + n_2^2 - n_1 n_2) - n_1 d_2 + n_2 d_1 + 1.
$$
Cleary $\mathcal{P}$ is irreducible of the same dimension as $\mathcal{N}_L$,
and since it is contained in $\mathcal{S}$ (like $\mathcal{N}_L$) it
must birationally equivalent to $\mathcal{N}_L$.
Notice that  If $\GCD(n,d_2)=1$, then 
$\mathcal{P}$ is a Picard sheaf. 
\end{proof}

Combining
Theorems \ref{thm:irreducibility-alpha0} and \ref{thm:existence-n1=n2}, and 
Proposition \ref{prop:nnNL} we have the following.

\begin{theorem}\label{thm:moduli-n1=n2}.
The moduli space $\mathcal{N}_L(n,n,d_1,d_2)$ is non-empty and irreducible. Moreover, it is
birationally  equivalent to a 
$\mathbb{P}^N$-fibration  over
$M^s(n,d_2)\times \Div^d(X)$, where the fiber dimension is
$N=n(d_1-d_2)-1$.
\end{theorem}

\begin{theorem}\label{thm:irreducibility-moduli-stable-triples-n1=n2}
 Let $ \alpha\ge 2g-2$.  Then 
\begin{itemize}
\item[$(1)$] 
The moduli space  $\mathcal{N}^s_\alpha$ is  
 birationally equivalent 
to $\mathcal{N}_L$  and it is hence non-empty and  irreducible.

\item[$(2)$] 
If in addition either
\begin{itemize}
\item [$\bullet$] $\GCD(n,2n,d_1+d_2)=1$ and $\alpha\ge 2g-2$ is  generic, or
\item [$\bullet$] $d_1-d_2<\alpha$,
\end{itemize}
then  $\mathcal{N}_{\alpha}(n,n,d_1,d_2)$ is birationally equivalent 
to $\mathcal{N}_L(n,n,d_1,d_2)$ and hence irreducible. 
\end{itemize}
\end{theorem}
\begin{proof}

(1) {}From Theorem \ref{thm:irreducibility-alpha0} we 
know that $\mathcal{N}_L$ is birationally equivalent to 
$\mathcal{N}_L^s$. The result  follows now from Corollary 
\ref{cor:birationality} and Theorem \ref{thm:moduli-n1=n2}.

(2)  For the first  part, we observe that 
 from (4) in Proposition \ref{triples-critical-range} one has that 
$\mathcal{N}_\alpha=\mathcal{N}_\alpha^s$ if
$\GCD(n,2n,d_1+d_2)=1$ and $\alpha$ is generic, and hence the result
follows from (1). The second  part  is a consequence of  Theorem \ref{thm:irreducibility-alpha0}.
\end{proof}

%%%%%%%%%%%%%%%%%%%%%%%%%%%%%%%%%%%%%%%%%%%%%%%%%%%%%%%
\section{Existence and connectedness for $\U(p,q)$ and
$\PU(p,q)$ moduli spaces}
\label{sec:main-results}
%%%%%%%%%%%%%%%%%%%%%%%%%%%%%%%%%%%%%%%%%%%%%%%%%%%%%%%

We now return to the representation spaces $\mathcal{R}(\PU(p,q))$
and
$\mathcal{R}_{\Gamma}(\U(p,q))$, defined in section
\ref{sec:background}. Recall that we identified components of
$\mathcal{R}(\PU(p,q))$ labeled by $[a,b]\in \Z\oplus\Z/(p+q)\Z$,
and similarly identified components of
$\mathcal{R}_{\Gamma}(\U(p,q))$ labeled by $(a,b)\in \Z\oplus\Z$.

By Proposition \ref{prop:principal-jac} $\mathcal{R}_{\Gamma}(a,b)$ is
a $\U(1)^{2g}$-fibration over $\mathcal{R}[a,b]$, and hence the number
of connected components of $\mathcal{R}_{\Gamma}(a,b)$ is greater than
or equal to that of $\mathcal{R}[a,b]$. By Proposition \ref{prop:R=M}
there is an homeomorphism between $\mathcal{R}_\Gamma(a,b)$ and the
moduli space $\mathcal{M}(a,b)$ of $\U(p,q)$-Higgs bundles and this
restricts to give a homeomorphism between the subspace
$\mathcal{R}^*_\Gamma(a,b)$ of irreducible elements in
$\mathcal{R}_\Gamma(a,b)$ and the subspace $\mathcal{M}^s(a,b)$ of
stable Higgs bundles in $\mathcal{M}(a,b)$. By Proposition
\ref{prop:topology-exercise} the number of connected components of
$\mathcal{M}(a,b)$ is determined by the number of connected components
in the subspace of local minima for the Bott-Morse function defined in
Section \ref{subs:Morse}. By Theorem \ref{thm:minima=triple-moduli} we
can identify the subspace of local minima as a moduli space of
$\alpha$-stable triples, with $\alpha=2g-2$.  Summarizing, we have:
\begin{align*}
|\pi_0(\mathcal{R}[a,b])| &\leq |\pi_0(\mathcal{R}_\Gamma(a,b))|\quad
&\mathrm{[Proposition\ \ref{prop:principal-jac}]}\\
 &\leq |\pi_0(\mathcal{M}(a,b))|\quad
&\mathrm{[Proposition\ \ref{prop:R=M}]}\\
 &\leq |\pi_0(\mathcal{N}(a,b))|\quad
&\mathrm{[Proposition\ \ref{prop:topology-exercise}]}\\
&=|\pi_0(\mathcal{N}_{2g-2}(n_1,n_2,d_1,d_2))|\quad
&\mathrm{[Proposition\ \ref{thm:minima=triple-moduli}]}
\end{align*}
where $|\pi_0(\cdot)|$ denotes the number of components, and (in the
notation of Section
\ref{sec:stable-triples}) the moduli space of triples which appears
in the last line is either
$\mathcal{N}_{2g-2}(p,q,a + p(2g-2),b)$ (if $a/p \leq b/q$) or
$\mathcal{N}_{2g-2}(q,p,b + q(2g-2),a)$ (if $a/p \geq b/q$).
Similarly, replacing Proposition \ref{prop:topology-exercise} with
Proposition \ref{prop:N-M-list}, we get that

$$ |\pi_0(\bar{\mathcal{R}}^*[a,b])|\le
|\pi_0(\bar{\mathcal{N}}^s_{2g-2}(n_1,n_2,d_1,d_2))|\ .
$$
In particular, if the moduli spaces of triples are connected, then so
are the spaces $\mathcal{R}[a,b]$ and
$\bar{\mathcal{R}}^*[a,b]$.

%Thus we see that the number of connected components in
%$\mathcal{R}(\PU(p,q))$ and $\mathcal{R}_{\Gamma}(\U(p,q))$ is
%completely determined by the number of connected components in
%$\mathcal{N}_{2g-2}(p,q,a + p(2g-2),b)$ or
%$\mathcal{N}_{2g-2}(q,p,b + q(2g-2),a)$.

We use the results of Sections
\ref{crossing-critical-values}-\ref{sec:n_1=n_2} to determine the
connectedness
of the spaces $\mathcal{N}_{2g-2}(p,q,a + p(2g-2),b)$ and
$\mathcal{N}_{2g-2}(q,p,b + q(2g-2),a)$.  Our main tool is Corollary
\ref{cor:birationality}, which allows us convert the problem into one
of counting the components for $\mathcal {N}_L$, the moduli space of
$\alpha$-stable triples for sufficiently large
$\alpha$. In the case $p\ne q$, Theorem \ref{thm:largealpha} and its dual
supply the requisite details. The results in Section
\ref{sec:n_1=n_2} cover the the case
$p=q$, which is not as well understood as the case $p\ne q$.

Before presenting specific results, we make some general
observations.

%%%%%%%%%%%%%%%%%%%%%%%%%%%%%%%%%%%%%%%%%%%%%%%%%%%%%%%
\subsection{Toledo invariant and the number of
components}\label{subs:ab-tau}
%%%%%%%%%%%%%%%%%%%%%%%%%%%%%%%%%%%%%%%%%%%%%%%%%%%%%%

As described in Section \ref{sec:topological-bounds}, it is useful to 
introduce the combination

$$\tau :=\tau(a,b)=\frac{2}{n}(qa-pb) \,$$
known  as the Toledo invariant. Throughout this section, for 
convenience, we set $n=p+q$.

The map $(a,b)\mapsto \tau(a,b)$ defines an invariant on each 
component 
$\mathcal{R}_\Gamma(a,b)$. Moreover, since the map factors through
$\Z\oplus \Z/(p+q)\Z$, it defines an invariant of 
$\mathcal{R}[a,b]$, where it takes the same value as on
$\mathcal{R}_\Gamma(a,b)$.

\begin{proposition}\label{prop:[a,b]vsTau} Suppose that $\GCD(p,q)=k$.
Then the map 
\begin{align*}
\tau:&\Z\oplus\Z/(p,q)\Z\longrightarrow \frac{2}{n}\Z\ \\
&[a,b]\mapsto \frac{2}{n}(aq-bp) 
\end{align*}
fits in an exact sequence
\begin{equation}
\begin{CD}
0@>>> &\Z/k\Z @>\sigma >> 
\Z\oplus\Z/(p,q)\Z@>\tau >>
\frac{2k}{n}\Z@>>> 0
\end{CD}
\end{equation}
where the map $\sigma$ is $[t]\mapsto [t\frac{p}{k},t\frac{q}{k}]$. 
In particular, $\tau$ is a $k:1$ map onto the subset 
$\frac{2k}{n}\Z\subset \frac{2}{n}\Z$. 
\end{proposition}

\begin{proof} The map $\sigma$ is clearly injective, and 
$\tau\circ\sigma=0$. To see that $\ker(\tau)=\im(\sigma)$, observe 
that if $\tau[a,b]=0$ then either $a=b=0$ or 
$\frac{a}{b}=\frac{p}{q}$, i.e.\ $[a,b]=[t\frac{p}{k},t\frac{q}{k}]$
for some $t\in\Z$.  Finally, if 
$a_0q-b_0p=k$, then for any $l\in\Z$ we have $\tau[la_0,lb_0]=\frac{2k}{n}l$.
Thus $\tau$ is surjective onto $\frac{2k}{n}\Z$.
\end{proof}

\begin{remark} Proposition \ref{prop:[a,b]vsTau} shows why we 
must\footnotemark use 
$[a,b]$ rather than $\tau$ to label the components of $\mathcal 
{R}(\PU(p,q))$ or of $\mathcal {R}_{\Gamma}(\U(p,q))$. 
\footnotetext{Unless $p$ and $q$ are
coprime, in which case the correspondence between $[a,b]$ and $\tau$ 
is 1-1.} There are, nevertheless, important features of 
$\mathcal{R}[a,b]$ and $\mathcal{R}_{\Gamma}(a,b)$ 
which depend on $\tau$ rather than on $[a,b]$ or $(a,b)$. The 
rigidity result described in Section \ref{subs:rigidity} is one such 
feature.  Others emerge in the analysis of the minimal submanifolds, 
i.e.\ in the moduli spaces 
$\mathcal{N}_{2g-2}(p,q,a + p(2g-2),b)$ and 
$\mathcal{N}_{2g-2}(q,p,b + q(2g-2),a)$. In particular, for triples of 
the indicated type, the key critical values of the stability 
parameter 
$\alpha$  (i.e.\ $\alpha_m, \alpha_0, \alpha_t,\alpha_M$, in the 
notation of Section \ref{sec:stable-triples}-\ref{sec:n_1=n_2}), as 
well as the all-important location of $2g-2$ within the interval 
$(\alpha_m,\alpha_M)$, are all determined by $\tau$.
\end{remark}

Proposition \ref{prop:[a,b]vsTau} allows us to count the number of 
components in the decomposition 
$\mathcal{R}(\PU(p,q))=\bigcup_{(p,q)}\mathcal {R}[a,b]$.
Recall (from Section 
\ref{sec:topological-bounds}) that if 
$p\le q$ then the Toledo invariant is bounded by $|\tau|\le 
\tau_M=2p(g-1)$. Thus the set of allowed values for $\tau$ is 
contained in  $[-\tau_M,\tau_M]\cap \frac{2}{n}\Z$. 

\begin{definition}\label{defn: C-N} Suppose that 
$\GCD(p,q)=k$. Define
\begin{equation}\label{eqn:C}
\mathcal{C}=\tau^{-1}([-\tau_M,\tau_M]\cap \frac{2k}{n}\Z)\ ,
\end{equation}
\noindent where $\tau$ is the map defined in Proposition \ref{prop:[a,b]vsTau}. 
\end{definition}

The following is then an immediate corollary of Proposition 
\ref{prop:[a,b]vsTau}.
 
\begin{corollary}\label{cor:number-of-[a,b]}
 Suppose that $\GCD(p,q)=k$ and 
$\mathcal{C}$ is as above. Then 
$\mathcal{C}$ is precisely the 
set of all the points in $\Z\oplus\Z/(p,q)\Z$ which label components 
$\mathcal{R}[a,b]$ in $\mathcal{R}(\PU(p,q))$.  The
cardinality of $\mathcal{C}$ is
\begin{align*}
|\mathcal {C}|&=2n\min\{p,q\}(g-1)+k\\ 
&=|([-\tau_M,\tau_M]\cap 
\frac{2}{n}\Z)|+\GCD(p,q)-1\ .
\end{align*} 
\end{corollary}

\begin{proof} The first statement is a direct consequence of 
Proposition 
\ref{prop:[a,b]vsTau} and the bound on $\tau$. Suppose for definiteness that
$\min\{p,q\}=p$. Then since 
$\tau_M=2\min\{p,q\}(g-1)=\frac{2k}{n}(n\frac{p}{k}(g-1))\in \frac{2k}{n}\Z$, 
the number of points in $[-\tau_M,\tau_M]\cap \frac{2k}{n}\Z$ is  
$2n\frac{p}{k}(g-1)+1$. The second statement now follows from the 
fact that $\tau$ is a $k:1$ map. The proof is similar if 
$\min\{p,q\}=q$.
\end{proof}

%%%%%%%%%%%%%%%%%%%%%%%%%%%%%%%%%%%%%%%%%%%%%%%%%%%%%%%%%%%%%%%%%%%%%%%%%%%
\subsection{Genericity of 2g-2 and coprimality conditions}\label{subs: 2g-2}
%%%%%%%%%%%%%%%%%%%%%%%%%%%%%%%%%%%%%%%%%%%%%%%%%%%%%%%%%%%%%%%%%%%%%%%%%%%%

In general the Higgs moduli spaces 
$\mathcal{M}(a,b)$ and also the triples moduli spaces 
$\mathcal{N}_{2g-2}(p,q,a + p(2g-2),b)$ and 
$\mathcal{N}_{2g-2}(q,p,b + q(2g-2),a)$, are not smooth. In 
both the Higgs and the triples moduli spaces, the singularities occur 
at points representing strictly semistable objects.  Apart from 
complications caused by such singularities, further difficulties can 
arise if 
$2g-2$ is not a generic value in the range of the stability parameter
for the triples.

We now show that both types of problems are avoided if we make the 
assumption that $p+q$ and $a+b$ are coprime, i.e.\ that 
$\GCD(p+q,a+b)=1$. 

\begin{proposition}
Suppose that $\GCD(p+q,a+b)=1$. Then 
\begin{enumerate}
\item $\mathcal{M}(a,b)$ is smooth,
\item $\alpha=2g-2$ is not a critical value for triples of
type $(p,q,a + p(2g-2),b)$ or $(q,p,b + q(2g-2),a)$, and
\item The moduli spaces $\mathcal{N}_{2g-2}(p,q,a + p(2g-2),b)$ and 
$\mathcal{N}_{2g-2}(q,p,b + q(2g-2),a)$ are smooth.
\end{enumerate}
\end{proposition}
\begin{proof} The first statement is simply a restatement of Remark 
\ref{remark:coprime}. The second and third statements follow from
Lemma~\ref{critical-integer}. Indeed, for triples of type 
$(p,q,a + p(2g-2),b)$, the coprime condition in that lemma becomes
\begin{equation}\label{not-critical}
\GCD(p+q,a+b +(2g-2)(p+q))=1\ .
\end{equation}
But (\ref{not-critical}) is equivalent to 
$\GCD(p+q,a+b)=1$. The proof is similar for triples of type 
$(q,p,b + q(2g-2),a)$.
\end{proof}

Recall that $p$ and $q$ are fixed, but $[a,b]$ runs over 
$\mathcal {C}$, the indexing set for the components
$\mathcal{R}[a,b]$. The coprime condition 
$\GCD(p+q,a+b)=1$ can thus be satisfied on some components but not on others.

\begin{proposition}\label{prop:good-bad}Fix $p$ and $q$ and let 
$\mathcal {C}\subset \Z\oplus\Z/(p+q)\Z$ be as in 
Definition \ref{defn: C-N}. Let 
$\mathcal {C}_{(1)}$ denote the subset of classes 
$[a,b]\in\mathcal {C}$ for which the condition $\GCD(p+q,a+b)=1$
is satisfied. Then both $\mathcal {C}_{(1)}$ and its complement in 
$\mathcal {C}$ are non-empty. 
\end{proposition} 

\begin{proof}
If $a=p$ and $b=q-1$ then $\GCD(p+q,a+b)=1$. Also, 
$\tau(p,q-1)=\frac{2p}{p+q}$, which is in 
$[-\tau_M,\tau_M]\cap \frac{2k}{n}\Z$. Thus
$[p,q-1]$ is in $\mathcal {C}_{(1)}$. It is similarly straightforward
to see that $(p,q)= (0,0)$ defines an element in $\mathcal {C}- 
\mathcal {C}_{(1)}$, as does $(p,q)=(p,-p)$ if $p\leq q$
or $(p,q)= (q,-q)$ if $q\leq p$.
\end{proof}
 
It seems somewhat complicated to go beyond this result and completely 
enumerate the elements in $\mathcal {C}_{(1)}$. The following result 
is, however, a step in that direction. 

\begin{definition} Let $\Omega\subset\R\oplus\R$ be the region  bounded by 
\begin{itemize}
\item the ray $ b=q$  and $a\le p$,
\item the ray $ a=p$ and $b\le q$,
\item the ray $ a=0$  and $b\le 0$,
\item the ray $ b=0$ and $a\le 0$,
\item the line $ aq-bp=np(g-1)$, and
\item the line $ aq-bp=-np(g-1)$,
\end{itemize} 
and including all the boundary lines except the first two rays. Let
$\Omega_{\Z}$ be the set of integer points in $\Omega$, i.e.\ 
$\Omega_{\Z}=\Omega\bigcap\Z\oplus\Z$. We refer to $\Omega_{\Z}$  as the 
{\bf fundamental region for $(p,q)$}. 
\end{definition}

\begin{proposition} Suppose that $p$ and $q$ are integers with
$\GCD (p,q)=k$ and $p\le q$\footnotemark.
\footnotetext{With an analogous Proposition for the case $p\ge q$.}
\begin{enumerate}
\item There is a bijection between $\mathcal{C}$ and 
$\Omega_{\Z}$. 
\item If $(a,b)$ lies in $\Omega_{\Z}$ then $d=a+b$ satisfies the 
bounds 
\begin{equation}
-n(g-1)\le d< n\ . 
\end{equation}
All values of $d$ in this range occur. 
\item Let $l_t$ denote the line $aq-bp=tk$. Then the points
on $l_t\bigcap\Omega_{\Z}$ define the locus of points $(a,b)$ for 
which $\tau(a,b)=t\frac{2k}{n}$.
\item The line $l_t$ intersects $\Omega_{\Z}$ for $-\frac{np}{k}(g-1)\le t
\le \frac{np}{k}(g-1)$ For each integer $t$ in this range, there are 
$k$ points on $l_t\bigcap\Omega_{\Z}$.
\item For a fixed $t$, all integer points $(a,b)\in l_t\bigcap\Omega_{\Z}$ have the same 
$\GCD(d,\frac{n}{k})$, where $d=a+b$.
\item If $\GCD(d,\frac{n}{k})\ne 1$ then $\GCD(d',n)\ne 1$ for all 
$(a',b')\in l_t\bigcap\Omega_{\Z}$.
%\item The condition for d to occur is that there is an allowed $\tau$ and 
%an allowed $b$ such that
%$pb$ is congruent to $-n\tau/2$ mod $q$ YIKES!
\end{enumerate}
\end{proposition}

\begin{proof}
(1) Suppose first that $\frac{a}{p}\le \frac{b}{p}$. Pick $l$ such 
that $0\le a+lp\le p$. Then $b+lq\le q$, so that $(a+lp, b+lq)$ is in 
the fundamental region. Similarly, if  $\frac{a}{p}\ge \frac{b}{p}$ 
then we pick $l$ such that $0\le b+lq\le q$ and see that $a+lp\le p$. 
In this way we get a well defined map from $\mathcal {C}$ to the 
fundamental region. The map is clearly injective. To see that it is 
surjective, notice that the boundary lines $ aq-bp=np(g-1)$, and 
$ aq-bp=-np(g-1)$ correspond to the conditions $\tau=\tau_M$ and 
$\tau=-\tau_M$ respectively.

(2) This is clear from a sketch of the fundamental region. In such a 
sketch, the loci of points with constant value of $d=a+b$ are 
straight lines of slope $-1$. Since $p\le q$, the extreme cases are 
those of the lines passing through the points $(0,-n(g-1)$ and 
$(p,q)$. Using the points $(0,b)$ with $-n(g-1)\le b<q$ we get points at 
which all values of $d$ in the range $-n(g-1)\le d<q$ are realized. 
We get values of $d$ in the range $g\le d<n$ at the points $(a, 
q-1)$, with $1\le a<p$. 

(3)-(4) This is simply a re-statement of Proposition 
\ref{prop:[a,b]vsTau}.

(5)-(6) Both follow from the fact that for any two points $(a,b)$ and 
$(a',b')$ on $l_t$, we get $d'=d+s\frac{n}{k}$ for some $s\in\Z$.
\end{proof}

\begin{remark} There is no converse to (6), i.e.\ it is possible to 
have $\GCD(d',n)\ne 1$ for some $(a',b')\in l_t\bigcap\Omega_{\Z}$, 
even if $\GCD(d,\frac{n}{k})= 1$. For example, take $p=2, q=4, a=-1, 
b=0, a'=0, b'=2 $, and $t=-2$. Then $\GCD(d',n)=2$ but 
$\GCD(d,\frac{n}{k})= 1$. 
\end{remark}

%%%%%%%%%%%%%%%%%%%%%%%%%%%%%%%%%%%%%%%%%%%%%%%%%%%%%%%%%%%%%
\subsection{Moduli spaces of $\U(p,q)$-Higgs bundles}
\label{sec:results-U(p,q)-Higgs}
%%%%%%%%%%%%%%%%%%%%%%%%%%%%%%%%%%%%%%%%%%%%%%%%%%%%%%%%%%%%%

In this section we give our main results on the moduli spaces of
$\U(p,q)$-Higgs bundles.  As noted several times before, the results
are conveniently stated using the Toledo invariant $\tau(a,b) =
2(qa-pb)/(p+q)$.  For completeness we recall from
Remark~\ref{rem:tau-M} that the maximal value of the Toledo invariant
is $\tau_M = \min\{p,q\}(2g-2)$.  In the case $p=q$ these definitions
simplify to $\tau = a-b$ and $\tau_M = p(2g-2)$.

Recall from Proposition~\ref{prop:smoothness-higgs} that, whenever the
moduli space $\mathcal{M}^s(a,b)$ of stable $\U(p,q)$-Higgs bundles
with invariants $(a,b)$ is non-empty, it is a smooth complex manifold
of dimension $1 + (p+q)^2(g-1)$.  We shall refer to this dimension as
the \emph{expected dimension} in the following.

\begin{theorem}\label{thm:tau=0-M(a,b)}
  Let $(a,b)$ be such that $\tau(a,b)=0$, then $\mathcal{M}(a,b)$ is
  non-empty and connected.
\end{theorem}
\begin{proof} It follows from (1) in Proposition \ref{prop:MW}, applied
to the triples in $\mathcal{N}(a,b)$, that if $\tau=0$ then 
$\frac{a}{p}=\frac{b}{q}$.  The result is thus simply a re-statement
of (3) in Theorems
\ref{thm:connected-triples-connected-higgs} and \ref{thm:minima=triple-moduli} 
(plus the non-emptiness of the moduli spaces $M(p,a)$ and $M(q,b)$).
\end{proof}

\begin{remark}
It is not clear if 
$\mathcal{M}^s(a,b)$ is non-empty.  Notice that the minimal
subvariety is $\mathcal{N}(a,b)=M(p,a)\times M(q,b)$, in which every 
element is reducible. To show that strictly stable points exist we 
must look beyond the minima.  We observe though, that if it is 
non-empty, then  $\mathcal{M}^s(a,b)$ is a smooth  connected 
 manifold of the expected dimension. 
\end{remark}

\begin{theorem}\label{thm:Ms(a,b)-general-tau}
  Let $(a,b)$ be such that $0<|\tau(a,b)|<\tau_M$. Then the
  closure $\bar{\mathcal{M}}^s(a,b)$ of the moduli space of stable
  $\U(p,q)$-Higgs bundles with fixed invariant $(a,b)$ is connected.
  Moreover, the stable locus $\mathcal{M}^s(a,b)$
  is a smooth non-empty manifold of the expected dimension.
\end{theorem}
\begin{proof}
  If  $p \neq q$,  from  (2) of
  Proposition~\ref{prop:MW} we have  that $2g-2 < \alpha_M$ and therefore
  Theorem~\ref{thm:irreducibility-moduli-stable-triples} shows that
  $\mathcal{N}^{s}_{2g-2}$ is non-empty and irreducible. If $p=q$,
Theorem \ref{thm:irreducibility-moduli-stable-triples-n1=n2} shows that
  $\mathcal{N}^{s}_{2g-2}$ is non-empty and irreducible. 
 From this we
  draw two conclusions. Firstly, we see from (1) and (2) of
  Theorem~\ref{thm:connected-triples-connected-higgs} that
  $\bar{\mathcal{M}}^s(a,b)$ is connected.  Secondly, it  follows that
  $\mathcal{M}^s(a,b)$ is a smooth non-empty manifold of the expected
  dimension.
\end{proof}

\begin{theorem}\label{connectednes-p=q-maximal-toledo}
  Let $p=q$ and let $(a,b)$ be such that $|\tau(a,b)|=\tau_M$. Then
  $\mathcal{M}(a,b)$ is non-empty and connected. Moreover,
  $\mathcal{M}^s(a,b)$ is non-empty, and smooth of the expected
  dimension.
\end{theorem}

\begin{proof}
  Since $p=q$, the Toledo invariant is $\tau=a-b$. Suppose for
  definiteness that $a>b$.  Since, by hypothesis,
  $a-b=\tau_M=p(2g-2)$, the moduli space of triples
  $\mathcal{N}_{2g-2}(p,p,b + p(2g-2),a)$ is simply
  $\mathcal{N}_{2g-2}(p,p,a,a)$, which by
  Corollary~\ref{cor:largealphamoduli} is irreducible.  Now, by (2) in
  Theorem~\ref{thm:connected-triples-connected-higgs}, we have that
  $\mathcal{M}(a,b)$ is connected.  A similar argument applies if
  $a<b$.
\end{proof}

If $p\ne q$ then Theorem \ref{prop:rigidity} gives a description of
$\mathcal{M}(a,b)$ when $(a,b)$ is such that $|\tau(a,b)|=\tau_M$.
For definiteness, assume that $p<q$ and $\tau(a,b)=\tau_M$. Then, using the
more precise notation $\mathcal{M}(p,q,a,b)$ for the moduli space of
semistable $\U(p,q)$-Higgs bundles with topological invariant $(a,b)$,
Proposition~\ref{prop:rigidity} gives
\begin{itemize}
\item $
    \mathcal{M}(p,q,a,b)\cong \mathcal{M}(p,p,a,a-p(2g-2))
\times M(q-p, b-a +p(2g-2))$,

where  $M(q-p, b-a +p(2g-2))$ is the moduli space of polystable 
bundles of degree $q-p$ and rank $b-a +p(2g-2)$. 
\item The dimension at a 
smooth point in 
$\mathcal{M}(p,q,a,b)$ is $2+ (p^2+5q^2-2pq)(g-1)$, and it is
hence strictly smaller than the expected dimension. 
\item Every element in
$\mathcal{M}(p,q,a,b)$ is strictly semistable, i.e.\ 
$\mathcal{M}^s(a,b)$ is empty.
\end{itemize}
Similar statements  hold for $p>q$ and/or $\tau=-\tau_M$. We can now 
add: 

\begin{theorem}\label{thm:tauM-M(a,b)}
  Let $p\neq q$ and let $(a,b)$ be such that $|\tau(a,b)|=\tau_M$.
  Then $\mathcal{M}(a,b)$ is non-empty and connected.
 \end{theorem}

\begin{proof}
This follows from (1) and (2) of
Theorem~\ref{thm:connected-triples-connected-higgs} and 
Corollary~\ref{cor:N-alpha-M}. 
\end{proof}

\begin{theorem}\label{thm:gcd=1-M(a,b)}
 Let $(a,b)$ be such that $|\tau(a,b)| \leq\tau_{M}$, and suppose
  moreover that $\GCD(p+q,a+b)=1$. Then $\mathcal{M}(a,b)$ 
  is a non-empty, connected,
  smooth manifold of the expected dimension.
\end{theorem}
\begin{proof}
Under the condition $\GCD(p+q,a+b)=1$ there are no strictly
semi-stable $\U(p,q)$-Higgs bundles and hence $\mathcal{M}(a,b)^s =
\mathcal{M}(a,b)$ (cf.\ Remark~\ref{remark:coprime}).  

Note also that the condition $\GCD(p+q,a+b)=1$ excludes the
possibility $\tau = 0$ and, when $p \neq q$, the possibility $|\tau| =
\tau_M$ (by Proposition~\ref{prop:rigidity}).  Hence we have, in fact,
that $0 < |\tau| \leq \tau_M$ and, when $p \neq q$, that $0 < |\tau| <
\tau_M$. 

Thus the statement is immediate from the preceding Theorems.
\end{proof}

In the case $p=q$ and $(p-1)(2g-2)< |\tau| \leq \tau_M$ we do not 
need to impose any coprimality conditions in order to prove 
connectedness of the moduli spaces $\mathcal{M}(a,b)$.  This follows 
from the fact that in this situation no flips are required to go from 
$\mathcal{N}_{L}$ to $\mathcal{N}_{2g-2}$.  Thus we have the 
following. 

\begin{theorem}\label{thm:MX-result}
  Assume that $p=q$ and that $\tau = a - b$ is such that 
  \begin{displaymath}
    (p-1)(2g-2) < |\tau| \leq \tau_M = p(2g-2),
  \end{displaymath}
  then $\mathcal{M}(a,b)$ is non-empty and  connected.
\end{theorem}

\begin{proof}
For definiteness let us assume that $a>b$ (an analogous argument 
applies to $a<b$). We must then prove the connectedness of the  
moduli space of ($2g-2$-polystable)  triples 
$\mathcal{N}_{2g-2}(p,p,a + p(2g-2),b)$
Recall from 
\eqref{alpha0-equal-ranks} that
\begin{displaymath}
  \alpha_{0} = d_{1} - d_{2} = a-b + p(2g-2).
\end{displaymath}
Since, by hypothesis,  $a-b> (p-1)(2g-2)$ we have that
$2g-2>\alpha_0$. We thus have that $\mathcal{N}_{2g-2}(p,p,a + p(2g-2),b)$
is non-empty and irreducible by Theorem 
\ref{thm:irreducibility-moduli-stable-triples-n1=n2}, 
and hence $\mathcal{M}(a,b)$
is connected by Theorem \ref{thm:connected-triples-connected-higgs}.
\end{proof}
\begin{remark}
The connectedness part of Theorem \ref{thm:MX-result} is a 
result previously proved by  Markman and Xia
\cite{markman-xia:2001}.
\end{remark}

%%%%%%%%%%%%%%%%%%%%%%%%%%%%%%%%%%%%%%%%%%%%%%%%%%%%%%%%%%%%%%%%%%%%%%
\subsection{Moduli spaces of $\U(p,q)$ and $\PU(p,q)$ representations}
%%%%%%%%%%%%%%%%%%%%%%%%%%%%%%%%%%%%%%%%%%%%%%%%%%%%%%%%%%%%%%%%%%%%%%

In this section we finally give our results on the moduli spaces of
representations. We begin by translating the results of 
Section~\ref{sec:results-U(p,q)-Higgs} into results about the moduli 
spaces of $\U(p,q)$ representations via Proposition~\ref{prop:R=M}. 
We denote by $\bar{\mathcal{R}}_\Gamma^*(a,b)$ the closure of 
$\mathcal{R}_\Gamma^*(a,b)$ in $\mathcal{R}_\Gamma(a,b)$.

\begin{theorem}\label{thm:results-R-Gamma}
  
Let $(a,b)$ be such that $|\tau(a,b)|\le\tau_M$. Unless further 
restrictions are imposed, let $(p,q)$ be any pair of positive 
integers. 

\begin{itemize}
  
\item[$(1)$] If either of the following sets of conditions apply, 
then the moduli space $\mathcal{R}_\Gamma^*(a,b)$ of irreducible 
semi-simple representations, is a non-empty, smooth manifold of the 
expected dimension, with connected closure 
$\bar{\mathcal{R}}_\Gamma^*(a,b)$:
\begin{enumerate}
\item[$(i)$] $0<|\tau(a,b)|<\tau_M$ ,
\item[$(ii)$] $|\tau(a,b)|=\tau_M$ and $p=q$
\end{enumerate}

\item[$(2)$] If any one of the 
following sets of conditions apply, then the moduli space 
$\mathcal{R}_\Gamma(a,b)$ of all semi-simple representations is 
non-empty and connected: 
\begin{enumerate}
\item[$(i)$] $\tau(a,b)=0$,
  
\item[$(ii)$] $|\tau(a,b)|=\tau_M$ and $p\neq q$. , 
  
\item[$(iii)$] $(p-1)(2g-2) < |\tau| \leq \tau_M = p(2g-2)$ and 
$p=q$,

\item[$(iv)$]  $\GCD(p+q,a+b)=1$ 
\end{enumerate}

\item[$(3)$] If  $|\tau(a,b)|=\tau_M$ and $p\neq q$ 
then any representation  in $\mathcal{R}_\Gamma(a,b)$ is reducible  
(i.e. $\mathcal{R}^*_\Gamma(a,b)$ is empty).  If $p < q$, then any 
such representation decomposes as a direct sum of a semisimple 
representation of 
$\Gamma$ in $\U(p,p)$ with maximal Toledo invariant and a semisimple 
representation in  $\U(q-p)$.  Thus, if $\tau = p(2g-2)$ then there 
is an isomorphism 
  \begin{displaymath}
    \mathcal{R}_\Gamma(p,q,a,b) \cong \mathcal{R}_\Gamma(p,p,a,a - p(2g-2))
    \times R_{\Gamma}(q-p,b-a + p(2g-2)),
  \end{displaymath}
  where the notation $\mathcal{R}_\Gamma(p,q,a,b)$ indicates
  the moduli space of representations of $\Gamma$ in $\U(p,q)$ with
  invariants $(a,b)$, and $R_{\Gamma}(n,d)$ denotes the moduli space of
  degree $d$ representations of $\Gamma$ in $\U(n)$. (A similar result
  holds if $p>q$, as well as if $\tau=-p(2g-2)$).
\item[$(4)$] If  $\GCD(p+q,a+b)=1$ then $\mathcal{R}_\Gamma(a,b)$ 
is a smooth manifold of the expected dimension. 
\end{itemize}
\end{theorem}

{}From Theorem~\ref{thm:results-R-Gamma} and 
Proposition~\ref{prop:principal-jac} we finally obtain the following 
theorem about the moduli spaces of $\PU(p,q)$ representations.  Note 
that the closure $\bar{\mathcal{R}}^*[a,b]$ in $\mathcal{R}[a,b]$ is 
the image of $\bar{\mathcal{R}}_\Gamma^*(a,b)$ under the map of 
Proposition~\ref{prop:principal-jac}, hence these two spaces have the 
same number of connected components. 

\begin{theorem}
\label{thm:results-R}
Let $(a,b)$ be such that $|\tau(a,b)|\le\tau_M$. Unless further 
restrictions are imposed, let $(p,q)$ be any pair of positive 
integers. 

\begin{itemize}
  
\item[$(1)$] If 
either of the following sets of conditions apply, then the moduli 
space $\mathcal{R}^*[a,b]$ of irreducible semi-simple 
representations, is non-empty, with connected closure 
$\bar{\mathcal{R}}^*[a,b]$: 
\begin{enumerate}
\item[$(i)$] $0<|\tau(a,b)|<\tau_M$ ,
\item[$(ii)$] $|\tau(a,b)|=\tau_M$ and $p=q$
\end{enumerate}

\item[$(2)$] If any one of the 
following sets of conditions apply, then the moduli space 
$\mathcal{R}[a,b]$ of all semi-simple representations is non-empty 
and connected: 
\begin{enumerate}
\item[$(i)$] $\tau(a,b)=0$,
  
\item[$(ii)$] $|\tau(a,b)|=\tau_M$ and $p\neq q$. , 
  
\item[$(iii)$] $(p-1)(2g-2) < |\tau| \leq \tau_M = p(2g-2)$ and 
$p=q$,

\item[$(iv)$]  $\GCD(p+q,a+b)=1$ 
\end{enumerate}

\item[$(3)$] If  $|\tau(a,b)|=\tau_M$ and $p\neq q$ 
then any representation  in $\mathcal{R}[a,b]$ is reducible  (i.e.\ 
$\mathcal{R}^*[a,b]$ is empty).  If $p < q$, then any such representation
reduces to a semisimple  representation of 
$\pi_1 X$  in  
  $\mathrm{P}(\U(p,p) \times \U(q-p))$, such that the representation in
  $\PU(p,p)$ induced via projection on the first factor has maximal
  Toledo invariant.  (A similar result
  holds if $p>q$, as well as if $\tau=-p(2g-2)$).
\end{itemize}
\hfill\qed
\end{theorem}

\begin{remark}
  As explained by Hitchin in \cite[Section 5]{hitchin:1987}, the
  moduli space of irreducible representations in the adjoint form of a
  Lie group is liable to acquire singularities, because of the
  existence of stable vector bundles which are fixed under the action
  of tensoring by a finite order linebundle.  For this reason we are
  not making any smoothness statements in Theorem~\ref{thm:results-R}.
\end{remark}

%%%%%%%%%%%%%%%%%%%%%%
%\section{References}
%%%%%%%%%%%%%%%%%%%%%%

%\bibliography{upq}
%\bibliographystyle{invmatplain}

\end{document}